\documentclass[11pt ,reqno]{amsart}
\usepackage{geometry}                
\usepackage{graphicx}
\usepackage{amssymb}
\usepackage{epstopdf}
\usepackage{xcolor}
\usepackage{amsthm}
\usepackage{amsmath}
\usepackage{hyperref}
\usepackage{todonotes}
\usepackage{comment}
\usepackage{tikz}
\newenvironment{proof1}[1]{\begin{trivlist} \item[] {\em Proof of #1:}}{\newline \textcolor{white}{.}\hfill $\Box$
                      \end{trivlist}}

\theoremstyle{plain}

\newtheorem{thm}{Theorem}
\newtheorem{lemma}{Lemma}[section]
\newtheorem{prop}{Proposition}[section]
\newtheorem{cor}{Corollary}

\newtheorem{conjecture}{Conjecture}[section]

\newtheorem{remark}{Remark}
\newtheorem{defn}{Definition}[section]

\newcommand{\pa}{\partial}

\newcommand{\R}{\mathbb{R}}

\newcommand{\w}{\bm{w}}

\newcommand{\norm}[1]{\left\lVert#1\right\rVert}

\newcommand{\sub}[1]{_{_{#1}}}

\author{Thomas Beck}
\email{tbeck7@fordham.edu}
\address{Mathematics Department, Fordham University, Bronx, NY 10458}
\author{Yaiza Canzani}
\email{canzani@email.unc.edu}
\address{Department of Mathematics, University of North Carolina, Chapel Hill, NC 27514}
\author{Jeremy L. Marzuola}
\email{marzuola@math.unc.edu}
\address{Department of Mathematics, University of North Carolina, Chapel Hill, NC 27514}

\vspace{-0.3in}
\title{Eigenfunction asymptotics and nodal domain estimates for dumbbell domains}
\date{\vspace{-0.4in}\today}

\begin{document}

\begin{abstract}
We study the nodal structure of Neumann eigenfunctions on planar dumbbell domains as the width of the neck joining two fixed end domains tends to zero. Under simplicity and nondegeneracy assumptions, we show that sufficiently thin necks force an increasing number of nodal domains as the limiting eigenvalue grows. At the bottom of the spectrum this forcing can attain Courant’s upper bound: in particular, the third Neumann eigenfunction has exactly three nodal domains and is Courant sharp. At the same time, nodal deficiency of an eigenfunction on an end domain can persist on the full dumbbell, giving eigenfunctions that are not Courant sharp. The main analytic difficulty is that the singular limit may vanish on one or both end domains even though these regions contribute nodal domains for every positive neck width. We overcome this by deriving refined asymptotics for the first nonzero profiles on the vanishing regions, expressed through Neumann Green's functions with poles at the neck attachment points and obtained by a two-dimensional matched asymptotic analysis with logarithmic terms.
\end{abstract}

\maketitle

\section{Introduction}

\begin{figure}[b]
    \centering
\begin{tikzpicture}[scale=1.8, thick]

\def\a{0.45}
\def\b{0.65}  
\def\w{0.14}

\fill[gray!4] (-\a,-1) rectangle (\a,1);

\fill[gray!30] (-\w,-1) rectangle (\w,1);
\draw (-\w,-1) rectangle (\w,1);
\node at (0.5,0) {$Q(w)$};

\fill[blue!10]
(-\b,1)
-- (\b,1)
-- (1.15,1.35)
.. controls (1.45,1.75) and (1.05,2.2) .. (0.35,2.35)
-- (-0.35,2.25)
.. controls (-1.05,2.1) and (-1.35,1.55) .. (-1.05,1.2)
-- (-\b,1)
-- cycle;

\draw[blue!70!black]
(-\b,1)
-- (\b,1)
-- (1.15,1.35)
.. controls (1.45,1.75) and (1.05,2.2) .. (0.35,2.35)
-- (-0.35,2.25)
.. controls (-1.05,2.1) and (-1.35,1.55) .. (-1.05,1.2)
-- (-\b,1);

\node at (0,1.75) {$\Omega\sub{T}$};

\fill[green!10]
(-\b,-1)
-- (\b,-1)
-- (1.05,-1.2)
.. controls (1.35,-1.35) and (1.2,-1.85) .. (0.55,-1.8)
-- (-0.25,-2.25)
.. controls (-1.15,-2.25) and (-1.45,-1.65) .. (-1.1,-1.25)
-- (-\b,-1)
-- cycle;

\draw[green!60!black]
(-\b,-1)
-- (\b,-1)
-- (1.05,-1.2)
.. controls (1.35,-1.35) and (1.2,-1.85) .. (0.55,-1.8)
-- (-0.25,-2.25)
.. controls (-1.15,-2.25) and (-1.45,-1.65) .. (-1.1,-1.25)
-- (-\b,-1);

\node at (0,-1.75) {$\Omega\sub{B}$};

\fill (0,1) circle (1.2pt);
\fill (0,-1) circle (1.2pt);

\node[above=8pt] at (0,0.9) {$p\sub{T}$};
\node[below=8pt] at (0,-0.9) {$p\sub{\!B}$};

\fill (-\a,1) circle (1pt);
\fill (\a,1) circle (1pt);
\fill (-\a,-1) circle (1pt);
\fill (\a,-1) circle (1pt);

\node[above left=2pt] at (-\a,0.6) {$(-a,1)$};
\node[above right=2pt] at (\a,0.6) {$(a,1)$};
\node[below left=2pt] at (-\a,-0.6) {$(-a,-1)$};
\node[below right=2pt] at (\a,-0.6) {$(a,-1)$};

\draw[<->] (-\w,-0.55) -- (\w,-0.55);
\node at (0,-0.72) {$2w$};

\end{tikzpicture}
\caption{The dumbbell domain $\Omega(w)=\Omega\sub{T}\cup Q(w)\cup\Omega\sub{B}$.}
\label{fig:dumbbell}
\end{figure}
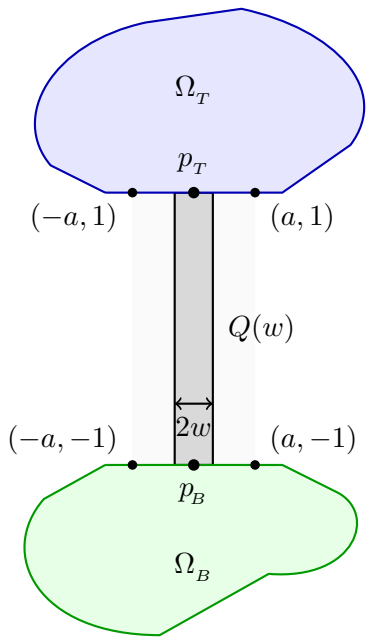

We study the nodal structure of the Neumann eigenfunctions of a planar dumbbell domain as the width of a thin neck joining two fixed end domains tends to zero (see Figure \ref{fig:dumbbell}).
Our goal is to understand how the singular degeneration of the necks constrains the number of nodal domains.

For an open, bounded domain $\Omega$ with piecewise smooth boundary, we write $0=\lambda_1<\lambda_2\leq\lambda_3\leq\cdots$ for the eigenvalues and $u_1$, $u_2$, $u_3,\ldots$  for the eigenfunctions of the Neumann problem
\begin{align*}
    (\Delta+\lambda)u= 0 \text{ in }\Omega, \qquad \pa_{n}u=0\text{ on }\pa\Omega.
\end{align*}
Here, $\pa_n$ is the outward pointing normal derivative, defined almost everywhere on $\pa\Omega$. For an eigenfunction $u$, let $\mathcal N(u)$ denote the closure of $\{x\in\Omega:u(x)=0\}$, and let $\nu(u)$ be the number of connected components of $\Omega\setminus\mathcal N(u)$, called the nodal domains of $u$. Courant's Nodal Domain Theorem \cite{courant} states that $\nu(u_j)\leq j$, and  $u_j$ is said to be Courant sharp whenever $\nu(u_j)=j$.

In dimensions two and higher, spectral position alone gives little lower control on the nodal count. For example, for the Dirichlet Laplacian on the square, there is an infinite sequence of eigenvalues tending to infinity whose eigenspaces contain eigenfunctions with exactly two nodal domains \cite{berardhelfferstern}; i.e., $\nu(u_{j_k})=2$ for some $j_k\to \infty$.

In one dimension, by contrast, every Dirichlet or Neumann eigenfunction is Courant sharp. Some of this one-dimensional rigidity persists for dumbbell domains with thin necks and Neumann boundary conditions. For sufficiently small neck widths, we prove that the neck forces the total number of nodal domains to be bounded below by a quantity of order the square root of the limiting eigenvalue. At the bottom of the spectrum, this forcing can be sharp: under the nondegeneracy assumptions of Corollary \ref{cor:3rd}, the third Neumann eigenfunction has exactly three nodal domains and is Courant sharp.

The mechanism behind this phenomenon is visible in the limiting spectrum. As the neck collapses, the Neumann spectrum of the dumbbell receives contributions from the Neumann spectra of the two end domains and from the Dirichlet spectrum of the interval corresponding to the neck. A limiting eigenvalue can therefore come from an end domain or from the neck itself. This spectral degeneration is well understood at leading order. Jimbo and Arrieta \cite{Jim,arrieta95} identify the limiting spectrum and the leading behavior of the eigenfunctions, while Gadyl'shin \cite{Gad} develops a matched-asymptotic analysis in three dimensions. These results, however, do not determine the nodal structure on an end domain where the limiting eigenfunction vanishes. Beck and Lyons \cite{beck2026nodal} obtained precise nodal information for symmetric dumbbells. Here we treat general, not necessarily symmetric, dumbbells and obtain quantitative nodal-domain bounds by identifying the first nonzero profiles on the vanishing end domains.

We now proceed to make the dumbbell geometry and the limiting spectral picture precise.

\begin{defn}[Dumbbell domain] \label{defn:dumbbell}
    Let $\Omega\sub{T}$ and $\Omega\sub{B}$ be disjoint and bounded curvilinear polygons, such that, setting $p\sub{T}=(0,1)$, $p\sub{\!B}=(0,-1)$, and for a fixed $a>0$, the boundaries of $\Omega\sub{T}$ and $\Omega\sub{B}$ satisfy
    \begin{align*}
        \{(x,1)\,:\,|x|<a\}\subset {\partial\Omega\sub{T}}, \qquad\quad  \{(x,-1)\,:\,|x|<a\}\subset {\partial\Omega\sub{B}}.
    \end{align*}
    We also assume that these domains are disjoint from the rectangle $[-a,a]\times[-1,1]$.

    For $0<w<a$, the dumbbell domain $\Omega(w)$ is defined by
    \begin{align*}
        \Omega(w)=\Omega\sub{T}\cup Q(w)\cup\Omega\sub{B},
    \end{align*}
    where the neck $Q(w)=(-w,w)\times[-1,1]$.
\end{defn}

For each $0<w<a$, we let $\{\lambda_j^w\}_{j=1}^{\infty}$ be the Neumann eigenvalues of $\Omega(w)$, with corresponding $L^2(\Omega(w))$-normalized eigenfunctions $\varphi_j^w$.

As $w\to0$, the Dirichlet problem decouples: the spectrum converges to that of $\Omega\sub{T}\cup\Omega\sub{B}$, and the corresponding eigenfunctions localize on the two end domains \cite{daners2003dirichlet}. The Neumann problem is qualitatively different: its limiting spectrum also contains modes determined by the neck length, and Neumann eigenfunctions need not vanish in the neck as it collapses \cite{Jim,arrieta95,Gad}.

Let $\{\mu_k(\Omega\sub{T})\}$ and $\{\mu_k(\Omega\sub{B})\}$  be the Neumann   spectra of the two end domains, with normalized  eigenfunctions $\psi_k^T$   and  $\psi_k^B$, respectively.  Let
\[
    \tau_n=\frac{\pi^2n^2}{4},\qquad n\geq1,
\]
 be the Dirichlet   spectrum  of $[-1,1]$, with   normalized eigenfunctions $\gamma_n(s)=\sin\!\left(\frac{n\pi(s+1)}{2}\right)$. The limiting spectrum is the ordered union
\[
    \{\mu_k(\Omega\sub{T})\}_{k\geq1}\cup\{\mu_k(\Omega\sub{B})\}_{k\geq1}\cup\{\tau_n\}_{n\geq1},
\] 
  which we denote by $\{\lambda_j\}$ and count with multiplicity.  
  
By Theorem 2.2   of \cite{arrieta95}, 
\[
\lambda_j^w\to\lambda_j \quad \text{ as} \; \;w\to 0.
\]
In particular, since $0$ is a simple Neumann eigenvalue of $\Omega(w)$, $\Omega\sub{T}$, and $\Omega\sub{B}$, we have $\lambda_1^w=0$, $\lambda_2^w>0$ with $\lim_{w\to0}\lambda_2^w=0$, and $\lim_{w\to0}\lambda_3^w=\lambda_3>0$. 

The limiting behavior of $\varphi_j^w$ also involves the eigenfunctions $\psi_k^T$, $\psi_k^B$, and $\gamma_n$. For example, the second eigenfunction $\varphi_2^w$ converges to constants of different signs in $H^1(\Omega\sub{T})$ and $H^1(\Omega\sub{B})$, or equivalently to multiples of the constant eigenfunctions $\psi_1^T$ and $\psi_1^B$.

The source of a limiting eigenvalue determines how the neck oscillations and the profiles on the end domains combine. Fix $j$ with $\lambda_{j-1}<\lambda_j<\lambda_{j+1}$. By spectral convergence, $\lambda_j^w$ is then simple for all sufficiently small $w>0$. We first consider the case in which $\lambda_j$ comes from an end domain. Relabeling the ends if necessary, write $\lambda_j=\mu_k(\Omega\sub{T})$. At leading order, $\varphi_j^w$ resembles a multiple of $\psi_k^T$ in $\Omega\sub{T}$ and a trigonometric function in $Q(w)$, while its profile in $\Omega\sub{B}$ is of lower order. We quantify the nodal contributions from $\Omega\sub{T}$ and $Q(w)$ by setting
 \begin{equation}
\label{def:mk}
    m _k: = \min\{ m  \geq1: \mu_k(\Omega\sub{T})<\tau_m \} ,
\qquad
    M_k:=  \nu(\psi_k^T).
 \end{equation}
Here, $m_k$ is the smallest index $m$ for which the $m$-th Dirichlet eigenvalue of the neck satisfies $\mu_k(\Omega_T)<\tau_m$, while $M_k$ is the number of nodal domains of the limiting eigenfunction on $\Omega\sub{T}$. 

\begin{thm} \label{thm:ends}
Suppose that $\lambda_{j-1}<\lambda_j<\lambda_{j+1}$ with $\lambda_j=\mu_k(\Omega\sub{T})$ and let $m_k$   and $M_k$ be as in \eqref{def:mk}.  Assume  that $\psi_k^T(p\sub{T})\neq0$ and   $\cos(2\sqrt{\mu_k(\Omega\sub{T})})\neq0$. Then, there exists $w_0>0$ such that, for all $0<w<w_0$,
 \begin{align}\label{eq:LUbound}
    m_k+2 \leq \nu(\varphi_j^w) \leq j-(k-M_k),
\end{align}
with exactly $m_k$ nodal domains intersecting $Q(w)$, at least two nodal domains intersecting $\Omega\sub{T}$, and at least two nodal domains intersecting $\Omega\sub{B}$.

If, in addition, $\nabla\psi_k^T\neq 0$ on the nodal set of $\psi_k^T$, then $w_0>0$ can be chosen so that, for $0<w<w_0$, $\nu(\varphi_j^w)\geq m_k+M_k$, with equality if $\mu_k(\Omega\sub{T})<\mu_2(\Omega\sub{B})$.
\end{thm}

In what follows, we will refer to the additional condition $\nabla\psi_k^T\neq0$ on the nodal set as the no-crossings condition. Under this assumption, the $M_k$ nodal domains of $\psi_k^T$ persist in the dumbbell and yield the stronger lower bound $m_k+M_k$.

The bounds in \eqref{eq:LUbound} describe how the neck oscillations and the nodal structure inherited from the end domain combine. Exactly $m_k$ nodal domains intersect the neck, while, if the nodal set of $\psi_k^T$ has no crossings, its $M_k$ nodal domains yield the sharper lower bound $m_k+M_k$. The upper bound shows that if $\psi_k^T$ falls $k-M_k$ domains short of Courant sharpness, then $\varphi_j^w$ falls at least as far short.

The condition $\psi_k^T(p\sub{T})\neq0$ ensures that the nodal set of $\psi_k^T$ is disjoint from a small neighborhood of the attachment point $p\sub{T}$; the proof shows that the same is then true for $\varphi_j^w$ for sufficiently small $w$. The condition $\cos(2\sqrt{\mu_k(\Omega\sub{T})})\neq0$ is the nondegeneracy needed to identify the first nonzero term in the asymptotic expansion on the neck and the opposite end.

Here and throughout, $w_0$ is independent of $w$, but depends on $j$, the geometry of $\Omega\sub{T}$ and $\Omega\sub{B}$, and the separation of $\lambda_j$ from the neighboring limiting eigenvalues. In particular, $w_0$ may need to be chosen smaller as $\min\{\lambda_j-\lambda_{j-1},\lambda_{j+1}-\lambda_j\}$ decreases.

We now turn to the complementary case, in which the limiting eigenvalue $\lambda_j=\tau_n$ comes from the neck. Note that the limiting neck eigenfunction $\gamma_n$ has $n$ nodal intervals in $[-1,1]$.

 \begin{thm} \label{thm:neck}
Suppose that $\lambda_{j-1}<\lambda_j<\lambda_{j+1}$, with $\lambda_j=\tau_n$. Then, there exists $w_0>0$ such that, for all $0<w<w_0$, $\nu(\varphi_j^w) \geq n+2$, with exactly $n$ nodal domains intersecting $Q(w)$, at least two nodal domains intersecting $\Omega\sub{T}$, and at least two nodal domains intersecting $\Omega\sub{B}$.

If, in addition, $\tau_n<\min\{\mu_2(\Omega\sub{T}), \mu_2(\Omega\sub{B})\}$, then $w_0>0$ can be chosen so that, for all $0<w<w_0$, the eigenfunction $\varphi_j^w$
\begin{enumerate}
    \item[i)] is Courant sharp;
    \item[ii)] has nodal set consisting of one curve in $\Omega\sub{T}$, one curve in $\Omega\sub{B}$, and $n-1$ curves in $Q(w)$.
\end{enumerate}
\end{thm}

Even though $\varphi_j^w$ converges to zero on both end domains as $w\to0$, for every sufficiently small positive $w$ at least two of its nodal domains intersect each end. In total, $\nu(\varphi_j^w)\geq n+2$. If $\tau_n<\min\{\mu_2(\Omega\sub{T}),\mu_2(\Omega\sub{B})\}$, then $\tau_n=\lambda_{n+2}$: up to this level, the limiting spectrum consists only of the two zero modes and the neck eigenvalues $\tau_1,\ldots,\tau_n$. Hence, the lower bound $n+2$ reaches Courant's upper bound.

Whether $\lambda_j$ comes from an end domain or from the neck, Theorems \ref{thm:ends} and \ref{thm:neck} give a lower bound depending only on $\lambda_j$. To state it uniformly, whenever $\lambda_j$ belongs to the spectrum of an end domain, we relabel the ends, if necessary, so that $\lambda_j=\mu_k(\Omega\sub{T})$.

\begin{cor} \label{cor:lower}
Suppose that $\lambda_{j-1}<\lambda_j<\lambda_{j+1}$, and that if $\lambda_j=\mu_k(\Omega\sub{T})$, then  $\psi_k^T(p\sub{T})\neq0$ and $\cos(2\sqrt{\mu_k(\Omega\sub{T})}) \neq0$. Then, there exists $w_0>0$ such that, for $0<w<w_0$, $\varphi_j^w$ satisfies 
\begin{align*}
\nu(\varphi_j^w)\geq \lceil2\pi^{-1}\sqrt{\lambda_j}+2\rceil\geq 3.
\end{align*}
\end{cor}

The number of nodal domains is therefore bounded below by a quantity of order $\sqrt{\lambda_j}$, whether the limiting eigenvalue comes from an end domain or from the neck. This makes precise the contrast with general domains described above.

At the first positive limiting eigenvalue, this lower bound reaches Courant's upper bound. If $\lambda_3<\lambda_4$, then $\lambda_3$ is either the smaller of the two first positive end-domain eigenvalues or the first neck eigenvalue $\tau_1$. After relabeling the ends if necessary, the main theorems give the following conclusion.

\begin{cor} \label{cor:3rd}
    Assume that $\lambda_3<\lambda_4$ and suppose that if $\lambda_3=\mu_2(\Omega\sub{T})$, then  $\psi_2^T(p\sub{T})\neq0$ and $\cos(2\sqrt{\mu_2(\Omega\sub{T})}) \neq0$. Then, there exists $w_0>0$ such that the eigenfunction $\varphi_3^w$ has exactly three nodal domains in $\Omega(w)$, for $0<w<w_0$. Moreover, the nodal set of $\varphi_3^w$ is disjoint from the neck $Q(w)$ and consists of two curves, one in $\Omega\sub{T}$ and the other in $\Omega\sub{B}$.
\end{cor}

The second and third eigenfunctions consequently have sharply different nodal geometry. The nodal set of the third eigenfunction reaches each end domain, whereas the nodal set of the second eigenfunction is expected to remain inside the neck; \cite{MS23} proves that it cannot extend far from it.

  The upper bound in Theorem \ref{thm:ends} also gives the following criterion for failure of Courant sharpness.

 \begin{cor} \label{cor:4th}
Suppose that $\lambda_{j-1}<\lambda_j <\lambda_{j+1}$, with $\lambda_j = \mu_k(\Omega\sub{T})$, and that $\psi_k^T(p\sub{T})\neq0$ and $\cos(2\sqrt{\mu_k(\Omega\sub{T})}) \neq0$. If $\psi_k^T$ is not Courant sharp in $\Omega\sub{T}$, then there exists $w_0>0$ such that $\varphi_j^w$ is not Courant sharp in $\Omega(w)$ for every $0<w<w_0$.
\end{cor}

A nodal deficit on an end domain can therefore persist after the neck is attached and prevent Courant sharpness on the full dumbbell. The fourth eigenfunction already displays both possibilities. If $\lambda_4=\tau_2<\min\{\mu_2(\Omega\sub{T}),\mu_2(\Omega\sub{B})\}$, Theorem \ref{thm:neck} gives  four nodal domains and Courant sharpness. If instead $\lambda_4=\mu_3(\Omega\sub{T})$, the hypotheses of Theorem \ref{thm:ends} hold, and  $\psi_3^T$ has exactly two nodal domains, then $\nu(\varphi_4^w)=3$ for all sufficiently small $w$, so $\varphi_4^w$  is not Courant sharp.  In the symmetric case $\Omega\sub{T}=\Omega\sub{B}$, if $\mu_2(\Omega\sub{T})<\tau_1$ and $\psi_2^T(p\sub{T})\neq0$, then \cite{beck2026nodal} shows that $\varphi_3^w$ and $\varphi_4^w$ have three and four nodal domains, respectively, for all sufficiently small $w>0$.

The Courant-sharp cases have an additional variational interpretation: their nodal domains form minimal spectral partitions for the maximum of the first mixed Dirichlet--Neumann eigenvalues of the pieces \cite{bonnaillie2015nodal,HH13}. For a fixed planar domain, only finitely many Neumann eigenfunctions can be Courant sharp \cite{pleijel1956,polterovich2009,lena-pleijel,BCM-pleijel}; moreover, such eigenfunctions have been identified only in special settings, often with substantial symmetry; see, for example, \cite{helffer2010spectral,lena2015courant,berard2016courant,BBF17}.

To prove Theorems \ref{thm:ends} and \ref{thm:neck}, we begin with Arrieta's estimates \cite{arrieta95}, which identify the leading-order behavior of the eigenfunctions: an end-domain mode converges to $\psi_k^T$ on its active end and to zero on the opposite end, while a neck mode converges, after the natural normalization, to $\gamma_n$ in the neck and to zero on both ends. These estimates, however, give no nodal information on an end domain where the limiting eigenfunction vanishes. Our main analytic step is to identify the first nonzero profile on each such end domain. We show that these profiles are given by Neumann Green's functions with poles at the attachment points $p\sub{T}$ and $p\sub{\!B}$. To obtain them, we develop a matched-asymptotic expansion near the attachment points, extending the three-dimensional analysis of Gadyl'shin \cite{Gad} to two dimensions. The logarithmic singularity of the two-dimensional Green's function introduces powers of $\log w$ that are absent in the three-dimensional expansion. Combining these refined asymptotics with nodal estimates for the Green's functions and the techniques of \cite{beck2026nodal} yields Theorems \ref{thm:ends} and \ref{thm:neck}.

\subsection{Nodal counts and the graph structure of chain domains}

Our longer-term goal is to understand how much of the geometry and topology of a
chain domain can be recovered from the nodal structure of its low-energy
Neumann eigenfunctions. The dumbbell considered in this paper is the simplest
nontrivial example. More generally, a chain domain consists of $M$ bounded
subdomains $\{D_i\}_{i=1}^M$ joined by $K$ thin necks
$\{Q_r(w)\}_{r=1}^K$. As the neck widths shrink, the geometry separates into
the individual subdomains together with the combinatorial information
describing how they are connected. The question is whether the nodal counts of
the continuum eigenfunctions retain this information.

Replacing the subdomains by vertices and the necks by edges produces an
associated graph. The weights in the corresponding graph Laplacian reflect the
geometry of the necks and subdomains, while the cycle structure records how the
pieces are connected. This graph already governs the low-energy spectrum: the
first $M$ eigenvalues converge to zero at rates determined by the weighted graph
Laplacian, and the corresponding graph eigenvectors determine the limiting
constant values of the continuum eigenfunctions on the subdomains \cite{MJim}. At high
energy, by contrast, we proved that Neumann eigenfunctions of chain domains are
not Courant sharp uniformly as the neck widths shrink \cite{BCM-pleijel}. The present results
suggest that the low-energy regime contains a different kind of rigidity, in
which nodal counts may detect the topology of the limiting graph.

The discrete nodal theory gives a precise model for what one should expect. For
a graph eigenvector with no zero entries, its nodal domains are the maximal
connected sets of vertices on which the entries have the same sign. For a path,
Gantmacher and Krein \cite{GK02} show that the $j$-th eigenvector has exactly
$j$ nodal domains. More generally, under the usual simplicity and nonvanishing
assumptions, Berkolaiko \cite{Ber07} gives the lower bound $j-\ell$, where
$\ell$ is the cycle rank of the graph; see also \cite{Biy03} for trees. Band
\cite{Band12} proves a converse: if every eigenvector attains the index bound,
then the graph is a tree. These results suggest that tree topology should be
visible in the nodal counts of the first continuum eigenfunctions.

Since the first $M$ continuum eigenfunctions are governed by this graph limit,
one expects them to inherit the discrete oscillation law and hence to be
Courant sharp when the limiting graph is a tree. Our dumbbell results show that
the phenomenon persists one step further: when $M=2$, the third eigenfunction
is also Courant sharp. At the next eigenvalue, however, the geometry and nodal
structure of the individual end domains begin to matter. This leads to the
following conjecture.

The nondegeneracy assumptions below refer to simplicity of the weighted graph
spectrum, nonvanishing of its eigenvectors, isolation of the first positive
limiting eigenvalue, and nondegeneracy of the corresponding data at the neck
attachment points.

\begin{conjecture}
\label{conj:1}
Suppose that the metric graph associated with a chain domain is a tree and that
the preceding nondegeneracy assumptions hold. Then, for sufficiently small neck
widths, the first $M+1$ Neumann eigenfunctions of the chain domain are Courant
sharp. Courant sharpness of the $(M+2)$-nd Neumann eigenfunction is not
determined by the tree topology alone, but also depends on the spectral and
nodal structure of the subdomains $\{D_i\}_{i=1}^{M}$.
\end{conjecture}

For $M=2$, Corollary \ref{cor:3rd} verifies the first assertion: the third
eigenfunction, with index $M+1$, is Courant sharp. The next eigenfunction can
exhibit either behavior. Corollary \ref{cor:4th} gives conditions under which
the fourth eigenfunction is not Courant sharp, while \cite{beck2026nodal} gives
symmetric examples in which it is Courant sharp.

The conjecture would also allow nodal counts to be used in the reverse
direction, to constrain the structure of the chain. Under the same
nondegeneracy assumptions, if the $N$-th Neumann eigenfunction is not Courant
sharp for all sufficiently small neck widths, then either the associated
metric graph is not a tree or the chain domain has at most $N-2$ subdomains.
Indeed, if the graph were a tree and the chain had $M\geq N-1$ subdomains, then
$N\leq M+1$, and Conjecture \ref{conj:1} would predict that the $N$-th
eigenfunction is Courant sharp. In this sense, failure of Courant sharpness at
low energy would give information not only about the presence of cycles, but
also about the number of subdomains in the chain.

Our numerical experiment with two domains joined by two necks illustrates the
role of the tree hypothesis (See Section \ref{sec:numExp}). The associated graph has a cycle, and the third
eigenfunction need not be Courant sharp. Together with the dumbbell results,
this gives the first evidence for the broader principle that low-energy nodal
counts can detect both the topology of the limiting graph and the number of
subdomains in the chain.

\subsection{Outline of the paper}

The remainder of the paper is organized as follows. First, in Section \ref{sec:Green}, we define the Green's functions for $\Omega\sub{T}$ and $\Omega\sub{B}$ that will play a role when approximating the eigenfunctions of the dumbbell. In Propositions \ref{prop:Green1} and \ref{prop:Green2}, we establish their behavior near the points $p\sub{T}$ and $p\sub{\!B}$ where the neck of the dumbbell is attached, including their logarithmic singularities, and obtain estimates on their number of nodal domains. Then, in Section \ref{sec:approx}, we give the leading-order approximation to the eigenfunction, both when the eigenvalue converges to a Neumann eigenvalue of $\Omega\sub{T}$ and when it converges to a Dirichlet eigenvalue coming from the neck (see Propositions \ref{prop:approx1} and \ref{prop:approx2}). We then, in Section \ref{sec:proof}, use these approximation results for eigenfunctions of $\Omega(w)$ to study the location of their nodal sets and prove the two theorems. Finally, in Section \ref{sec:matched}, we prove the approximation results of Propositions \ref{prop:approx1} and \ref{prop:approx2}. These proofs involve a delicate matched asymptotics argument, adapting the work of Gadyl'shin \cite{Gad} to the two-dimensional case, where the expansions involve powers of $\ln(w)$ due to the logarithmic singularity of the Green's functions. The final section, Section \ref{sec:numExp}, presents numerical experiments illustrating the dumbbell results and the chain-domain conjecture.

\subsection*{Acknowledgments} 
 Y.C. was supported by NSF CAREER Grant DMS-2045494. J.L.M. acknowledges support from the NSF through NSF FRG grant DMS-2152289 and NSF Applied Math Grant DMS-2307384.

\section{The Green's function} \label{sec:Green}

We will approximate the eigenfunctions $\varphi_j^w$ in $\Omega\sub{T}$ and $\Omega\sub{B}$ in terms of the limit of a sequence of Green's functions in $\Omega\sub{T}$ and $\Omega\sub{B}$ as the pole approaches the boundary at $p\sub{T}$ and $p\sub{\!B}$. Therefore, it will be important to obtain information about the number of nodal domains of these limiting Green's functions. We start by giving their precise definition. 
\begin{defn}
For $y\in \Omega\sub{T}$, and $\mu>0$ not in the Neumann spectrum of $\Omega\sub{T}$, we define $G^T(x,y,\mu)$ to be the Neumann Green's function for the operator $\Delta + \mu$ in $\Omega\sub{T}$.

For $p\sub{T} = (0,1)\in\pa\Omega\sub{T}$, we let, for all large enough $j$, $p_j = (0,1+j^{-1})\in\Omega\sub{T}$, and define $G^T(x,p\sub{T},\mu)$ as the pointwise limit
\begin{align*}
    G^T(x,p\sub{T},\mu) = \lim_{j\to\infty}G^T(x,p_j,\mu).
\end{align*}
{Note, this limit is well-defined pointwise for each $x \in \Omega\sub{T}$ as it is a limit of smooth, bounded functions in a small neighborhood of each such $x$.}
We also define 
$$\pa_{y_1}^\ell G^T(x,p\sub{T},\mu) = \lim_{j\to\infty}\pa_{y_1}^\ell G^T(x,y,\mu)|_{y=p_j}$$ 
for each $\ell\geq1$. The functions $G^B(x,y,\mu)$, $G^B(x,p\sub{\!B},\mu)$, and $\pa_{y_1}^\ell G^B(x,p\sub{\!B},\mu)$ are defined analogously.  
\end{defn}

 Since we will use $G^T(x,p\sub{T},\mu)$ for $\mu$ close (but not equal) to a Neumann eigenvalue $\mu_k(\Omega\sub{T})$, it will be important to have bounds with explicit dependence on $\mu-\mu_k(\Omega\sub{T})$. We will denote $D_r(p\sub{T})$ to be the disc centered at $p\sub{T}$ of radius $r>0$. Note that by Definition \ref{defn:dumbbell}, for $r<a$, the set $\pa\Omega\sub{T} \cap D_r(p\sub{T})$ is a horizontal line segment. In Lemma 3.1 in \cite{Gad}, the three dimensional case is considered, and the behavior of the limiting Green's function near the pole is established. In our two dimensional case, the limiting Green's function $G^T(x,p\sub{T},\mu)$ has the following properties (and analogously for $G^B(x,p\sub{\!B},\mu)$). 
\begin{prop} \label{prop:Green1}
Let $\psi_1, \ldots,\psi_n$ be an orthonormal basis of Neumann eigenfunctions of $\Omega\sub{T}$ with eigenvalue $\mu_k(\Omega\sub{T})$, and let $c^*_k$ be the distance of $\mu_k(\Omega\sub{T})$ from the rest of the Neumann spectrum of $\Omega\sub{T}$. Then, for $|\mu-\mu_k(\Omega\sub{T})|<\tfrac{1}{2}c^*_k$, the limiting Green's function $G^T(x,p\sub{T},\mu)$ has the following properties (and analogously for $G^B(x,p\sub{\!B},\mu)$).
    \begin{enumerate}
        \item[1)] $(\Delta+\mu)\pa_{y_1}^\ell G^T(x,p\sub{T},\mu) = 0$ for $x\in \Omega\sub{T}$;
        \item[2)] $\pa_{\nu}\pa_{y_1}^\ell G^T(x,p\sub{T},\mu) = 0$ on $\pa\Omega\sub{T}\backslash\{p\sub{T}\}$, whenever the outward pointing unit normal $\nu$ is defined;
        \item[3)] $\pa_{y_1}^\ell G^T(x,p\sub{T},\mu)$ can be written as
        \begin{align} \label{eqn:Green1-a}
            \pa_{y_1}^\ell G^T(x,p\sub{T},\mu) = \sum_{m=1}^{n} \frac{\psi_m(x)(\pa_{y_1}^\ell \psi_m)(p\sub{T})}{\mu-\mu_k(\Omega\sub{T})} - \frac{(-1)^\ell i}{2}\pa_{x_1}^\ell H_0(\sqrt{\mu}|x-p\sub{T}|) + g^T_{\ell}(x,\mu).
        \end{align}
        Here $H_0$ is the zero-th Hankel function of the first kind. The functions $g^T_{\ell}(x,\mu)$ are smooth for $x$ in the interior of $\Omega\sub{T}$, and smooth up to the boundary in $D_{a/2}(p\sub{T})\cap \bar{\Omega}_T$, with bounds depending only on $c_k^*$, an upper bound on $\mu$, and the geometry of $\Omega\sub{T}$. On $\pa\Omega\sub{T}\cap D_{a/2}(p\sub{T})$, the functions $g^T_{\ell}(x,\mu)$ satisfy $\pa_{\nu}g_\ell^T(x,\mu) = -\pa_{x_2}g_\ell^T(x,\mu) = 0$.
    \end{enumerate}
\end{prop}
\begin{remark} \label{rem:Green}
    Since the Hankel function can be written in terms of Bessel functions as $H_0(z) = J_0(z)+iY_0(z)$, with $J_0$ non-singular at $z=0$ and $Y_0(z)$ with a $\tfrac{2}{\pi}\ln|z|$ singularity at $z=0$, the singularity of $G^T(x,p\sub{T},\mu)$ at $x=p\sub{T}$ is given by $\tfrac{1}{\pi}\ln|x-p\sub{T}|$. In particular, there exists $r>0$ such that $G^T(x,p\sub{T},\mu)<0$ in $\Omega\sub{T}\cap D_r(p\sub{T})$.

     Working in complex notation (with $\emph{Re}(x)=x_1$ and $\emph{Im}(x)=x_2$), we have $\ln|x-p\sub{T}| = \emph{Re}(\ln(x-p\sub{T}))$, and as a function of $x$, $G_0(x) = \ln(x-p\sub{T})$ and, for $\ell\geq1$, $G^{(\ell)}_0(x) = \frac{(-1)^{\ell-1}(\ell-1)!}{(x-p\sub{T})^\ell}$ are analytic for Im$(x-p\sub{T})>0$. Therefore, 
    \begin{align*}
        \pa_{x_1}^{\ell}\ln|x-p\sub{T}| =   \pa_{x_1}^{\ell}\emph{Re}(\ln(x-p\sub{T})) = \emph{Re}\left(\frac{(-1)^{\ell-1}(\ell-1)!}{(x-p\sub{T})^\ell}\right).
    \end{align*}
    We have the analogous statement for derivatives of $\ln|x-p\sub{\!B}|$.
\end{remark}
\begin{remark} \label{rem:Green2}
    If $\mu$ is a distance $c>0$ from the whole Neumann spectrum of $\Omega\sub{T}$, then one has the same decomposition for $\pa_{y_1}^\ell G^T(x,p\sub{T},\mu)$ as in \eqref{eqn:Green1-a}, now without the sum on the right hand side, with bounds on  $g_\ell^T(x,\mu)$ depending on $c$, an upper bound on $\mu$, and the geometry of $\Omega\sub{T}$.
\end{remark}
As the Green's functions will be used to approximate the Neumann eigenfunctions $\varphi_j^w$ of $\Omega(w)$ for small $w>0$, an estimate on their nodal domain count will be crucial.
\begin{prop} \label{prop:Green2}
The Green's function $G^T(x,p\sub{T},\mu)$ satisfies the following properties (and analogously for $G^B(x,p\sub{\!B},\mu)$).
\begin{enumerate}
  \item[1)] If $0<\mu <\mu_k(\Omega\sub{T})$ and $\mu$ is not in the Neumann spectrum of $\Omega\sub{T}$, then $G^T(x,p\sub{T},\mu)$ has at most $k$ nodal domains in $\Omega\sub{T}$.
  \item[2)] $G^T(x,p\sub{T},\mu)$ satisfies
  \begin{align*}
      \int_{\Omega\sub{T}} G^T(x,p\sub{T},\mu)\,dx = \frac{1}{\mu}>0,
  \end{align*}
  and always has at least 2 nodal domains
\end{enumerate}
In particular, if $0<\mu<\mu_2(\Omega\sub{T})$, then $G^T(x,p\sub{T},\mu)$ has exactly 2 nodal domains.
\end{prop}
We end this section by proving the two propositions. To prove Proposition \ref{prop:Green1}, we use a similar strategy to the proof of Lemma 3.1 in \cite{Gad}, by first isolating the singular part of the Green's function via the fundamental solution, and then obtaining elliptic estimates on the remainder. 
\begin{proof1}{Proposition \ref{prop:Green1}}
The first two properties of the proposition follow immediately from the definition of $\pa_{y_1}^{\ell}G^T(x,p\sub{T},\mu)$, since, for sufficiently large $j$, $G^T(x,p_j,\mu)$ is smooth for $x$ in any compact subset of the interior of $\Omega\sub{T}$, and smooth up to the boundary in any compact subset of $(\bar{\Omega}_T\cap D_a(p\sub{T}))\backslash \{p\sub{T}\}$, where {$a$ is taken as in Definition \ref{defn:dumbbell}.}

To prove 3), we adapt the proof of Lemma 3.1 in \cite{Gad} to the two dimensional case. For ease of notation, we drop the T's throughout, and after a translation, we assume that $p = p\sub{T}$ is at the origin. Throughout the proof, $\mu$ will satisfy $|\mu-\mu_k(\Omega)|<\tfrac{1}{2}c_k^*$, and any bounds will depend only on $c_k^*$, an upper bound on $\mu$, and the geometry of $\Omega$.  A fundamental solution for the Helmholtz operator $\Delta + \mu$ in two dimensions is given by $-\frac{i}{4}H_0(\sqrt{\mu}|x-y|)$, and $\{\psi_m\}_{m=1}^{n}$ is an orthonormal basis for the eigenspace of $\mu_k(\Omega)$. Therefore, for $x,y\in\Omega$, we can write
    \begin{align*}
        G(x,y,\mu) = \sum_{m=1}^{n}\frac{\psi_m(x)\psi_m(y)}{\mu-\mu_k(\Omega)} -\frac{i}{4}H_0(\sqrt{\mu}|x-y|) + g(x,y,\mu),
    \end{align*}
    with $g(x,y,\mu)$ a smooth function in $x,y\in\Omega$, and analytic in $\mu$. Moreover, since the boundary of $\Omega$ coincides with a horizontal line segment in a neighborhood of the origin, and setting $y^* = (y_1,-y_2)$, we write
    \begin{align} \label{eqn:Green1a}
          G(x,y,\mu) = \sum_{m=1}^{n}\frac{\psi_m(x)\psi_m(y)}{\mu-\mu_k(\Omega)} -\frac{i}{4}\left(H_0(\sqrt{\mu}|x-y|)+H_0(\sqrt{\mu}|x-y^*|)\right)\left(1-\chi(2|x|/a)\right) + \tilde{G}(x,y,\mu),
    \end{align}
    for $y\in D_a(0)$. Here $\chi(\cdot)$ is a smooth cut-off function, with $\chi(t)=1$ for $t>2$ and $\chi(t) = 0$ for $t<1$. Note that for $x\in\pa\Omega$ with $1-\chi(2|x|/a) \neq0$, we have $x = (x_1,0)$ and $\pa_{\nu} = -\pa_{x_2}$, and, for $\ell\geq0$, 
    \begin{align*}
-\pa_{x_2}\pa_{y_1}^{\ell}\left((H_0(\sqrt{\mu}|x-y|)+H_0(\sqrt{\mu}|x-y^*|))\left(1-\chi(2|x|/a)\right)\right)\big|_{x=(x_1,0)} = 0.
    \end{align*}

    Note also that $\psi_m(x)\psi_m(y)$ is smooth, uniformly up to the boundary for $x,y\in \Omega\cap D_a(0)$. Therefore, for each such $y$, and for $\ell\geq0$,
\begin{align*}
    (\Delta_x+\mu) \pa_{y_1}^{\ell}\tilde{G}(x,y,\mu)=\tilde{F}_\ell(x,y,\mu) \, \text{ for }x\in\Omega, \quad \pa_{\nu_x}\pa_{y_1}^{\ell}\tilde{G}(x,y,\mu) = 0 \,\text{ for }x\in\pa \Omega,
\end{align*}
where $\tilde{F}_\ell(x,y,\mu)$ is smooth in $x\in \Omega$ and $y\in \Omega\cap D_a(0)$, with bounds uniform up to the boundary in $x, y\in \Omega\cap D_a(0)$, and analytic in $\mu$. 

Therefore, letting $y=p_j = (0,1/j)$ be points in $\Omega\cap D_a(0)$, converging to $0$, and letting $\tilde{G}_{\ell}(x,\mu)$ be the pointwise limit $\lim_{j\to\infty}\pa_{y_1}^\ell\tilde{G}(x,y,\mu)|_{y=p_j}$, we obtain functions $\tilde{G}_\ell(x,\mu)$ satisfying
\begin{align*}
    (\Delta+\mu) \tilde{G}_\ell(x,\mu)=\tilde{F}_\ell(x,\mu) \text{ in }\Omega, \quad \pa_{\nu_x}\tilde{G}_\ell(x,\mu) = 0 \text{ on }\pa \Omega.
\end{align*}
Here the functions $\tilde{F}_\ell(x,\mu)$ are smooth for $x\in \Omega$, with bounds uniform up to the boundary for $x\in \Omega\cap D_a(0)$, and analytic in $\mu$.  We can now conclude  that the functions $\tilde{G}_\ell(x,\lambda)$ are smooth for $x\in\Omega$, analytic in $\mu$, with bounds uniform up to the boundary for $x\in \Omega\cap D_a(0)$,  exactly as for the functions $\tilde{G}_j$ in the second half of the proof of Lemma 3.1 on page 280 in \cite{Gad}.

We finally apply $\pa_{y_1}^{\ell}$ to \eqref{eqn:Green1a}, set $y=p_j$, and let $j\to\infty$ (noting that $y^* = p_j^*\to 0$ also). This, together with the above estimates on $\tilde{G}_\ell(x,\lambda)$ gives \eqref{eqn:Green1-a} in 3) of the proposition, and the required bounds on the functions $g_\ell(x,\mu)$.
\end{proof1}

To prove Proposition \ref{prop:Green2}, we first use Proposition \ref{prop:Green1} to remove the unique nodal domain of $G^T(x,p\sub{T},\mu)$ containing $p\sub{T}$ in its closure, and then use the proof of the Courant nodal domain theorem to obtain an upper bound on the number of nodal domains remaining.
\begin{proof1}{Proposition \ref{prop:Green2}}
    For ease of notation, we drop the $T$'s throughout. We start with the first claim. The min-max principle states that
$$
\mu_m(\Omega)
= \inf_{\substack{E \subset H^1(\Omega)\\ \dim E = m}}
\; \sup_{\substack{ u \in E\\u\ne 0}}\; Q_\Omega(u),
\qquad
Q_\Omega(u) := \frac{\int_\Omega |\nabla u|^2}{\int_\Omega |u|^2}.
$$
Fix $k \ge 2$ and $\mu \in (0,\mu_k(\Omega))$. Let
$
\Omega_0, \Omega_1, \dots, \Omega_{N-1}
$
denote the connected nodal domains of $G(\cdot,p,\mu)$ in $\Omega$, where $\Omega_0$ is the nodal domain satisfying $p\in\bar{\Omega}_0$. (By Remark \ref{rem:Green}, this is the unique $\Omega_\ell$ for which $p\in\bar{\Omega}_\ell$.)  Thus $G(\cdot,p,\mu)$ does not change sign in each $\Omega_\ell$ and vanishes on $\partial \Omega_\ell \cap \Omega$.

Set
$
\widetilde{\Omega} := \overline{\Omega} \setminus \Omega_0
$
and, with $\Gamma := \partial \Omega_0 \cap \Omega$,  consider on $\widetilde{\Omega}$ the mixed problem
$$
\begin{cases}
-\Delta u = \lambda u & \text{in } \widetilde{\Omega},\\
u = 0 & \text{on } \Gamma,\\
\partial_\nu u = 0 & \text{on } \partial \Omega \cap \partial \widetilde{\Omega}.
\end{cases}
$$
Let
$
0 < \lambda_1(\widetilde{\Omega}) \le \lambda_2(\widetilde{\Omega}) \le \cdots
$
denote its eigenvalues. 
Let 
$
\operatorname{Tr}: H^1(\widetilde \Omega) \to H^{1/2}(\partial \widetilde \Omega)
$
be the trace operator and define
$
H^1_\Gamma(\widetilde \Omega)
:= \{\, u \in H^1(\widetilde \Omega) : \operatorname{Tr}(u) = 0 \text{ on } \Gamma \,\}.
$
This is the natural form domain for the mixed problem on $\widetilde \Omega$. Then, 
$$
\lambda_m(\widetilde \Omega)
= \inf_{\substack{E \subset H^1_\Gamma(\widetilde \Omega)\\ \dim E = m}}
\ \sup_{\substack{ u \in E\\u\ne 0}}
Q_{\widetilde \Omega}(u).
$$
Every $u \in H^1_\Gamma(\widetilde \Omega)$ can be extended  to a function $\tilde u$ on $\Omega$ by setting it to be $0$ on $\Omega_0$.
Then, $\tilde u \in H^1(\Omega)$ 
and 
$Q_\Omega(\tilde u) = Q_{\widetilde \Omega}(u)$
for all $u \in H^1_\Gamma(\widetilde \Omega)\setminus\{0\}.$
From this observation it follows that
$
\lambda_m(\widetilde \Omega) \ge \mu_m(\Omega)
$
for all $m \ge 1.$
In particular, our assumption $\mu < \mu_k(\Omega)$ yields
\begin{equation}\label{eq:mu-less-than-lambda-k}
\mu < \mu_k(\Omega) \le \lambda_k(\widetilde{\Omega}).
\end{equation}

Next, for each $\ell \ge 1$ define
$$
v_\ell(x) :=
\begin{cases}
G(x,p,\mu), & x \in \Omega_\ell,\\
0, & x \in \widetilde{\Omega} \setminus \Omega_\ell.
\end{cases}
$$
Then $v_\ell \in H^1(\widetilde \Omega)$. Since $G(\cdot,p,\mu)$ vanishes on the interior nodal set $\partial\Omega_\ell \cap \Omega$ and we define $v_\ell$ by patching $G$ with $0$ across this set, the trace of $v_\ell$ on $\partial\Omega_\ell \cap \Omega$ is zero. On the other hand, on $\partial\Omega_\ell \cap \partial \Omega$ we have $v_\ell = G(\cdot,p,\mu)$, so $\partial_\nu v_\ell = \partial_\nu G(\cdot,p,\mu) = 0$ there by the Neumann boundary condition for $G$.

Moreover, the singularity of $G(\cdot,p,\mu)$ lies in $\bar{\Omega}_0$, so in each $\Omega_\ell$ with $\ell \ge 1$ we have
$$
(\Delta + \mu)v_\ell = 0 \quad \text{in } \Omega_\ell.
$$
Multiplying this equation by $v_\ell$ and integrating over $\Omega_\ell$, and using the two boundary facts above, we obtain
$$
\int_{\Omega_\ell} |\nabla v_\ell|^2 \, dx = \mu \int_{\Omega_\ell} v_\ell^2 \, dx.
$$

Since $v_\ell$ is supported in $\Omega_\ell \subset \widetilde{\Omega}$, we have
$
Q_{\widetilde{\Omega}}(v_\ell)= \mu.
$

Next, assume for contradiction that $N > k$, so there are at least $k$ nodal domains disjoint from $\Omega_0$. Then $v_1, \dots, v_k$ have pairwise disjoint supports and are therefore linearly independent in $H^1(\widetilde{\Omega})$. Let
$$
E_k := \operatorname{span}\{v_1, \dots, v_k\}.
$$
For $v = \sum_{\ell=1}^k a_\ell v_\ell \in E_k$, the disjointness of the supports gives
$$
Q_{\widetilde \Omega}(v)
= \frac{\displaystyle \int_{\widetilde \Omega} |\nabla v|^2}{\displaystyle \int_{\widetilde \Omega} |v|^2}
= \frac{\displaystyle \sum_{\ell=1}^k a_\ell^2 \int_{\widetilde \Omega} |\nabla v_\ell|^2}{\displaystyle \sum_{\ell=1}^k a_\ell^2 \int_{\widetilde \Omega} v_\ell^2}
= \frac{\displaystyle \sum_{\ell=1}^k a_\ell^2 \mu \int_{\widetilde \Omega} v_\ell^2}{\displaystyle \sum_{\ell=1}^k a_\ell^2 \int_{\widetilde \Omega} v_\ell^2}
= \mu.
$$
Hence $\sup_{\substack{v \in E_k\\v\neq 0}} Q_{\widetilde \Omega}(v) = \mu$, and the min–max principle yields $\lambda_k(\widetilde \Omega) \le \mu$, contradicting \eqref{eq:mu-less-than-lambda-k}. Thus $N \le k$, so $G(\cdot,p,\mu)$ has at most $k$ nodal domains in $\Omega$, proving item 1).

We now prove item 2). For $p_j = (0,1+j^{-1})\in\Omega$, the spectral expansion of $G(\cdot,p_j,\mu)$ reads
$$
G(x,p_j,\mu) = \sum_{n=1}^\infty \frac{1}{\mu - \mu_n(\Omega)} \, \psi_n(x)\,\psi_n(p_j)
$$
in $L^2(\Omega)$. The logarithmic singularity at $x=p_j$ is integrable, so we may integrate term by term to obtain
$$
\int_\Omega G(x,p_j,\mu) \, dx
= \sum_{n=1}^\infty \frac{1}{\mu - \mu_n(\Omega)} \, \psi_n(p_j)\,\int_\Omega \psi_n(x)\,dx.
$$
Since $\{\psi_n\}$ is an orthonormal basis in $L^2(\Omega)$ and $\psi_1 = (\text{Area}(\Omega))^{-1/2}$ is constant, we have
$
\int_\Omega \psi_n(x)\,dx = 0
$
for all $n \ge 2.$
Hence,
$$
\int_\Omega G(x,p_j,\mu) \, dx
= \frac{1}{\mu - \mu_1(\Omega)} \, \psi_1(p_j)\,\int_\Omega \psi_1(x)\,dx.
$$
Because $\mu_1(\Omega) = 0$ and $\psi_1 = (\text{Area}(\Omega))^{-1/2}$, this gives 
$$
\int_\Omega G(x,p_j,\mu) \, dx = \frac{1}{\mu} > 0.
$$
By the decomposition \eqref{eqn:Green1a} and the uniform bounds on its regular part, there is a constant $C$, independent of $j$, such that $|G(x,p_j,\mu)|\leq C(1+|\ln|x-p_j||)$. Since $p_j=(0,1/j)$, we have $|x-p_j|\geq |x_1|$. Hence, near the singularity, $|\ln|x-p_j||\leq |\ln|x_1||$, while away from the singularity the logarithmic terms are uniformly bounded. Therefore, $|G(x,p_j,\mu)|\leq C(1+|\ln|x_1||)$, where the right-hand side belongs to $L^1(\Omega)$. The dominated convergence theorem then gives 
$$
\int_\Omega G(x,p,\mu) \, dx = \lim_{j\to\infty}\int_\Omega G(x,p_j,\mu) \, dx = \frac{1}{\mu} > 0,
$$
as required.

As noted in Remark \ref{rem:Green}, $G(x,p,\mu)<0$ in $\Omega\cap D_r(p)$ for sufficiently small $r>0$, but we have now shown
 that $G(\cdot,p,\mu)$ has positive average over $\Omega$. It follows that $G(\cdot,p,\mu)$ must be positive at some point in $\Omega$, so it changes sign in $\Omega$, completing the proof of 2).

Finally, if $0 < \mu < \mu_2(\Omega)$, then item 1) with $k = 2$ shows that $G(\cdot,p,\mu)$ has at most two nodal domains in $\overline{\Omega} \setminus \{p\}$, while item 2) shows that it has at least two. Hence $G(\cdot,p,\mu)$ has exactly two nodal domains in $\overline{\Omega} \setminus \{p\}$.
\end{proof1}

\section{An approximation to eigenfunctions of $\Omega(w)$} \label{sec:approx}

The purpose of this section is to identify the first nonzero profile of $\varphi_j^w$ in each part of the dumbbell. This refinement is needed because the leading-order convergence theory becomes trivial on at least one end domain, precisely where nodal information remains to be recovered. We begin by recalling the leading-order approximation of Arrieta \cite{arrieta95}. If $\lambda_{j-1}<\lambda_j<\lambda_{j+1}$, then $\varphi_j^w$ is approximated in the three components of $\Omega(w)$ as follows.
\begin{enumerate}
    \item[] \textbf{Case 1.} $\lambda_j = \mu_k(\Omega\sub{T})$. The eigenfunction $\varphi_j^w$ satisfies
    \begin{align*}
        \norm{\varphi_j^w - \psi_k^T}_{H^1(\Omega\sub{T})} +  \norm{\varphi_j^w}_{H^1(\Omega\sub{B})} +  \norm{\varphi_j^w - \xi}_{H^1(Q(w))} = o(w^{1/2}). 
    \end{align*}
    Here $\xi(x_2)$ satisfies $\xi'' + \mu_k(\Omega\sub{T})\xi=0$, with $\xi(1) = \psi_k^T(p\sub{T})$, $\xi(-1) = 0$.
    
    \item[] \textbf{Case 2.} $\lambda_j = \tau_n$. The eigenfunction $\varphi_j^w$ satisfies
    \begin{align*}
        \norm{\varphi_j^w}_{H^1(\Omega\sub{T})} +  \norm{\varphi_j^w}_{H^1(\Omega\sub{B})} +  \norm{\varphi_j^w - (2w)^{-1/2}\gamma_n}_{H^1(Q(w))} = o(w^{1/2}|\ln(w)|^{1/2}). 
    \end{align*}
\end{enumerate}
Here, and throughout, the implicit constants in $o$ or $O$ depend on $j$, the geometry of $\Omega\sub{T}$ and $\Omega\sub{B}$, and $\min\{\lambda_j-\lambda_{j-1},\lambda_{j+1}-\lambda_j\}$, but are independent of $w$.

These estimates identify the leading profile, but in each case the eigenfunction converges to zero in at least one end domain. Hence, they do not determine the sign structure there and are not sufficient for counting nodal domains. We therefore refine Arrieta's approximation by identifying the first nonzero term precisely on those components where the leading profile vanishes.

To state these refined approximations, we assume that $\lambda_{j-1}<\lambda_j<\lambda_{j+1}$, and split into two cases depending on whether $\lambda_j$ is a Neumann eigenvalue of (wlog) $\Omega\sub{T}$ or is in $\{\tau_n\}_{n=1}^{\infty}$. 

\textbf{Case 1:} $\lambda_j = \mu_k(\Omega\sub{T})$. We use the following to approximate $\lambda_j^w$ and $\varphi_j^w$ (where to simplify notation, we write $\mu_k(\Omega\sub{T})$ as $\mu_k$).
\begin{defn} \label{defn:approx1}
    We define $\mu = \mu(w)$, and functions $U_T(x)$, $U_B(x)$, $U_Q(x)$, defined in $\Omega\sub{T}$, $\Omega\sub{B}$, $Q(w)$ respectively by:
    \begin{align*}
        \mu & = \mu_k + 2\psi_k^T(p\sub{T})a_1w, \quad
        U_Q(x) = \psi_k^T(p\sub{T})\sin(\sqrt{\mu_k}(x_2 + 1))/\sin(2\sqrt{\mu_k}), \\
        U_T(x) & = 2a_1wG^T(x,p\sub{T},\mu), \quad
        U_B(x)  = 2\tilde{a}_1wG^B(x,p\sub{\!B},\mu).
    \end{align*}
    Here $a_1 = (\pa_{x_2}U_Q)(p\sub{T})$ and $\tilde{a}_1 = -(\pa_{x_2}U_Q)(p\sub{\!B})$.
\end{defn}
Note that if $\psi_k^T(p\sub{T})$ is non-zero, then we have  $\tilde{a}_1\neq0$. Additionally, the assumption $\cos(2\sqrt{\mu_k})\neq0$ from Theorem \ref{thm:ends},  guarantees that the coefficient $a_1$ is also  non-zero.

We will show that the quantities in Definition \ref{defn:approx1} approximate $\lambda_j^w$ and $\varphi_j^w$ to the following precision.
\begin{prop} \label{prop:approx1}
Suppose that $\psi_k^T(p\sub{T})\neq0$ and $a_1\neq0$. As $w\to0$, the eigenvalue $\lambda_j^w$ and eigenfunction $\varphi_j^w$ of $\Omega(w)$ satisfy
\begin{align*}
   & \lambda_j^w - \mu  = O(w^2|\ln(w)|), \quad \norm{\varphi_j^w - c_wU_Q}_{H^1(Q(w)\backslash (D_\eta(p\sub{T})\cup D_\eta(p\sub{\!B})))} = O(w^{3/2}|\ln(w)|),  \\
   & \norm{\varphi_j^w - c_wU_T}_{H^1(\Omega\sub{T}\backslash D_\eta(p\sub{T}))}  = O(w^2|\ln(w)|), \quad \norm{\varphi_j^w - c_wU_B}_{H^1(\Omega\sub{B}\backslash D_\eta(p\sub{\!B}))} = O(w^2|\ln(w)|),
\end{align*}
for a constant $c_w = 1+o(1)$. Here $D_\eta(p\sub{T})$, $D_\eta(p\sub{\!B})$ are discs of radius $\eta>0$ centered at $p\sub{T}$, $p\sub{\!B}$, and the implicit constants depend on $\eta$, $j$, the geometry of $\Omega\sub{T}$, $\Omega\sub{B}$, and $\min\{\lambda_j-\lambda_{j-1},\lambda_{j+1}-\lambda_j\}$.
\end{prop}

\textbf{Case 2:} $\lambda_j = \tau_n$. We use the following to approximate $\lambda_j^w$ and $\varphi_j^w$.
\begin{defn} \label{defn:approx2}
    We define $\tau = \tau(w)$, and functions $X_T(x)$, $X_B(x)$, $X_Q(x)$, defined in $\Omega\sub{T}$, $\Omega\sub{B}$, $Q(w)$ respectively by:
    \begin{align*}
        \tau & = \tau_n + n^2\pi w\ln(w), \quad
        X_Q(x) = (2w)^{-1/2}\gamma_n(x_2), \\
        X_T(x) & = b_1(2w)^{1/2}G^T(x,p\sub{T},\tau), \quad
        X_B(x)  = \tilde{b}_1(2w)^{1/2}G^B(x,p\sub{\!B},\tau).
    \end{align*}
    Here $\gamma_n$ is the $L^2([1,-1])$-normalized $n$-th Dirichlet eigenfunction on $[1,-1]$,  $b_1 = \gamma_n'(1)$ and $\tilde{b}_1 = -\gamma_n'(-1)$. For $n=1$, we, without loss of generality, take $\gamma_1(x_2)>0$. 
\end{defn}
Note that $b_1\neq0$ and $\tilde{b}_1\neq0$. In this case, we will show that these quantities approximate $\lambda_j^w$ and $\varphi_j^w$ to the following precision.
\begin{prop} \label{prop:approx2}
As $w\to0$, the eigenvalue $\lambda_j^w$ and eigenfunction $\varphi_j^w$ of $\Omega(w)$ satisfy
\begin{align*}
   & \lambda_j^w - \tau  = O(w), \quad \norm{\varphi_j^w - c_wX_Q}_{H^1(Q(w)\backslash (D_\eta(p\sub{T})\cup D_\eta(p\sub{\!B})))} = O(w|\ln(w)|),  \\
 &   \norm{\varphi_j^w - c_wX_T}_{H^1(\Omega\sub{T}\backslash D_\eta(p\sub{T}))}  = O(w^{3/2}|\ln(w)|), \quad \norm{\varphi_j^w - c_wX_B}_{H^1(\Omega\sub{B}\backslash D_\eta(p\sub{\!B}))} = O(w^{3/2}|\ln(w)|),
\end{align*}
for a constant $c_w = 1+o(1)$, with implicit constants depending on $\eta>0$, $j$, the geometry of $\Omega\sub{T}$, $\Omega\sub{B}$,  and $\min\{\lambda_j-\lambda_{j-1},\lambda_{j+1}-\lambda_j\}$.
\end{prop}

We will prove these propositions in Section \ref{sec:matched}. In fact, in Section \ref{sec:matched}, we will construct asymptotics for $\varphi_j^w$ up to an arbitrarily high power of $w$, with leading order terms agreeing with the functions given in Definitions \ref{defn:approx1} and \ref{defn:approx2}. We will do this by building functions that satisfy the eigenfunction equation up to a large power of $w$, and then convert this to an approximation to $\lambda_j^w$ and $\varphi_j^w$ using the following lemma. 
\begin{lemma}[Lemma 2.6 analogue in \cite{Gad}] \label{lem:approx}
   Let $\lambda$ be a non-repeated element in $\{\lambda_{j}\}_{j=1}^{\infty}$. Suppose that there exists $\lambda(w)$ and a function $u_w\in H^1(\Omega(w))$ satisfying
    \begin{align*}
        (\Delta + \lambda(w))u_w = f_w \text{ in } \Omega(w), \quad \pa_\nu u_w = 0\text{ on }\pa\Omega(w),
    \end{align*}
    with $\lambda(w) = \lambda+o(1)$, $\norm{u_w}_{L^2(\Omega(w))} = 1+ o(1)$ and $\norm{f_w}_{L^2(\Omega(w))} = O(w^N)$ for some $N$. Then, letting $\lambda^w$ be the eigenvalue of $\Omega(w)$ satisfying $\lim_{w\to0}\lambda^w = \lambda$, it has a $L^2(\Omega(w))$-normalized eigenfunction $\varphi^w$, with
    \begin{align*}
        |\lambda(w) - \lambda^w| &= O(w^N) \\
       \norm{\varphi^w - c_wu_w}_{H^1(\Omega(w))} & = O(w^N),
    \end{align*}
    for a constant $c_w = 1+o(1)$.
\end{lemma}
\begin{proof1}{Lemma \ref{lem:approx}}
Let $v_w\in H^1(\Omega(w))$, with $(\Delta + \lambda(w))v_w = g_w\in L^2(\Omega(w))$, and $\pa_\nu v_w=0$ on $\pa\Omega(w)$. We first obtain a $L^2(\Omega(w))$ estimate on $v_w$. Recall $\{\varphi_{k}^{w}\}_{k=1}^{\infty}$ is an orthonormal basis of Neumann eigenfunctions for $\Omega(w)$, with corresponding eigenvalues $\lambda_{k}^{w}$. Then, we can write
    \begin{align*}
        g_w = \sum_{k=1}^{\infty}a_{k}\varphi_{k}^{w}
    \end{align*}
    for coefficients $a_k$, satisfying $\norm{g_w}^2_{L^2(\Omega(w))} = \sum_{k=1}^{\infty}a_k^2$. Then, provided $\lambda(w)$ is not in the Neumann spectrum of $\Omega(w)$,
    \begin{align*}
        v_w = \sum_{k=1}^{\infty}\frac{a_{k}}{\lambda(w)-\lambda_k^{w}}\varphi_{k}^w,
    \end{align*}
    and so
    \begin{align*}
        \norm{v_w}^2_{L^2(\Omega(w))} = \sum_{k=1}^{\infty}\frac{a_k^2}{(\lambda(w)-\lambda_k^{w})^2} \leq \frac{1}{\min_{k}(\lambda(w)-\lambda_{k}^{w})^2} \norm{g_w}^2_{L^2(\Omega(w))}.
    \end{align*}
    Since $\lambda(w) =\lambda + o(1)$, for $w>0$ sufficiently small, this minimum must be attained by $\lambda_k^w = \lambda^w$, and so
    \begin{align} \label{eqn:approx1}
         \norm{v_w}^2_{L^2(\Omega(w))} \leq \frac{1}{(\lambda(w)-\lambda^{w})^2} \norm{g_w}^2_{L^2(\Omega(w))}.
    \end{align}
    Moreover, if $v_w$ is orthogonal to $\varphi^w$, then the corresponding coefficient $a_k$ must be zero, and so in this case, for sufficiently small $w>0$,
    \begin{align} \label{eqn:approx2}
 \norm{v_w}^2_{L^2(\Omega(w))} = \sum_{k=1, \lambda^w_k\neq \lambda^w}^{\infty}\frac{a_k^2}{(\lambda(w)-\lambda_k^{w})^2} \leq \frac{1}{\min_{\lambda^w_k\neq \lambda^w}(\lambda(w)-\lambda_{k}^{w})^2} \norm{g_w}^2_{L^2(\Omega(w))} \leq C\norm{g_w}^2_{L^2(\Omega(w))}.
    \end{align}
    This final inequality holds because the limiting eigenvalue $\lambda$ is isolated: every eigenvalue branch other than $\lambda^w$ converges either to a limiting eigenvalue separated from $\lambda$ or eventually lies above a fixed spectral threshold. Hence, for small $w>0$, we have a uniform lower bound on $|\lambda(w) - \lambda_k^w|$ for all $k$ with $\lambda_k^w\neq\lambda^w$ that is coming from the distance of $\lambda$ from the rest of $\{\lambda_j\}_{j=1}^{\infty}$.

    Either $\lambda^w = \lambda(w)$, or else we apply \eqref{eqn:approx1} to $v_w = u_w$ and $g_w =f_w$ to get
    \begin{align*}
        1+o(1)= \norm{u_w}_{L^2(\Omega(w))} \leq \frac{O(w^N)}{|\lambda(w)-\lambda^w|}.
    \end{align*}
    This ensures that $\lambda^w = \lambda(w) + O(w^N)$, as required.
    
    Choose the sign of $\varphi^w$ so that $(u_w,\varphi^w)_{L^2(\Omega(w))}\geq0$. Now let $\tilde{u}_w = u_w - (u_w,\varphi^w)\varphi^w$. Then, $\tilde{u}_w$ is orthogonal to $\varphi^{w}$, satisfies Neumann boundary conditions, and 
\begin{align*}
    (\Delta + \lambda(w))\tilde{u}_{w} = f_{w}  -(\lambda(w)-\lambda^w)(u_w,\varphi^w)_{L^2(\Omega(w))}\varphi^{w}.
\end{align*}
Applying \eqref{eqn:approx2} with $v_w = \tilde{u}_w$ and $g_w$ the right hand side of the above equation, gives
\begin{align*}
    \norm{\tilde{u}_w}_{L^2(\Omega(w))} = O(w^N).
\end{align*}
This in particular ensures that $(u_w,\varphi^w)_{L^2(\Omega(w))}=1 + o(1)$, and so setting $c_w = (u_w,\varphi^w)_{L^2(\Omega(w))}^{-1}$, gives
\begin{align*}
    \norm{\varphi^w - c_w{u}_w}_{L^2(\Omega(w))} = O(w^N).
\end{align*}

Finally, set $z_w=\varphi^w-c_wu_w$. The preceding estimates give $\|z_w\|_{L^2(\Omega(w))}=O(w^N)$ and
\[
(\Delta + \lambda(w))z_w = h_w \text{ with }\|h_w\|_{L^2(\Omega(w))}=O(w^N).
\]
Multiplying the equation $(\Delta + \lambda(w))z_w = h_w$ by $z_w$, integrating by parts, and using the Neumann boundary condition yields
\[
\|\nabla z_w\|_{L^2(\Omega(w))}^2
\leq \lambda(w)\|z_w\|_{L^2(\Omega(w))}^2
 +\|(\Delta+\lambda(w))z_w\|_{L^2(\Omega(w))}\|z_w\|_{L^2(\Omega(w))}
=O(w^{2N}).
\]
Consequently, $\|z_w\|_{H^1(\Omega(w))}=O(w^N)$.
\end{proof1}

To construct the appropriate function $u_w$ in the two cases and prove Propositions \ref{prop:approx1} and \ref{prop:approx2}, we will use matched asymptotics to match  Ansatze for an approximate eigenfunction in the regions $\Omega\sub{T}$ and $\Omega\sub{B}$ with one in the neck $Q(w)$. We will carry this out in Section \ref{sec:matched}.

\section{A nodal domain count} \label{sec:proof}

In this section, we will prove Theorems \ref{thm:ends} and \ref{thm:neck}. We will do this by using Propositions \ref{prop:approx1} and \ref{prop:approx2} to relate the number of nodal domains of $\varphi_j^w$ in $\Omega(w)$ to the number of nodal domains of $G^T(x,p\sub{T},\mu)$ and $G^B(x,p\sub{\!B},\mu)$ in $\Omega\sub{T}$ and $\Omega\sub{B}$, and to the number of nodal domains of $U_Q(x)$ and $X_Q(x)$ in the neck $Q(w)$. We use similar techniques to those in \cite{beck2026nodal} where the nodal sets of eigenfunctions of symmetric dumbbells were studied. First, the following lemma uses the estimates of Propositions \ref{prop:approx1} and \ref{prop:approx2} to obtain information about the nodal sets of $\varphi_j^w$ in the neck.
\begin{lemma}[Analogue of Lemma 5.4 in \cite{beck2026nodal}] \label{lem:nodal-neck}
    For $\lambda_{j-1}<\lambda_j=\mu_k(\Omega\sub{T})<\lambda_{j+1}$, assume that $\psi_k^T(p\sub{T})\neq0$ and $\cos(2\sqrt{\mu_k(\Omega\sub{T})})\neq0$, and let $N_k$ denote the number of values of $x_2$ in $(-1,1)$ for which $U_Q(x)$ from Definition \ref{defn:approx1} is equal to $0$. Then, $N_k = m-1$, where $\tau_m$ is the first of the eigenvalues $\{\tau_n\}_{n=1}^{\infty}$ larger than $\mu_k(\Omega\sub{T})$. Moreover, there exists a constant $\delta_0>0$, such that for each $0<\delta<\delta_0$, we can find $w_0>0$ so that, if $0<w<w_0$, then the nodal set of $\varphi_j^w$ in $(-w,w)\times(-1+\delta,1-\delta)$ consists of $N_k$ disjoint curves, each intersecting the left and right boundaries of $Q(w)$. Further, $\varphi_j^w(x)$  has the same sign as $U_Q(x)$, for $x= (x_1,1-\delta)$ and $x=(x_1,-1+\delta)$.

    For $\lambda_{j-1}<\lambda_j=\tau_n<\lambda_{j+1}$, we get the same properties in terms of the $n-1$ values of $x_2$ in $(-1,1)$ for which $X_Q(x)$ from Definition \ref{defn:approx2} is equal to $0$. 
\end{lemma}
Before proving this lemma, we write down some further properties of the nodal set of $\varphi_j^w$. To accurately count the number of nodal domains of $\varphi_j^w$, we will need to rule out nodal domains of small area. Away from the neck $Q(w)$, this will be done using the following lemma.
\begin{lemma}[Lemma 5.5 in \cite{beck2026nodal}] \label{lem:small-nodal}
 Given a Neumann eigenfunctions $\varphi_j^w$ of $\Omega(w)$, there exist constants $c^*>0$, $w_0>0$, independent of $w>0$, such that if $0<w<w_0$ and if $D$ is a nodal domain of $\varphi_j^{w}$ with $\norm{\varphi_j^{w}}_{L^2\left(D\cap(\Omega\sub{T}\cup\Omega\sub{B})\right)}\geq \tfrac{1}{2}\norm{\varphi_j^{w}}_{L^2(D)}$, then $\emph{Area}(D)\geq c^*.$
\end{lemma}
We also need to rule out small nodal domains containing $p\sub{T}$ and $p\sub{\!B}$. If $\lambda_j = \mu_k(\Omega\sub{T})$ and  $\psi_k^T(p\sub{T})\neq0$, $\cos(2\sqrt{\mu_k(\Omega\sub{T})})\neq0$, so that in particular the nodal set of $\psi_k^T$ within $\Omega\sub{T}$ is disjoint from a small neighborhood of $p\sub{T}$, then the following lemma ensures that we have the same property for the nodal set of $\varphi_j^w$. The same property also holds in the $\lambda_j = \tau_n$ case.
\begin{lemma}[Analogue of Lemma 5.8 in \cite{beck2026nodal}] \label{lem:nodal-join}
Suppose $\lambda_{j-1}<\lambda_j = \mu_k(\Omega\sub{T})<\lambda_{j+1}$, with $\psi_k^T(p\sub{T})\neq0$ and $\cos(2\sqrt{\mu_k(\Omega\sub{T})})\neq0$. Then, there exist constants $\delta<\delta_0$ (with $\delta_0$ as in Lemma \ref{lem:nodal-neck}) and $w_0>0$ so that for $0<w<w_0$, $\varphi_j^w$ is of one sign in  $\Omega(w)\cap D_{2\delta}(p\sub{T})$ and $\Omega(w)\cap D_{2\delta}(p\sub{\!B})$. 

This same property holds in the case $\lambda_{j-1}<\lambda_j = \tau_n<\lambda_{j+1}$.
\end{lemma}

\begin{remark} \label{rem:nodal-join}
   In the course of the proof of Lemma \ref{lem:nodal-join}, we will see that  the sign of $\varphi_j^w$ in $\Omega(w)\cap D_{2\delta}(p\sub{T})$ and $\Omega(w)\cap D_{2\delta}(p\sub{\!B})$ in the lemma is the same as $U_Q(x_1,1-\delta)$ and $U_Q(x_1,-1+\delta)$ when $\lambda_j=\mu_k(\Omega\sub{T})$, and the same as $X_Q(x_1,1-\delta)$, $X_Q(x_1,-1+\delta)$ when $\lambda_j=\tau_n$. 
\end{remark}
We first use these lemmas together with Propositions \ref{prop:approx1} and \ref{prop:approx2} to prove Theorems \ref{thm:ends} and \ref{thm:neck}, and then we will prove Lemmas \ref{lem:nodal-neck} and \ref{lem:nodal-join}.
\begin{proof1}{Theorem \ref{thm:ends}}
We fix $\delta$ as in Lemma \ref{lem:nodal-join}, and then choose $w_0>0$ so that the estimates from Lemmas \ref{lem:nodal-neck} - \ref{lem:nodal-join} apply for $0<w<w_0$. Since $\lambda_j= \mu_k(\Omega\sub{T})$ and $m_k$ is the smallest index for which $\lambda_j<\tau_{m_k}$, from Lemmas \ref{lem:nodal-neck} and \ref{lem:nodal-join},  for $w>0$ sufficiently small, the eigenfunction $\varphi_j^w$ has nodal set consisting of $m_k-1$ curves in the neck $Q(w)$, and hence exactly $m_k$ nodal domains intersecting $Q(w)$.

To complete the proof of the theorem, we will estimate the number of nodal domains of $\varphi_j^w$ in $\Omega\sub{T}$ and $\Omega\sub{B}$. We start with $\Omega\sub{T}$. The nodal set of $\psi_k^T$, denoted by $\mathcal{N}_{\psi_k^T}$, is disjoint from a neighborhood of $p\sub{T}$, since $\psi_k^T(p\sub{T})>0$. The nodal set partitions $\Omega\sub{T}$ into $M_k$ nodal domains. Moreover, the nodal set is the union of a finite number of smooth curves that may intersect or self-intersect, and at any such intersection point, the tangent lines to these nodal curves meet at equal angles (see, for example, Proposition 2.6 and Remarks 2.7 and 2.8 in \cite{bonnaillie2015nodal}). We now let $\Omega\sub{T}(\eta)$ be the subset of $\Omega\sub{T}$ given by
\begin{align*}
    \{x\in\Omega\sub{T}\,:\, \text{dist}(x,\pa\Omega\sub{T})>\eta \text{ and } \text{dist}(x,\mathcal{N}_{\psi^T_k})>\eta\},
\end{align*}
which for $\eta>0$ sufficiently small has exactly $M_k$ connected components. We will now show that $\varphi_j^w$ and $\psi_k^T$ have the same sign on in $\Omega\sub{T}(\eta)$: There exists a constant $c_1(\eta)>0$ such that $|\psi_k^T|>c_1(\eta)$ in $\Omega\sub{T}(\eta)$. Moreover, we have
\begin{align} \label{eqn:thm-ends1}
    \Delta(\varphi_j^w-\psi_k^T) = -\lambda_j^w(\varphi_j^w-\psi_k^T) - (\lambda_j^w - \mu_k(\Omega\sub{T}))\psi_k^T
\end{align}
which by Propositions \ref{prop:Green1} and \ref{prop:approx1} is $o(1)$ in  $L^2(\Omega\sub{T}(\eta))$. Therefore, applying interior elliptic regularity estimates in $\Omega\sub{T}(\eta)$, for $w$ sufficiently small (depending on $c_1(\eta)$), $\varphi_j^w$ and $\psi_k^T$ have the same sign in $\Omega\sub{T}(\eta)$. In particular, since $\psi_k^T$ takes on positive and negative values in $\Omega\sub{T}(\eta)$, the eigenfunction $\varphi_j^w$ has at least $2$ nodal domains in $\Omega\sub{T}$.  By Lemmas \ref{lem:small-nodal} and \ref{lem:nodal-join}, by taking $\eta>0$ sufficiently small (relative to $c^*$ and $\delta$ in those lemmas), $\varphi_j^w$ has no nodal domain strictly contained in $\Omega\sub{T}\backslash\Omega\sub{T}(\eta)$ for $w>0$ sufficiently small. Therefore, $\varphi_j^w$ has at least 2 and at most $M_k$ nodal domains intersecting $\Omega\sub{T}$.

We use a similar strategy to estimate the number of nodal domains of $\varphi_j^w$  intersecting $\Omega\sub{B}$. Let $k_B$ be the smallest index such that $\mu_k(\Omega\sub{T})<\mu_{k_B}(\Omega\sub{B})$. Then, by Proposition \ref{prop:Green2}, the Green's function $G^B(x,p\sub{\!B},\mu_k(\Omega\sub{T}))$ has at least two and at most $k_B$ nodal domains in $\Omega\sub{B}$. Proposition \ref{prop:Green1} ensures that its nodal set, denoted by $\mathcal{N}_{G^B_{\mu_k}}$ is disjoint from a neighborhood of $p\sub{\!B}$. Since $G^B(x,p_B,\mu_k(\Omega_T))$ satisfies the eigenfunction equation with Neumann boundary conditions away from $p_B$, this ensure that this nodal set consists of a finite number of smooth curves, with the same structure as for that of $\psi_k^T$ above. Therefore, letting $\Omega\sub{B}(\eta) = \{x\in\Omega\sub{B}:\text{dist}(x,\pa\Omega\sub{B})>\eta \text{ and }\text{dist}(x,\mathcal{N}_{G^B_{\mu_k}})>\eta\}$, as for $\Omega\sub{T}$ above, we will show that $\varphi_j^w$ and $2\tilde{a}_1wG^B(x,p\sub{\!B},\mu_k(\Omega\sub{T})$ have the same sign in $\Omega\sub{B}(\eta)$. To show this, we use
\begin{align*}
    & \Delta(\varphi_j^w-2\tilde{a}_1wG^B(x,p\sub{\!B},\mu_k(\Omega\sub{T})) = \\
    & \hspace{2cm}-\mu_k(\Omega\sub{T})(\varphi_j^w-2\tilde{a}_1wG^B(x,p\sub{\!B},\mu_k(\Omega\sub{T})))  - (\lambda_j^w - \mu_k(\Omega\sub{T}))\varphi_j^w,
\end{align*}
which by Propositions \ref{prop:Green1} and \ref{prop:approx1} is $o(w)$ in  $L^2(\Omega\sub{B}(\eta))$. Elliptic regularity estimates again therefore ensures that $\varphi_j^w$ and $2\tilde{a}_1wG^B(x,p\sub{\!B},\mu_k(\Omega\sub{T}))$ have the same sign in $\Omega\sub{B}(\eta)$ for $w>0$ sufficiently small. Again $\varphi_j^w$ has no nodal domain contained strictly in $\Omega\sub{B}\backslash \Omega\sub{B}(\eta)$. Therefore,  for small $w>0$, $\varphi_j^w$ has at least $2$ and at most $k_B$ nodal domains intersecting $\Omega\sub{B}$.

Putting everything together, for $w>0$ sufficiently small, the eigenfunction $\varphi_j^w$ has at least $m_k + 2 + 2-2 = m_k+2$ and at most $m_k+M_k + k_B-2$ nodal domains in $\Omega(w)$. (Here we have subtracted 2 to avoid double counting the nodal domain that intersects both $\Omega\sub{T}$ and $Q(w)$ and the nodal domain that intersects both $\Omega\sub{B}$ and $Q(w)$.) By the way that $m_k$ and $k_B$ have been defined, we have $j = k + (m_k-1) + (k_B-1)$, and so these bounds can be written as
\begin{align*}
    m_k+2 \leq \nu(\varphi_j^w) \leq j-(k-M_k),
\end{align*}
as required.

Finally, if additionally the nodal set of $\psi_k^T$ has no crossings in $\Omega\sub{T}$, then the nodal set $\mathcal{N}_{\psi_k^T}$ consists of $M_k-1$ disjoint curves. Moreover, a small $\eta>0$ can be chosen such that the set $\{x\in\Omega\sub{T}:\text{dist}(x,\mathcal{N}_{\psi_k^T})=\eta\}\cap \Omega\sub{T}$ consists of $M_k-1$ pairs of curves, with $\psi_k^T$ positive and negative on one of each of these pairs of curves, and with any intersection of these curves with $\pa\Omega\sub{T}$ located on the smooth sides of $\pa\Omega\sub{T}$. We can then again use \eqref{eqn:thm-ends1} and elliptic regularity  to ensure that $\varphi_j^w$ and $\psi_k^T$ have the same sign on these pairs of curves, for $w>0$ sufficiently small. Since by Lemma \ref{lem:small-nodal}, for small enough $\eta$, $\varphi_j^w$ cannot have a nodal domain contained entirely between one of these pairs, this ensures that, for small enough $w>0$, the nodal set of $\varphi_j^w$ in $\Omega\sub{T}$ consists of $M_k-1$ disjoint curves, and hence $\varphi_j^w$ has exactly $M_k$ nodal domains intersecting $\Omega\sub{T}$. So, $\varphi_j^w$ has at least $m_k+M_k+2-2= m_k+M_k$ and at most $m_k+M_k+k_B-2$ nodal domains in $\Omega(w)$. In particular, if $\mu_k(\Omega\sub{T})<\mu_2(\Omega\sub{B})$, then $k_B=2$, and $\varphi_j^w$ has exactly $m_k+M_k$ nodal domains, completing the proof of the theorem. 
\end{proof1}
\begin{proof1}{Theorem \ref{thm:neck}}
The proof of this theorem follows the same strategy as for Theorem \ref{thm:ends} above. Firstly, from Lemmas \ref{lem:nodal-neck} and \ref{lem:nodal-join}, $\varphi_j^w$ has exactly $n$ nodal domains intersecting $Q(w)$ for small $w$, with nodal set in $Q(w)$ consisting of exactly $n-1$ curves.

To estimate the number of nodal domains in $\Omega\sub{T}$ and $\Omega\sub{B}$, we use the same proof strategy as for $\Omega\sub{B}$ in the proof of Theorem \ref{thm:ends}: Proposition \ref{prop:Green1} ensures that the nodal set of $G^T(x,p\sub{T},\tau_n)$, denoted by $\mathcal{N}_{G^T_{\tau_n}}$, is disjoint from a neighborhood of $p\sub{T}$. We now let $\Omega\sub{T}(\eta) = \{x\in\Omega\sub{T}:\text{dist}(x,\pa\Omega\sub{T})>\eta \text{ and }\text{dist}(x,\mathcal{N}_{G^T_{\tau_n}})>\eta\}$. Following the proof above when estimating nodal domains in $\Omega\sub{B}$, this time using  Proposition \ref{prop:approx2} in place of Proposition \ref{prop:approx1}, ensures that, for small $w>0$, $\varphi_j^w$ and $b_1(2w)^{1/2}G^T(x,p\sub{T},\tau_n)$ have the same sign in $\Omega\sub{T}(\eta)$. Moreover, it has no nodal domains entirely contained in $\Omega\sub{T}\backslash \Omega\sub{T}(\eta)$. Therefore, $\varphi_j^w$ has at least 2 nodal domains intersecting $\Omega\sub{T}$, and exactly $2$ nodal domains if $\tau_n<\mu_2(\Omega\sub{T})$. Analogously, for small $w>0$, $\varphi_j^w$ has at least 2 nodal domains intersecting $\Omega\sub{B}$, and exactly $2$ nodal domains if $\tau_n<\mu_2(\Omega\sub{B})$. Putting everything together, $\varphi_j^w$ has at least $n+2+2-2 = n+2$ nodal domains in $\Omega(w)$, for small $w>0$. If $\tau_n<\min\{\mu_2(\Omega\sub{T}),\mu_2(\Omega\sub{B})\}$, then $j = n+2$, and so $\varphi_j^w$ must have exactly $n+2$ nodal domains and hence is Courant sharp, with the required nodal set structure given in the statement of the theorem.
\end{proof1}

We now prove the two remaining lemmas.
\begin{proof1}{Lemma \ref{lem:nodal-neck}}
\textbf{Case 1.} $\lambda_j = \mu_k(\Omega\sub{T})$. We first establish the formula for the number of zeros of $U_Q(x_2)$ in $(-1,1)$, denoted by $N_k$. Recall that $\gamma_n$ is the $n$-th Dirichlet eigenfunction of $[1,-1]$. If $m=1$, then $0<\mu_k(\Omega\sub{T})<\tau_1=\pi^2/4$, so $U_Q$ has no zeros in $(-1,1)$, and therefore $N_k=0=m-1$; hence we may assume that $m\geq2$.  Since $\mu_k(\Omega\sub{T})<\tau_m$, and $U_Q(-1)=0$, applying the Sturm comparison theorem to $U_Q$ and $\gamma_m$, we find that $\gamma_m$ has a zero in each of the $N_k$ intervals between consecutive zeros of $U_Q$. Since $\gamma_m$ has exactly $m-1$ zeros, this means that $N_k\leq m-1$. Similarly, since $\tau_{m-1}<\mu_k(\Omega\sub{T})$ and $\gamma_{m-1}(\pm 1)=0$, we find that $U_Q$ has a zero in each of the $m-1$ intervals between consecutive zeros of $\gamma_{m-1}$, and so $N_k \geq m-1$. Combining the two estimates gives $N_k=m-1$.

We now study the nodal set of $\varphi_j^w$ in $Q(w)$. Following exactly the proof of Lemma 5.5 in \cite{beck2026nodal}, we get the following: Let $h(x_2)$ be a solution of a Sturm-Liouville problem in $[-1,1]$, with $N$ zeros in $(-1,1)$, and let $g(x;w)$ be a family of functions satisfying
\begin{align*}
    \norm{g(\cdot;w) - h}^2_{H^1(Q(w))} = o(w)\norm{h}^2_{L^2([-1,1])}.
\end{align*}
Then, for each $\delta>0$ sufficiently small, there exists $w_0>0$ such that for $0<w<w_0$, the function $g(\cdot;w)$ has exactly $N$ nodal curves in $(-w,w)\times[-1+\delta,1-\delta]$, disjoint and touching the left and right boundaries of the neck. Moreover, $g(x_1,x_2;w)$ has the same sign as $h(x_2)$ for $x_2 = -1+\delta$, $1-\delta$.

Applying this to $h = U_Q$ and $g = \varphi_j^w$, using Proposition \ref{prop:approx1}, gives the desired nodal structure for $\varphi_j^w$.

\textbf{Case 2.} $\lambda_j = \tau_n$. Since $X_Q(x)$ is a $n$-th Dirichlet eigenfunction on $[-1,1]$, it has $n-1$ zeros in $(-1,1)$. The desired nodal structure for $\varphi_j^w$, then follows as above, by setting $h=X_Q$, $g=\varphi_j^w$.
\end{proof1}
\begin{proof1}{Lemma \ref{lem:nodal-join}}
\textbf{Case 1.} $\lambda_j = \mu_k(\Omega\sub{T})$. We start by establishing the sign of $\varphi_j^w$ in a neighborhood of $p\sub{T}$. Assuming, without loss of generality, that $\psi_k^T(p\sub{T})>0$, using Proposition \ref{prop:approx1} and elliptic regularity, there exists $\delta_0>0$ such that for each $0<\delta<\delta_0$, provided $w>0$ is sufficiently small, the eigenfunction $\varphi_j^w$ is positive on the semi-circle $\pa D_{2\delta}(p\sub{T})\cap \Omega\sub{T}$. Reducing $\delta_0$ if necessary, by Lemma \ref{lem:nodal-neck}, we also have that $\varphi_j^w$ is positive on $\pa D_{2\delta}(p\sub{T})\cap Q(w)$.

Therefore, we need to rule out a nodal domain $D$ of $\varphi_j^w$ contained entirely within $\Omega(w)\cap D_{2\delta}(p\sub{T})$. The proof of this below holds for any eigenfunction $\varphi_j^w$, and does not use that $\lambda_j =\mu_k(\Omega\sub{T})$. By Lemma \ref{lem:small-nodal}, and reducing $\delta_0$ so that $D_{2\delta_0}(p\sub{T})$ has area smaller than the constant $c^*$ in that lemma, such a nodal domain must satisfy
\begin{align} \label{eqn:nodal-join1}
    \norm{\varphi_j^w}_{L^2(D\cap Q(w))} >\tfrac{1}{2}\norm{\varphi_j^w}_{L^2(D)}.
\end{align}
Suppose such a nodal domain exists, and let $v$ be a multiple of $\varphi_j^w$, with $\norm{v}_{L^2(D)} = 1$. Then, by \eqref{eqn:nodal-join1}, there exists $y^*\in[1-2\delta,1]$ such that 
\begin{align} \label{eqn:nodal-join2}
    \int_{V(y^*)} v(x_1,y^*)^2\,dx_1 >\tfrac{1}{8}\delta^{-1}.
\end{align} 
Here $V(y^*)$ is the cross-section of $D$ at $x_2 = y^*$. 
Let $\tilde v$ be equal to $v$ on $D$ and to $0$ on $Q(w)\setminus D$. Since $D$ is a nodal domain of $v$, we have $\tilde v\in H^1(Q(w))$. Hence, for almost every $x_1\in(-w,w)$, the function $\tilde v(x_1,\cdot)$ belongs to $H^1(1-2\delta,1)$.
Moreover,
\begin{align} \label{eqn:nodal-join3}
    \int_{V(y^*)}\int_{1-2\delta}^{y^*}|\pa_{x_2} \tilde{v}(x_1,x_2)|^2\,dx_2\,dx_1 \leq \lambda_j^w.
\end{align}
This final inequality comes from $v$ being an eigenfunction of eigenvalue $\lambda_j^w$ and $\norm{v}_{L^2(D)} = 1$. Writing
\begin{align*}
    \tilde{v}(x_1,x_2) = \tilde{v}(x_1,y^*) - \int_{x_2}^{y^*}\pa_{t}\tilde{v}(x_1,t)\,dt
\end{align*}
and then using Young's inequality gives
\begin{align*}
     \tilde{v}(x_1,x_2)^2 \geq \tfrac{1}{2}\tilde{v}(x_1,y^*)^2 - \left(\int_{x_2}^{y^*}\pa_{t}\tilde{v}(x_1,t)\,dt\right)^2.
\end{align*}
Integrating this with respect to $x_1\in V(y^*)$, and using \eqref{eqn:nodal-join2}, \eqref{eqn:nodal-join3}, and Cauchy-Schwarz on the second term on the right hand side, we therefore obtain
\begin{align*}
    \int_{V(y^*)} \tilde{v}(x_1,x_2)^2\,dx_1 > \tfrac{1}{16} \delta^{-1} - |y^*-x_2|\int_{V(y^*)}\int_{x_2}^{y^*}|\pa_t\tilde{v}(x_1,t)|^2\,dt\,dx_1 \geq \tfrac{1}{16}\delta^{-1}  - |y^*-x_2|\lambda_j^w.
\end{align*}

As $\lim_{w\to0}\lambda_j^w=\lambda_j$, the eigenvalues $\lambda_j^w$ are uniformly bounded for sufficiently small $w$. Taking $x_2=1-2\delta$ and reducing $\delta_0$ if necessary, the right-hand side of the preceding inequality is positive for every $0<\delta<\delta_0$. Hence the support of $\tilde v$ meets the line $x_2=1-2\delta$ in a set of positive measure and therefore contains a point outside $D_{2\delta}(p\sub{T})$.
This gives a contradiction as the support of $\tilde{v}$ is $D\subset D_{2\delta}(p\sub{T})$. Therefore, $\varphi_j^w$ is positive on $\pa D_{2\delta}(p\sub{T})\cap \Omega(w)$ and has no nodal domain contained within $D_{2\delta}(p\sub{T})$, and hence must be positive in the whole of $\Omega(w)\cap D_{2\delta}(p\sub{T})$. 

The sign of $\varphi_j^w$ in a neighborhood of $p\sub{\!B}$ can be established in the same way once we have shown that $\varphi_j^w$ is of one sign on $\pa D_{2\delta}(p\sub{\!B})\cap \Omega(w)$. For small enough $w$, from Lemma \ref{lem:nodal-neck}, the eigenfunction $\varphi_j^w$ has the same sign as $U_Q(x_1,-1+2\delta)$ on $\pa D_{2\delta}(p\sub{\!B})\cap Q(w)$, and this sign is the opposite to $\tilde{a}_1 = -(\pa_{x_2}U_Q)(p\sub{\!B})$. From Proposition \ref{prop:Green1}, we know that $U_B(x)$ has the opposite sign to $\tilde{a}_1$ in a small neighborhood around $p\sub{\!B}$ in $\Omega\sub{B}$. Therefore,  using Proposition  \ref{prop:approx1} and elliptic regularity, the eigenfunction $\varphi_j^w$ is of one sign on the $\pa D_{2\delta}(p\sub{\!B})\cap \Omega(w)$ for sufficiently small $w$. The same proof as above rules out a nodal domain contained entirely within $\Omega(w)\cap D_{2\delta}(p\sub{\!B})$, finishing the proof of this case.

\textbf{Case 2.} $\lambda_j=\tau_n$. The proof of this case works in the same way as how the sign of $\varphi_j^w$ was established in $D_{2\delta}(p\sub{\!B})\cap\Omega(w)$ in Case 1: This time we can use Lemma \ref{lem:nodal-neck} and Propositions \ref{prop:Green1} and \ref{prop:approx2} to establish that $\varphi_j^w$ is of one sign on $\pa D_{2\delta}(p\sub{T})\cap \Omega(w)$ and on $\pa D_{2\delta}(p\sub{\!B})\cap \Omega(w)$. The proof in Case 1 then ensures there is no nodal domain contained within $\Omega(w)\cap D_{2\delta}(p\sub{T})$ or $\Omega(w)\cap D_{2\delta}(p\sub{\!B})$, finishing the proof in this case.
\end{proof1}

We conclude the section with the short deductions of the three corollaries stated in the Introduction.

\begin{proof1}{Corollary \ref{cor:lower}}
If $\lambda_j=\mu_k(\Omega\sub{T})$, then we can apply Theorem \ref{thm:ends}. In this case, $\lambda_j<\tau_{m_k} = \pi^2m_k^2/4$, and so $m_k >2\pi^{-1}\sqrt{\lambda_j}$. Using the lower bound on $\nu(\varphi_j^w)$ from Theorem \ref{thm:ends} therefore gives
\begin{align*}
    \nu(\varphi_j^w)\geq m_k+2 > 2\pi^{-1}\sqrt{\lambda_j} + 2.
\end{align*}
Since $\nu(\varphi_j^w)$ is an integer, this bound implies the bound in the corollary.

If instead $\lambda_j = \tau_n = \pi^2n^2/4$, then applying Theorem \ref{thm:neck} gives
\begin{align*}
    \nu(\varphi_j^w)\geq n+2 = 2\pi^{-1}\sqrt{\lambda_j} + 2,
\end{align*}
again implying the required lower bound.
\end{proof1}

\begin{proof1}{Corollary \ref{cor:3rd}}
If $\lambda_3 = \mu_2(\Omega\sub{T})$, then the assumptions of the corollary ensure that Theorem \ref{thm:ends} applies, with $m_2 = 1$. Therefore, $\varphi_3^w$ has at least three (and hence by Courant exactly three) nodal domains for small $w>0$. Moreover, it has exactly one nodal domain intersecting $Q(w)$ (ensuring the nodal set of $\varphi_3^w$ is disjoint from $Q(w)$) and at least two nodal domains intersect each of $\Omega\sub{T}$ and $\Omega\sub{B}$ (which, together with the total count of three and the single nodal domain meeting $Q(w)$, yields one nodal curve in each end domain).

If instead $\lambda_3=\tau_1$, then Theorem \ref{thm:neck} (with $n=1$) immediately implies the result of the corollary.
\end{proof1}
\begin{proof1}{Corollary \ref{cor:4th}}
If $\psi_k^T$ is not Courant sharp in $\Omega\sub{T}$, then $k-M_k>0$, and so applying Theorem \ref{thm:ends}, we have
\begin{align*}
    \nu(\varphi_j^w) \leq j-(k-M_k)<j,
\end{align*}
for small $w>0$, ensuring that $\varphi_j^w$ is not Courant sharp in $\Omega(w)$.
\end{proof1}

\begin{remark} \label{rem:4th}
In the symmetric case $\Omega\sub{T}=\Omega\sub{B}$, if $\lambda_3=\lambda_4=\mu_2(\Omega\sub{T})=\mu_2(\Omega\sub{B})<\tau_1$ and $\psi_2^T(p\sub{T})\neq0$, Theorem 1.5 of \cite{beck2026nodal} gives three nodal domains for $\varphi_3^w$ and four for $\varphi_4^w$ for all sufficiently small $w>0$.
\end{remark}

\section{Matched asymptotics} \label{sec:matched}

In this section, we will use matched asymptotics to prove Propositions \ref{prop:approx1} and \ref{prop:approx2}. In both cases, we will construct $\lambda(w) = \lambda_j + o(1)$, and a function $u_w$ with $\norm{u_w}_{L^2(\Omega(w))}=1+o(1)$, satisfying Neumann boundary conditions, such that
\begin{align*}
    (\Delta + \lambda(w))u_w = f_w \text{ in }\Omega(w),
\end{align*}
with $\norm{f_w}_{L^2(\Omega(w))} = O(w^N)$ for some large $N$. The strategy of proof is the same as the three dimensional case studied in \cite{Gad}, where Ansatze for the  eigenfunction in $\Omega\sub{T}$, $\Omega\sub{B}$, and $Q(w)$ are written down, and then determined by using matched asymptotics in a neighborhood of where the necks are joined. However, because our two dimensional Green's function has a logarithmic singularity (rather than a power singularity as in three dimensions), the Ansatze will now be constructed as a series involving powers of both $w$ and $\ln(w)$ (rather than  just powers of $w$ as in three dimensions). The first terms in these Ansatze will agree with the quantities given in Definitions \ref{defn:approx1} and \ref{defn:approx2}. We will then apply Lemma \ref{lem:approx} to show that $\lambda_j^w$ and $\varphi_j^w$ are well-approximated by $\lambda(w)$ and a multiple of $u_w$, and use this to prove the propositions.

\textbf{Case 1. $\lambda_j=\mu_k(\Omega\sub{T})$ and the proof of Proposition \ref{prop:approx1}.}  To simplify notation, we will write $$\lambda = \mu_k(\Omega\sub{T}),$$ and denote a $L^2(\Omega\sub{T})$-normalized Neumann eigenfunction of eigenvalue $\lambda$ in $\Omega\sub{T}$ by $\psi$. 

\subsection{Eigenvalue and Eigenfunction Ansatze:}

We first define Ansatze for the eigenvalue, $\lambda(w)$, and the eigenfunction in the neck, which we will denote by $V_Q(x_2)$. These are written as a series in terms of powers of $w$ and $\ln(w)$ (due to the logarithmic singularity in the Green's function), and with leading order terms consistent with the expected leading order behavior of the eigenvalue and eigenfunction. That is, the first term in $\lambda(w)$ is given by $\lambda$, and the first term in $V_Q(x_2)$ is the trigonometric function $U_Q(x)$ given in Definition \ref{defn:approx1}.

\begin{defn} \label{defn:evalue-ansatz}
    Fixing a large integer $M$, we write
    \begin{align*}
         \lambda(w)  = \lambda_{0,0} + \sum_{k=1}^{M}\sum_{j=0}^{k-1}w^k(\ln(w))^j\lambda_{k,j} , \qquad
    V_Q(x_2)  = \sum_{k=0}^{M}\sum_{j=0}^{k}w^k(\ln(w))^jv_{k,j}(x_2),
    \end{align*}
    for $\lambda_{k,j}\in\R$ and $v_{k,j}\in C^{\infty}([-1,1])$. We set $\lambda_{0,0}=\lambda$ and
    \begin{align*}
        v_{0,0}(x_2) = \psi(p\sub{T})\sin(\sqrt{\lambda}(x_2+1))/\sin(2\sqrt{\lambda}).
    \end{align*}
     We call the function $v_{k,j}(x_2)$ admissible if the coefficient of $w^k(\ln(w))^j$ in $V_Q''(x_2)+\lambda(w)V_Q(x_2)$ is equal to $0$, and call $V_Q(x_2)$ admissible if $v_{k,j}(x_2)$ is admissible for all $0\leq k \leq M$, $0\leq j \leq k$.
    \end{defn}
\begin{remark} \label{rem:evalue-ansatz}
   For the function $v_{k,j}(x_2)$ to be admissible, it will need to satisfy 
    \begin{align*}
        v_{k,j}''(x_2)+\lambda v_{k,j}(x_2) = F_{k,j}(x_2)
    \end{align*}
    for a function $F_{k,j}(x_2)$ depending on $v_{m,\ell}$, $\lambda_{m,\ell}$ for $m\leq k-1$, and $\lambda_{k,j}$. In particular, an admissible $v_{k,j}(x_2)$ is uniquely determined by these quantities, and its boundary values at $x_2=\pm1$.

    Moreover, an admissible $V_Q(x_2)$ satisfies $(\Delta +\lambda(w))V_Q(x_2) = o(w^{M})$, where here and throughout, the implicit constant depends on the $\lambda_{k,j}$ and $v_{k,j}(x_2)$.
\end{remark}
Next we write down Ansatze for the eigenfunction in $\Omega\sub{T}$ and $\Omega\sub{B}$.
\begin{defn} \label{defn:efn-ansatz}
We write, for $x\in\Omega\sub{T}\backslash\{p\sub{T}\}$,
\begin{align*}
    V_T(x) =\,& \psi(p\sub{T})^{-1}(\lambda(w)-\lambda)G^T(x,p\sub{T},\lambda(w)) +   \sum_{i=2}^{101M}\sum_{j=0}^{i-2}w^{i}(\ln(w))^ja^T_{i,j}G^T(x,p\sub{T},\lambda(w)) \\
    & +   \sum_{m=1}^{100M}\sum_{k=1}^{M}\sum_{\ell=0}^{k-1}w^{k+m}(\ln(w))^{\ell}a^T_{k,\ell,m}\pa_{y_1}^mG^T(x,p\sub{T},\lambda(w)) .
\end{align*}
In the triple sum in $V_T(x)$, the functions of $w$ that appear are of the form $w^{i}(\ln(w))^j$, with $2\leq i \leq 101M$, $0\leq j \leq i-2$. We will uniquely define the coefficients $a_{i,j}^T$ in terms of $a^T_{k,j,m}\in\R$ with $k+m=i$ (and so in particular $k,m \leq i-1$), such that
\begin{align*}
    \sum_{i=2}^{101M}\sum_{j=0}^{i-2}w^{i}(\ln(w))^ja^T_{i,j}\psi(p\sub{T})   + \sum_{m=1}^{100M}\sum_{k=1}^{M}\sum_{\ell=0}^{k-1}w^{k+m}(\ln(w))^{\ell}a^T_{k,\ell,m}\pa_{y_1}^m\psi(p\sub{T}) = 0.
\end{align*}
Under this definition of $a_{i,j}^T$ for all $2\leq i\leq 101M$, $0 \leq j \leq i-2$, we will then call the function $V_T(x)$ admissible.

We also write, for $x\in \Omega\sub{B}\backslash\{p\sub{\!B}\}$,
\begin{align*}
     V_B(x) = \sum_{m=0}^{100M}\sum_{k=1}^{M}\sum_{\ell=0}^{k-1}w^{k+m}(\ln(w))^\ell a^B_{k,\ell,m}\pa_{y_1}^mG^B(x,p\sub{\!B},\lambda(w)),
\end{align*}
with $V_B(x)$ admissible for any $a_{k,\ell,m}^B\in\R$.
\end{defn}
\begin{remark} \label{rem:efn-ansatz}
By construction, we have $(\Delta+\lambda(w))V_T=0$ in $\Omega\sub{T}$, with Neumann boundary conditions on $\pa\Omega\sub{T}\backslash\{p\sub{T}\}$, and analogously for $V_B$.
\end{remark}
Provided the leading coefficients are choosing appropriately, via the following lemma, we can reduce the proof of Proposition \ref{prop:approx1} to showing that $\lambda_j^w$ is sufficiently well approximated by $\lambda(w)$, and $\varphi_j^w$ is sufficiently well approximated by $V_Q$, $V_T$, and $V_B$.
\begin{lemma} \label{lem:approx1}
    Suppose that $V_Q$, $V_T$, $V_B$ are admissible with the coefficient $\lambda_{1,0}$ in $\lambda(w)$ equal to $2\psi(p\sub{T})a_1$, and the coefficient $a_{1,0,0}^B$ in $V_B(x)$ equal to $2\tilde{a}_1$. Here $a_1$ and $\tilde{a}_1$ are as in Definition \ref{defn:approx1}. Then,
    \begin{align*}
        & |\lambda(w) - \mu|  = O(w^2|\ln(w)|), \quad \norm{V_Q - U_Q}_{H^1(Q(w)\backslash (D_\eta(p\sub{T})\cup D_\eta(p\sub{\!B})))} = O(w^{3/2}|\ln(w)|),  \\
   & \norm{V_T - U_T}_{H^1(\Omega\sub{T}\backslash D_\eta(p\sub{T}))}  + \norm{V_B - U_B}_{H^1(\Omega\sub{B}\backslash D_\eta(p\sub{\!B}))} = O(w^2|\ln(w)|),
    \end{align*}
   where the implicit constants depend on $\eta$.
\end{lemma}
\begin{proof1}{Lemma \ref{lem:approx1}}
The bound on $\lambda(w)-\mu$ follows immediately from Definitions \ref{defn:approx1} and \ref{defn:evalue-ansatz}. Since the function $v_{0,0}(x_2)$ is equal to $U_Q(x)$, we obtain a pointwise $O(w|\ln(w)|)$ bound on $V_Q(x_2)-U_Q(x)$ and its derivatives. Combining this with the area of the neck $Q(w)$ being $4w$ gives the desired bounds for $V_Q-U_Q$.
 
The definition of the coefficients $a_{i,j}^T$ ensures that the total contribution to the double and triple sum in $V_T(x)$ from the $\frac{\psi(x)\psi(p\sub{T})}{\lambda(w)-\mu}$ term in $G^T(x,p\sub{T},\lambda(w))$ coming from \eqref{eqn:Green1-a} in Proposition \ref{prop:Green1} vanishes. Therefore, the remaining terms in $G^T(x,p\sub{T},\lambda(w))$ from \eqref{eqn:Green1-a} guarantees that all but the first term in $V_T(x)$ are $O(w^2)$ in $H^1(\Omega\sub{T}\backslash D_\eta(p\sub{T}))$. Moreover, we can write
\begin{align*}
&(\lambda(w)-\lambda)G^T(x,p\sub{T},\lambda(w)) - \psi(p\sub{T})U_T(x) \\
&= (\lambda(w)-\mu) \tilde{G}^T(x,p\sub{T},\lambda(w)) + (\mu - \lambda)(\tilde{G}^T(x,p\sub{T},\lambda(w)) - \tilde{G}^T(x,p\sub{T},\mu)),
\end{align*}
where $\tilde{G}^T(x,p\sub{T},\cdot)$ is the Green's function with the contribution from the first term in \eqref{eqn:Green1-a} removed. Using the bounds on $\lambda(w)-\mu$, and the smoothness of $\tilde{G}^T(x,p\sub{T},\cdot)$ as a function of its final variable in $\Omega\sub{T}\backslash D_\eta(p\sub{T})$, we therefore have the required bound on $V_T-U_T$.

All but the $m=0$, $k=1$, $\ell=0$ term in $V_B(x)$ are $O(w^2|\ln(w)|)$ in $H^1(\Omega\sub{B}\backslash D_\eta(p\sub{\!B}))$. Since the coefficient, $a_{1,0,0}^B$ of $G^B(x,p\sub{\!B},\lambda(w))$ in $V_B(x)$ is chosen to agree with the coefficient of $G^B(x,p\sub{\!B},\mu)$ in $U_B(x)$, we can therefore again use the smoothness of $\tilde{G}^B(x,p\sub{\!B},\cdot)$ to obtain the necessary bound on $V_B-U_B$.
\end{proof1}

\subsection{Matching functions in neighborhoods of $p\sub{T}$ and $p\sub{\!B}$}

In order to use the above functions to build an approximation $u_w$ to the actual eigenfunction $\varphi_j^w$, we will carry out matched asymptotics of $V_Q(x)$ with $V_T(x)$ and $V_B(x)$ in neighborhoods of $p\sub{T}$ and $p\sub{\!B}$. This will done using functions of $\xi$ defined in 
\begin{align*}
    \mathcal{V} = \{\xi = (\xi_1,\xi_2)\in\R^2:\xi_2>0\}\cup (-1,1)\times(-\infty,0].
\end{align*}
We write this union as $\mathcal{V} = \mathcal{V}_+\cup\mathcal{V}_-$, with $\mathcal{V}_+$ the half-plane and $\mathcal{V}_-$ a semi-infinite strip of width 2. Note that, since $\pa\Omega\sub{T}$, $\pa\Omega\sub{B}$ are flat near $p\sub{T}$, $p\sub{\!B}$, the set $\mathcal{V}$ can be viewed as a blow-up of $\Omega(w)$ centered at $p\sub{T}$ or $p\sub{\!B}$. A key step in the proof of the matched asymptotics is to construct functions $P(\xi)$ satisfying $\Delta P = F$ in $\mathcal{V}$, $\pa_\nu P = 0$ on $\pa\mathcal{V}$, for given $F$, and given asymptotics $P_+(\xi)$ for large $|\xi|$ in $\mathcal{V}_+$ and $P_-(\xi_2)$ for large $|\xi|$ in $\mathcal{V}_-$. In the below and throughout, $\chi$ is a smooth cut-off function, with $\chi(t)$ equal to 1 for $t>2$ and equal to 0 for $t<1$. Where convenient, we will also view $\xi$ as an element in $\mathbb{C}$ with Re$(\xi)=\xi_1$, Im$(\xi)=\xi_2$.

\begin{prop}[Two dimensional analogue of Lemma 4.3 in \cite{Gad}] \label{prop:match}
     Let $P_+(\xi_1,\xi_2)$ be a smooth function, satisfying $\pa_{\xi_2}P_+ = 0$ on $\{\xi_2=0\}$, and let $P_{-}(\xi_2)$ be a smooth function. Setting $G(\xi) := \Delta (\chi(|\xi|)P_+(\xi) + \chi(-\xi_2)P_-(\xi_2))$, let $F(\xi)$ be any smooth function such that \begin{align*}
    |F(\xi)-G(\xi)|& \leq C(1+|\xi|)^{-N} \text{ in }\mathcal{V}_+ \text{ for some integer } N\geq5, \\
    e^{-\mu\xi_2}|\pa_{\xi_1}^p(F(\xi)-G(\xi))|& \leq C_\mu \text{ in }\mathcal{V}_- \text{ for all } \mu<\pi/2,
\end{align*}
for constants $C$ and $C_{\mu}$, and $p=0,1$. Then, setting
\begin{align*}
    \int_{\mathcal{V}}F(\xi)-G(\xi) \,d\xi= c_F,
\end{align*} 
there exists a unique function $P$, which satisfies $\Delta P = F$ in $\mathcal{V}$, with Neumann boundary conditions on $\pa\mathcal{V}$, is in $H^1$ for any compact subset of $\mathcal{V}$, and satisfies, 
\begin{align*}
   & \left|P(\xi) -\chi(|\xi|)P_+(\xi) - c_F\pi^{-1}\chi(|\xi|)\ln(|\xi|) - \sum_{m=1}^{N-4}C_{F,m}\text{Re}(1/\xi^m)\right| \leq C'(1+|\xi|)^{-N+3} \text{ in } \mathcal{V}_+,  \\
       & |e^{-\mu\xi_2}\chi(-\xi_2)(P(\xi)-(P_-(\xi_2)+c_P))|\leq C_\mu' \text{ in } \mathcal{V}_- \text{ for all }\mu<\pi/2.
    \end{align*}
    The same estimates hold for the gradient of the quantities on the left hand sides. Here $C_{F,m}$, $C'$, $C_\mu'$ and $c_P$ are constants.
\end{prop}
    
Two special cases of this proposition are the following:
\begin{enumerate}
    \item[1)] For $P_+(\xi) = \chi(|\xi|)\xi_1$, $P_-(\xi_2) = 0$, and $F(\xi) = 0$. In this case, $G(\xi)$ is an odd function of $\xi_1$, and so  we have $c_F = 0$. We denote the resulting function $P$ from Proposition \ref{prop:match} by $Y_1(\xi)$.
    \item[2)] For $P_+(\xi) = 0$, $P_-(\xi_2) = \chi(-\xi_2)\xi_2$, and $F(\xi) = 0$. In this case, using Green's identity on the intersection of $\mathcal{V}$ with a square of side length $2R$ centered at the origin gives
    \begin{align*}
        \int_\mathcal{V} -G(\xi)\,d\xi = \lim_{R\to\infty}\int_{-1}^{1}\pa_{\xi_2}P_-(\xi_2)|_{\xi_2=-R}\,d\xi_1 = 2,
    \end{align*}
    and so we have $c_F = 2$. We denote the resulting function $P$ from Proposition \ref{prop:match} by $Y_2(\xi)$.

\end{enumerate}
\begin{proof1}{Proposition \ref{prop:match}}
 We adapt the proof of Lemma 4.3 in \cite{Gad} to the two dimensional case, now taking into account the logarithmic singularity of the fundamental solution of the Laplacian in two dimensions. First, we let $P_1(\xi) = \chi(|\xi|)P_+(\xi) + \chi(-\xi_2)P_-(\xi_2) + c_F\pi^{-1}\ln|\xi|\chi(|\xi|)$. Then, $P_1$ satisfies Neumann boundary conditions on $\pa\mathcal{V}$ and defining $F_1:=\Delta P_1$, we have that $F_1-G$ is compactly supported. A direct calculation yields
    \begin{align*}
       & \int_{\mathcal{V}}(F_1 -F) = \int_{\mathcal{V}}(G-F)+\int_{\mathcal{V}}(F_1-G)  \\
        &= -c_F + \lim_{R\to\infty}\int_{\mathcal{V}\cap D_R}c_F\pi^{-1}\Delta(\ln|\xi|\chi(|\xi|)) = -c_F + c_F\pi^{-1}\lim_{R\to\infty}\int_{\pa(\mathcal{V}\cap D_R)}\pa_{\nu}(\ln|\xi|\chi(|\xi|)) = 0.
    \end{align*}
Let $d_2(\xi)$ be a smooth, and strictly positive function in $\mathcal{V}$, with $d_2(\xi) =- \xi_2$ for $\xi_2<-4$ and $d_2(\xi) = |\xi|\ln|\xi|$ for $\xi_2>0$, $|\xi|>4$. Then, we define $\mathcal{H}(\mathcal{V})$ to be the completion of $C_0^{\infty}(\bar{\mathcal{V}})$ in the norm
\begin{align*}
    \norm{u}_{\mathcal{H}(\mathcal{V})} = \norm{\nabla u}_{L^2(\mathcal{V})} + \norm{d_2^{-1}u}_{L^2(\mathcal{V})}.
\end{align*}
The function $-F_1+F$ is equal to $F-G$ plus the compactly supported function $-F_1+G$, and so by the assumptions on $F-G$, it satisfies $d_2(-F_1+F)\in L^2(\mathcal{V})$. It also integrates to zero. The proof of Lemma 4.1 in \cite{nazarov99} ensures the existence of a solution $\tilde{P}$ in $\mathcal{H}(\mathcal{V})$ to
\begin{align*}
    \Delta \tilde{P} = -F_1+F \text{ in }\mathcal{V}, \qquad \pa_{\nu}\tilde{P} = 0 \text{ on }\pa\mathcal{V},
\end{align*}
and that this solution is unique up to a constant. 

We now find the asymptotic behavior of $\tilde{P}$. First, $\tilde{P}_{-}(\xi):=\chi(-\xi_2)\tilde{P}(\xi)$ satisfies
\begin{align*}
    \Delta \tilde{P}_{-} = g_{-} \text{ in }\mathcal{V}_-, \qquad \tilde{P}_{-} = 0 \text{ on }(-1,1)\times\{0\}, \qquad \pa_{\nu}\tilde{P}_{-} = 0 \text{ on }\pa (-1,1)\times(-\infty,0),
\end{align*}
where, using the estimates on $F-G$ in $\mathcal{V}_-$, we have that $g_{-}$ satisfies $e^{-\mu\xi_2}|\pa_{\xi_1}^pg_{-}(\xi)|\leq C'_\mu$, for some constant $C'_\mu$, for each $\mu<\pi/2$, and $p=0,1$. Therefore, since $\tilde{P}_{-}$ is in $\mathcal{H}(\mathcal{V})$, and the first two Neumann eigenvalues of the cross-section, $(-1,1)$, of $\mathcal{V}_-$ are $0$ and $\pi^2/4$, there exist constants $c_{\tilde{P}_-}$ and $C''_\mu$ such that
\begin{align} \label{eqn:P-minus}
    e^{-\mu\xi_2}|\pa_{\xi_1}^p\pa_{\xi_2}^q(\tilde{P}_{-}(\xi)-c_{\tilde{P}_-})| \leq C_\mu'',
\end{align}
for each $\mu<\pi/2$, $0\leq p+q \leq 1$, in $\mathcal{V}_-$. This can be seen by writing
\begin{align*}
    \tilde{P}_-(\xi) = \sum_{k=0}^{\infty}v_k(\xi_2)\cos(k\pi(\xi_1+1)/2) ,\qquad g_{-}(\xi) = \sum_{k=0}^{\infty}g_k(\xi_2)\cos(k\pi(\xi_1+1)/2)
\end{align*}
with $v_k$ satisfying
\begin{align*}
    v_k''(\xi_2) -\frac{k^2\pi^2}{4}v_k(\xi_2) = g_k(\xi_2), \qquad v_k(0) = 0, \quad v_k(\xi_2)/\xi_2 \text{ in }L^2((-\infty,-1]).
\end{align*}
By the estimates on $\pa_{\xi_1}^pg_{-}(\xi)$ for $p=0,1$, the functions $g_k(\xi_2)$ satisfy 
\begin{align} \label{eqn:gk-bounds1}
    e^{-\mu\xi_2}|g_k(\xi_2)| \leq C_\mu'k^{-1}  \text{ for each } \mu<\pi/2 .
\end{align}
This ensures that $v_0(\xi_2)$ is given by
\begin{align*}
    v_0(\xi_2) = A_0 - \int_{-\infty}^{\xi_2} tg_0(t)\,dt + \xi_2\int_{-\infty}^{\xi_2} g_0(t)\,dt,
\end{align*}
which is a constant $A_0$ plus a function $w_0(\xi_2)$ satisfying $e^{-\mu\xi_2}|w_0(\xi_2)| \leq C_\mu''$ for each $\mu<\pi/2$. The derivative $v_0'(\xi_2)$ satisfies this same exponential decay property. For $k\geq1$, we have
\begin{align*}
    v_k(\xi_2) = A_ke^{k\pi\xi_2/2} - \tfrac{1}{k\pi}e^{k\pi\xi_2/2}\int_{0}^{\xi_2} e^{-k\pi t/2}g_k(t)\,dt + \tfrac{1}{k\pi}e^{-k\pi\xi_2/2}\int_{-\infty}^{\xi_2} e^{k\pi t/2}g_k(t)\,dt,
\end{align*}
with $A_k = - \tfrac{1}{k\pi}\int_{-\infty}^{0} e^{k\pi t/2}g_k(t)\,dt$ chosen to ensure that $v_k(0) = 0$. By the exponential decay properties of $g_k(t)$ from \eqref{eqn:gk-bounds1}, and since $k\geq1$, we therefore have $|v_k(\xi_2)| \leq C_\mu''e^{\mu\xi_2}k^{-3}$ for each $\mu<\pi/2$. This allows us to sum over $k$, and obtain the estimate in \eqref{eqn:P-minus} for $p=q=0$. Taking a derivative in $\xi_2$ gives $|v_k'(\xi_2)| \leq C_\mu''e^{\mu\xi_2}k^{-2}$, and so after taking a derivative with respect to $\xi_1$ or $\xi_2$, we can still sum over $k$, and obtain the rest of the estimates in \eqref{eqn:P-minus}.

Next, $\tilde{P}_{+}(\xi):=\chi(|\xi|)\tilde{P}(\xi)$ satisfies
\begin{align} \label{eqn:P-plus}
    \Delta \tilde{P}_{+} = g_{+} \text{ in }\mathcal{V}_+, \qquad \pa_{\nu}\tilde{P}_{+} = 0 \text{ on }\{\xi_2=0\},
\end{align}
where $|g_{+}(\xi)| \leq C(1+|\xi|)^{-N}$ in $\mathcal{V}_+$. Since $\tilde{P}\in \mathcal{H}(\mathcal{V})$, we have that $\norm{\nabla\tilde{P}_{+}}_{L^2(\mathcal{V}_+)}$ is finite. Therefore, the integral of $g_+$ must be zero, and so
\begin{align*}
 \Phi(\xi) = \frac{1}{2\pi}\int_{\mathcal{V}_+}\left(\ln|\xi-t|+ \ln|\xi-t^*|\right)g_+(t)\,dt,
\end{align*}
satisfies \eqref{eqn:P-plus} with $\norm{\nabla\Phi}_{L^2(\mathcal{V}_+)}$ and $\norm{d_2^{-1}\Phi}_{L^2(\mathcal{V}_+)}$ finite.  The only solution to the homogeneous $\Delta \Psi=0$ in $\mathcal{V}_+$ with $\pa_{\nu}\Psi=0$ on $\{\xi_2=0\}$ and $\norm{d_2^{-1}\Psi}_{L^2(\mathcal{V}_+)}$ finite is the constant solution. Therefore, $\tilde{P}_+$ must be of the form
\begin{align*}
    \tilde{P}_+(\xi) = \tilde{c}_{P^+} + \frac{1}{2\pi}\int_{\mathcal{V}_+}\left(\ln|\xi-t|+ \ln|\xi-\bar{t}|\right)g_+(t)\,dt,
\end{align*}
with $t=(t_1,t_2)$ and $\bar{t} = (t_1,-t_2)$, and $\tilde{c}_{P^+}$ a constant. Since $\tilde{P}$ is defined up to a constant, we choose this constant so that $c_{\tilde{P}_+} = 0$, which now uniquely determines $\tilde{P}$.

Setting $P(\xi) = P_1(\xi)+\tilde{P}(\xi)$ therefore gives a function satisfying $\Delta P = F$ in $\mathcal{V}$, with Neumann boundary conditions on $\pa\mathcal{V}$. The estimates in \eqref{eqn:P-minus} then give the required bound on $P$ in $\mathcal{V}_-$.

We are left to prove the bounds on $P$ in $\mathcal{V}_+$. To do this, we work in complex notation. For large $|\xi|$, we can expand the $\ln|\xi-t|+\ln|\xi-\bar{t}|$ term appearing in $\tilde{P}_+(\xi)$ as
\begin{align} \label{eqn:P-plus1}
 &   2\ln|\xi| + \ln|1-t/\xi|+\ln|1-t/\bar{\xi}| \\
    & = 2\ln|\xi| -\text{Re}\left(\sum_{n=1}^{\infty}\frac{t^n(1/\xi^n + 1/\bar{\xi}^n)}{n}\right) = 2\ln|\xi| -2\sum_{n=1}^{\infty}\frac{\text{Re}(t^n)}{n}\text{Re}(1/\xi^n).
\end{align}
This allows for a further expansion of $\tilde{P}(\xi)$ 
for $|\xi|$ large in $\mathcal{V}_+$: Since $|g_{+}(t)| \leq C(1+|t|)^{-N}$,  we have that $|t^ng_+(t)|$ is integrable for $n\leq N-3$, giving an expansion of $\tilde{P}_{+}(\xi)$ in $\mathcal{V}_{+}$ up to an error of $O((1+|\xi|)^{-N+3})$. Using \eqref{eqn:P-plus1}, and that $g_+(t)$ integrates to $0$, this agrees with the expansion for $P(\xi)-\chi(|\xi|)P_+(\xi)-c_F\pi^{-1}\chi(|\xi|)\ln(|\xi|)$ given in the statement of the proposition. Taking the derivative with respect to $\xi_1$ or $\xi_2$ we also get the desired expansion for $\nabla_\xi\tilde{P}_+(\xi)$. Therefore, $P(\xi) = P_1(\xi)+\tilde{P}(\xi)$ satisfies the required  bounds in $\mathcal{V}_+$, and this completes the proof of the proposition.
\end{proof1}

 To do the matched asymptotics with $x$ in a small neighborhood of $p\sub{T}$ or $p\sub{\!B}$, we will set $\xi = (x-p\sub{T})/w$ or $(p\sub{\!B}-x)/w$.
\begin{defn} \label{defn:join-ansatz}
    We write, for $\xi \in\mathcal{V}$,
    \begin{align*}
    P_T(\xi) = \sum_{k=0}^{M}\sum_{j=0}^k w^k(\ln(w))^jp^T_{k,j}(\xi), \qquad P_B(\xi) = \sum_{k=0}^{M}\sum_{j=0}^kw^k(\ln(w))^jp^B_{k,j}(\xi).
\end{align*}
We call the functions $p^T_{k,j}(\xi)$, $p^B_{k,j}(\xi)$ admissible if they are $O((1+|\xi|)^M)$, satisfy Neumann boundary conditions on $\pa\mathcal{V}$, and if the coefficients of $w^k(\ln(w))^j$ in
\begin{align*}
    (\Delta_{\xi}+w^2\lambda(w)) P_T(\xi), \qquad  (\Delta_{\xi}+w^2\lambda(w)) P_B(\xi)
\end{align*}
are equal to $0$. The functions $P_T(\xi)$ and $P_B(\xi)$ are called admissible if $p^T_{k,j}(\xi)$, $p^B_{k,j}(\xi)$ are admissible for all $0 \leq k \leq M$, $0 \leq j \leq k$.
\end{defn}
\begin{remark} \label{rem:join-ansatz}
    In order to be admissible, we require $\Delta_{\xi}p^T_{k,j}(\xi) = G^T_{k,j}(\xi)$, for a function $G^T_{k,j}$ depending on $p^T_{m,\ell}(\xi)$ and $\lambda_{m,\ell}$ for $m \leq k-2$. In particular, $p^T_{k,j}(\xi)$ is harmonic for $k \leq 1$. The analogous is true for $p^B_{k,j}(\xi)$.

    Moreover, we have $(\Delta_{\xi}+w^2\lambda(w)) P_T(\xi)$, $(\Delta_{\xi}+w^2\lambda(w)) P_B(\xi) = o((1+|\xi|)^Mw^M)$ in $\mathcal{V}$,  for admissible $P_T(\xi)$, $P_B(\xi)$.
\end{remark}
The functions $p_{k,j}^T(\xi)$ will be defined using Proposition \ref{prop:match} once we know their asymptotic behavior as $|\xi|\to\infty$ in $\mathcal{V}$. To do this, we will work in a small neighborhood of $p\sub{T}$ (in which the boundary of $\Omega\sub{T}$ consists of a horizontal segment). Setting $\xi = (x-p\sub{T})/w$, we will write $V_Q(x_2)$ in terms of $\xi$ up to an error which is small when $w|\ln(w)||\xi|$ is small, and use this to define the asymptotic behavior of $p_{k,j}^T(\xi)$ in $\mathcal{V}_-$. We will also write $V_T(x_2)$ in terms of $\xi$, using this to define the asymptotic behavior of $p_{k,j}^T(\xi)$ in $\mathcal{V}_+$.

 The asymptotic behavior of $p^B_{k,j}(\xi)$ will be determined similarly by working in a small neighborhood of $p\sub{\!B}$, and setting $\xi= (p\sub{\!B}-x)/w$.

 We first expand $V_Q(x_2)$ in terms of $\xi$.
\begin{lemma} \label{lem:VQ}
    There exist smooth functions $\tilde{v}_{k,j}^T(\xi_2)$ such that, setting $\tilde{V}^T_Q(\xi_2) = V_Q(1+w\xi_2)$ for $-1/w<\xi_2<-1$, we have, for $p=0,1$,
    \begin{align*}
       \pa_{\xi_2}^p \tilde{V}^T_Q(\xi_2) = \sum_{k=0}^{M}\sum_{j=0}^kw^k(\ln(w))^j\pa_{\xi_2}^p\tilde{v}_{k,j}^T(\xi_2) + O(w^{M+1}(\ln(w))^{M+1}|\xi_2|^{M+1}).
    \end{align*}
        The functions $\tilde{v}_{k,j}^T(\xi_2)$ can be written as
    \begin{align*}
        \tilde{v}_{k,j}^T(\xi_2) = \tilde{v}^{T,0}_{k,j} + \tilde{v}^{T,*}_{k,j}(\xi_2),
    \end{align*}
    where $\tilde{v}^{T,0}_{k,j} = v_{k,j}(1)$ (a constant), and $\tilde{v}^{T,*}_{k,j}(\xi_2)$ is a polynomial of degree at most $k$, with coefficients determined in terms of $v_{i,j}(x_2)$ for $i \leq k-1$. In particular, $\tilde{v}^{T,*}_{k,k}(\xi_2) = 0$.

    There exist analogous functions $\tilde{v}^B_{k,j}(\xi_2)$, via setting $\tilde{V}^B_Q(\xi_2) = V_Q(-1 -w\xi_2)$, with $-1/w < \xi_2<-1$. 

    For $\lambda_{k,j}$ given for $k \leq K$, $0\leq j \leq k-1$, the functions $\tilde{v}^T_{k,j}(\xi_2)$ and $\tilde{v}^B_{k,j}(\xi_2)$ uniquely determine admissible functions $v_{i,\ell}(x_2)$ from Definition \ref{defn:evalue-ansatz} for all $0\leq i \leq K$, $0 \leq \ell\leq i$.
    \end{lemma}
\begin{remark} \label{rem:VQ}
      Note that, using Remark \ref{rem:evalue-ansatz}, if $V_Q(x_2)$ is admissible, then 
      \begin{align*}
          (\tilde{V}_Q^{T})''(\xi_2) = w^2V_Q''(x_2) = -w^2\lambda(w)\tilde{V}_Q^T(\xi_2) + o(w^{M+2}).
      \end{align*} 
      Equating the coefficients of $w^k(\ln(w)^j)$ for $0 \leq k \leq M$ gives an equation for $(\tilde{v}^T_{k,j})''(\xi_2) = (\tilde{v}^{T,*}_{k,j})''(\xi_2)$ in terms of $\tilde{v}^T_{i,\ell}(\xi_2)$ for $0 \leq i \leq k-2$. We have the analogous statement for $(\tilde{v}^B_{k,j})''(\xi_2)$.

      We will see in the proof of Lemma \ref{lem:VQ} that 
      \begin{align*}
          \tilde{v}^T_{0,0}(\xi_2) = \psi(p\sub{T}),\quad \tilde{v}^{T,*}_{1,0}(\xi_2) = v_{0,0}'(1)\xi_2 =a_1\xi_2,\quad \tilde{v}^B_{0,0}(\xi_2) = 0, \quad \tilde{v}^{B,*}_{1,0}(\xi_2) = -v_{0,0}'(-1)\xi_2=\tilde{a}_1\xi_2.
          \end{align*}
\end{remark}
\begin{proof1}{Lemma \ref{lem:VQ}}
 We first expand each $w^k(\ln(w))^jv_{k,j}(1+w \xi_2)$, as a Taylor series about $w\xi_2=0$, in powers of $w\xi_2$. The first $M-k-1$ terms in such a series  will equal $w^k(\ln(w))^jv_{k,j}(1+w\xi_2)$ up to an error of size $O(w^M(\ln(w))^{j}|\xi_2|^{M-k})$, and so we use only these terms to define the $\tilde{v}_{i,\ell}^T(\xi_2)$. The only functions of $w$ that will then appear are $w^k(\ln(w))^j$ for $0 \leq k \leq M$, $0\leq j \leq k$, and so we get the expression for $\tilde{V}_Q^T(\xi_2)$ and its derivative, given in the statement of the lemma.

  In particular, we have $\tilde{v}^T_{0,0}(\xi_2) = v_{0,0}(1) = \psi(p\sub{T})$, $\tilde{v}^T_{1,0}(\xi_2) = v_{0,0}'(1)\xi_2 + v_{1,0}(1)$. In general, the function $\tilde{v}^T_{k,j}(\xi_2)$ multiplying $w^k(\ln(w))^j$ is a polynomial of degree at most $k$, with constant term $v_{k,j}(1)$, and other terms involving  $v_{i,j}^{(k-i)}(1)\xi_2^{k-i}$, for $0 \leq i \leq k-1$ (and $j \leq i$). In the case $j=k$, this means that there are no other terms, and $\tilde{v}^T_{k,k}(\xi_2)$ is the constant $v_{k,k}(1)$. This ensures that the functions $\tilde{v}^T_{k,j}(\xi_2)$ have the form given in the statement of the lemma.

  We can do the analogous construction for $\tilde{V}^B_Q(\xi_2)$, with $\tilde{v}^B_{0,0}(\xi_2) = v_{0,0}(-1)=0$, and $\tilde{v}^B_{1,0}(\xi_2) = -v'_{0,0}(-1)\xi_2 + v_{1,0}(-1)$.

  Admissible functions $v_{k,j}(x_2)$ are uniquely determined by their values at $x_2=\pm1$, together with the differential equations they satisfy. The constant terms of $\tilde{v}^T_{k,j}(\xi_2)$ and $\tilde{v}^B_{k,j}(\xi_2)$ are given by $v_{k,j}(\pm1)$, and the right hand side of the equation that $v_{k,j}(x_2)$ satisfies depends on $\lambda_{k,j}$, as well as $\lambda_{i,\ell}$ and $v_{i,\ell}(x_2)$ for $i \leq k-1$. Therefore, for $\lambda_{k,j}$ given up to $k=K$, as the functions $\tilde{v}^T_{k,j}(\xi_2)$ and $\tilde{v}^B_{k,j}(\xi_2)$ determine $v_{k,j}(\pm 1)$, they can be iteratively used to define all the $v_{i,\ell}(x_2)$ with $i \leq K$.
\end{proof1}
 Next we expand the functions $V_T(x)$ and $V_B(x)$ in terms of $\xi$, using the structure of the Green's functions $G^T(x,p\sub{T},\lambda(w))$ and $G^B(x,p\sub{\!B},\lambda(w))$ near to $p\sub{T}$ and $p\sub{\!B}$, given by Proposition \ref{prop:Green1}.
\begin{lemma} \label{lem:VT}
Let $V_T(x)$, $V_B(x)$ be admissible functions. Then, there exist smooth functions $\tilde{u}^T_{k,j}(\xi)$ such that, setting $\tilde{V}_T(\xi) = V_T(p\sub{T}+w\xi)$ for $1<|\xi|<c/w$, $\xi_2\geq0$, we have, for $0\leq p+q\leq 1$,
\begin{align*}
   \pa_{\xi_1}^p\pa_{\xi_2}^q\tilde{V}_T(\xi) = \sum_{k=0}^{M}\sum_{j=0}^kw^k(\ln(w))^j\pa_{\xi_1}^p\pa_{\xi_2}^q\tilde{u}_{k,j}^T(\xi) + O(w^{M+1}(\ln(w))^{M+1}|\xi|^{M+1}).
\end{align*}
For $1 \leq k \leq M$, $0\leq j \leq k-1$, the functions $\tilde{u}_{k,j}^T(\xi)$ can be written as
    \begin{align*}
        \tilde{u}_{k,j}^T(\xi) = \tilde{u}^{T,0}_{k,j}(\xi) + \tilde{u}^{T,*}_{k,j}(\xi),
    \end{align*}
    where 
    \begin{align*}
      \tilde{u}^{T,0}_{k,j}(\xi) = \psi(p\sub{T})^{-1}\pi^{-1}\lambda_{k,j}\ln|\xi|+A_{k,j}^T -\pi^{-1}\sum_{m=1}^{100M}(m-1)!a^T_{k,j,m}\text{Re}(1/\xi^m) 
    \end{align*}
 and $\tilde{u}^{T,*}_{k,j}(\xi)$ is $O(|\xi|^k)$ for $|\xi|>1$, and is determined in terms of $\lambda_{i,\ell}$ and $a_{i,\ell,m}^T$ for $i \leq k-1$. Here $A_{k,j}^T$ is a constant, depending on $\lambda_{i,\ell}$ for $i\leq k$.

 For $1\leq k \leq M$, $\tilde{u}_{k,k}^T(\xi) = \tilde{u}^{T,0}_{k,k}(\xi) = \psi(p\sub{T})^{-1}\pi^{-1}\lambda_{k,k-1}$, and $\tilde{u}^{T}_{0,0}(\xi) = \tilde{u}^{T,0}_{0,0}(\xi) = \psi(p\sub{T})$.
 
    There exist analogous functions $\tilde{u}^B_{k,j}(\xi)$, via setting $\tilde{V}_B(\xi) = V_B(p\sub{\!B} -w\xi)$, with $1 < |\xi|<c/w$. For $1 \leq k \leq M$, $0\leq j \leq k-1$, the functions $\tilde{u}_{k,j}^B(\xi)$ can be written as
    \begin{align*}
        \tilde{u}_{k,j}^B(\xi) = \tilde{u}^{B,0}_{k,j}(\xi) + \tilde{u}^{B,*}_{k,j}(\xi),
    \end{align*}
    where 
    \begin{align*}
      \tilde{u}^{B,0}_{k,j}(\xi) = \pi^{-1}a_{k,j,0}^B\ln|\xi|+A_{k,j}^B -\pi^{-1}\sum_{m=1}^{100M}(m-1)!a^B_{k,j,m}\text{Re}(1/\xi^m) 
    \end{align*}
 and $\tilde{u}^{B,*}_{k,j}(\xi)$ is $O(|\xi|^k)$ for $|\xi|>1$, and is determined in terms of $\lambda_{i,\ell}$ and $a_{i,\ell,m}^B$ for $i \leq k-1$. Here $A_{k,j}^B$ is a constant, depending on $a_{i,\ell,0}^B$ for $i\leq k$. 

 For $1\leq k \leq M$, $\tilde{u}_{k,k}^B(\xi) = \tilde{u}^{B,0}_{k,k}(\xi) = \pi^{-1}a^B_{k,k-1,0}$, and $\tilde{u}_{0,0}^B(\xi) = 0$. 
\end{lemma}
\begin{remark} \label{rem:VT}
      Note, from the equation satisfied by $V_T(x)$, we have
      \begin{align*}
          \Delta_\xi\tilde{V}_T(\xi) = w^2\Delta V_T(x) = -w^2\lambda(w)V_T(x) = -w^2\lambda(w)\tilde{V}_T(\xi).
      \end{align*} 
      Therefore, equating the coefficients of $w^k(\ln(w)^j)$ for $0 \leq k \leq M$ gives an equation for $\Delta_\xi\tilde{u}^T_{k,j}(\xi)$ in terms of $\tilde{u}^T_{i,\ell}(\xi)$ for $0 \leq i \leq k-2$. Since $\tilde{u}_{k,j}^{T,0}(\xi)$ is harmonic away from $\xi=0$ (it is the real part of an analytic function),  $\Delta\tilde{u}^{T,*}_{k,j}(\xi)$ is equal to $\Delta_\xi\tilde{u}^T_{k,j}(\xi)$. The analogous statement holds for $\Delta \tilde{u}^B_{k,j}(\xi)$.

      We will see in the proof of Lemma \ref{lem:VT} that $\tilde{u}_{1,0}^{T,*}(\xi) = \pa_{x_1}\psi(p\sub{T})\xi_1$, and, $\tilde{u}^{B,*}_{1,0}(\xi)= 0$.
   \end{remark}
\begin{proof1}{Lemma \ref{lem:VT}}
We first study the contribution to $\tilde{V}_T(\xi)$ from the first term in $V_T(x)$, given by $\psi(p\sub{T})^{-1}(\lambda(w)-\lambda)G^T(x,p\sub{T},\lambda(w))$. From \eqref{eqn:Green1-a} in Proposition \ref{prop:Green1}, and the asymptotics of the Hankel function $H_0$, we have
\begin{align} \label{eqn:GT-0}
    G^T(x,p\sub{T},\lambda(w))  - \frac{1}{\lambda(w)-\lambda}\psi(x)\psi(p\sub{T}) = \frac{1}{\pi}\ln|x-p\sub{T}|(1+|x-p\sub{T}|^2f_0(x,\lambda(w))) + f_{0,1}(x,\lambda(w)),
\end{align}
where, in a neighborhood of $p\sub{T}$, both $f_0$ and $f_{0,1}$ are smooth functions of $x$ and $\lambda(w)$. For $x = p\sub{T} + w\xi$, we write out a Taylor series for $\psi(x)$, $f_0(x,\lambda(w)))$ and $f_{0,1}(x,\lambda(w))$, centered at $w\xi=0$. Here and throughout the proof, we consider the case where $|\xi|>1$. Keeping only those terms containing a factor of $w^k(\ln(w))^j$ with $k\leq M$, we have
that the function $\psi(p\sub{T})^{-1}(\lambda(w)-\lambda)G^T(x,p\sub{T},\lambda(w))$ equals
\begin{align} \nonumber
&\psi(p\sub{T})+\psi(p\sub{T})^{-1}\pi^{-1}\sum_{k=1}^{M}\sum_{j=0}^{k-1}w^k(\ln(w))^j\lambda_{k,j}\ln|\xi| + \psi(p\sub{T})^{-1}\pi^{-1}\sum_{k=1}^{M}\sum_{j=0}^{k}w^k(\ln(w))^{j}\lambda_{k,j-1} \\ \label{eqn:VT-1}
   & +\psi(p\sub{T})^{-1}\sum_{k=1}^{M}\sum_{j=0}^{k-1}w^k(\ln(w))^j\lambda_{k,j}f_{0,1}(p\sub{T},\lambda)+ \sum_{k=1}^{M}\sum_{j=0}^{k-1}w^{k}(\ln(w))^jc_{k,j,0}(\xi)
\end{align}
up to an error of size at most $O(w^{M+1}|\ln(w)|^{M+1}|\xi|^{M+1})$. Here and below, we set $\lambda_{k,-1}=0$. The function $c_{k,j,0}(\xi)$ depends on $\lambda_{i,\ell}$ for $i \leq k-1$, and can be written as the sum of a polynomial of degree $k$ in $\xi$ together with a polynomial of degree $k-1$ in $\xi$ multiplied by $\ln|\xi|$. Therefore, the first four terms contribute to $\tilde{u}^{T,0}_{k,j}(\xi)$, while the final term contributes to $\tilde{u}^{T,*}_{k,j}(\xi)$. Note that the first term is the only term independent of $w$, and that $\psi(p\sub{T})^{-1}\pi^{-1}\lambda_{k,k-1}$ is the only coefficient of $w^k(\ln(w))^k$ for $1 \leq k \leq M$. Also, the function $c_{1,0,0}(\xi)$ is given by $(\pa_{x_1}\psi)(p\sub{T})\xi_1$. In particular, this part of $\tilde{V}_T(\xi)$ has all of properties required in the lemma. 

For the contribution from the double and triple sum in  $V_T(x)$, we can ignore the contributions from the $\tfrac{1}{\lambda(w)-\lambda}\psi(x)\psi(p\sub{T})$ term in $G^T(x,p\sub{T},\lambda(w))$ from \eqref{eqn:GT-0}, due to how the coefficients $a^T_{i,j}$ are chosen in Definition \ref{defn:efn-ansatz}.

For the double sum, we again use the Taylor series for $f_0(x,\lambda(w)))$and $f_{0,1}(x,\lambda(w))$ in powers of $w\xi$, keeping only those terms leading to a factor of $w^k(\ln(w))^j$ with $k \leq M$. This ensures that the contribution from $\sum_{i=2}^{101M}\sum_{j=0}^{i-2}w^i(\ln(w))^ja^T_{i,j}G^T(x,p\sub{T},\lambda(w))$, can be written as
\begin{align} \label{eqn:VT-2}
\sum_{k=2}^{M}\sum_{j=0}^{k-1}w^k(\ln(w))^jc_{k,j,1}(\xi)
\end{align}
up to an error of size at most $O(w^{M+1}(\ln(w))^{M+1}|\xi|^M)$. Here $c_{k,j,1}(\xi)$ depends on $a_{i,\ell}^T$ for $i \leq k$ (and hence on the coefficients $a_{i,\ell,m}^T$ only for $i \leq k-1$), and on $\lambda_{i,\ell}$ for $i \leq k-2$. The function $c_{k,j,1}(\xi)$, which contributes only to $\tilde{u}^{T,*}_{k,j}(\xi)$, can be written as a sum of a polynomial and $\ln|\xi|$ multiplied by a polynomial, both of degree at most $k-2$, and so again this part of $\tilde{V}_T(\xi)$ has all of the properties  required in the lemma.

We finally study the contribution from the triple sum in $V_T(x)$. From \eqref{eqn:GT-0} and Remark \ref{rem:Green}, for $m\geq 1$, we can write (in complex notation)
\begin{align*}
    \pa_{y_1}^mG_T(x,y,\lambda(w))|_{y=p\sub{T}} = -\tfrac{(m-1)!}{\pi}\text{Re}(1/(x-p\sub{T})^m) + f_m(x,\lambda(w)),
\end{align*}
where $f_{m}(x,\lambda(w))$ has a $O(|x-p\sub{T}|^{-m+1})$ singularity at $x=p\sub{T}$. Writing $f_{m}(x,\lambda(w))$ in terms of $\xi$, and keeping only those terms leading to  a factor of $w^k(\ln(w))^j$ with $k \leq M$, we can therefore write
\begin{align*}
&\sum_{m=1}^{100M}\sum_{k=1}^{M}\sum_{\ell=0}^{k-1}w^{k+m}(\ln(w))^\ell a^T_{k,\ell,m}\pa_{y_1}^mG^T(x,y,\lambda(w))|_{y=p\sub{T}}
\end{align*}
as
\begin{align} \label{eqn:VT-3}
-\frac{1}{\pi}\sum_{m=1}^{100M}\sum_{k=1}^{M}\sum_{j=0}^{k-1}w^{k}(\ln(w))^j (m-1)!a^T_{k,j,m}\text{Re}(1/\xi^m) + \sum_{k=2}^{M}\sum_{j=0}^{k-1}w^{k}(\ln(w))^jc_{k,j,2}(\xi),
\end{align}
up to an error of at most $O(w^{M+1}(\ln(w))^{M+1}|\xi|^M)$. The first term gives the final contribution to $\tilde{u}^{0,T}_{k,j}(\xi)$. The function $c_{k,j,2}(\xi)$ depends on $a_{i,\ell,m}^T$ for $i \leq k-1$, and on $\lambda_{i,\ell}$ for $i \leq k-2$. The function $c_{k,j,2}(\xi)$, which again only contributes to $\tilde{u}^{T,*}_{k,j}(\xi)$, can be written as a sum of a rational function and $\ln|\xi|$ multiplied by a rational function, with both rational functions $O(|\xi|^{k-2})$ for $|\xi|>1$. 
Therefore, combining \eqref{eqn:VT-1}, \eqref{eqn:VT-2}, and \eqref{eqn:VT-3},  the function $\tilde{V}_T(\xi)$ has all of the required properties. 

The proof for $\tilde{V}_B(\xi)$, with $x = p\sub{\!B} -w\xi$, works in the same way: From Proposition \ref{prop:Green1}, we have
\begin{align} \label{eqn:GB-0}
    G^B(x,p\sub{\!B},\lambda(w)) = \frac{1}{\pi}\ln|x-p\sub{\!B}|(1+|x-p\sub{\!B}|^2g_0(x,\lambda(w))) + g_{0,1}(x,\lambda(w)),
\end{align}
for a constant $c^B_{0,1}$, and smooth functions $g_0$ and $g_{0,1}$. Using the part of $G^B(x,p\sub{\!B},\lambda(w))$ coming from $\frac{1}{\pi}\ln|x-p\sub{\!B}| + g_{0,1}(p\sub{\!B},\lambda)$ in the definition of $V_B(x)$ gives the contribution to $\tilde{V}_B(\xi)$ of the functions $\tilde{u}_{k,j}^{B,0}(\xi)$ (with in particular $\tilde{u}_{k,k}^{B,0}(\xi) = \pi^{-1}a_{k,k-1,0}^B$ for $k\geq1$, and $\tilde{u}_{0,0}^{B,0}(\xi) = 0$). The remaining parts of the right hand side of \eqref{eqn:GB-0} determine functions $\tilde{u}_{k,j}^{B,*}(\xi)$ with the required properties (including $\tilde{u}_{1,0}^{B,*}(\xi)=0$).
\end{proof1}
We can now combine Lemmas \ref{lem:VQ} and \ref{lem:VT} with Proposition \ref{prop:match} to inductively construct the functions $p_{k,j}^T(\xi)$ and $p_{k,j}^B(\xi)$. The asymptotics of these functions in $\mathcal{V}_+$ and $\mathcal{V}_-$ will match with those of the functions in Lemmas \ref{lem:VQ} and \ref{lem:VT} up to a high negative power of $\xi$.
\begin{prop} \label{prop:asymp}
    There exist values $\lambda_{k,j}$, and admissible functions $V_Q(x_2)$, $V_T(x)$, and $V_B(x)$ (with $\lambda_{0,0}=\lambda$, $v_{0,0}(1)=\psi(p\sub{T})$, $v_{0,0}(-1)=0$) such that, there exist admissible functions $P_T(\xi)$ and $P_B(\xi)$ with the following properties:

    There exists a constant $C_M$ such that, for all, $0 \leq k \leq M$, $0 \leq j \leq k$, $\mu<\pi/2$, $0 \leq p+q\leq 1$,
    \begin{align*}
        e^{-\mu\xi_2}|\pa_{\xi_1}^p\pa_{\xi_2}^q(\tilde{v}^T_{k,j}(\xi_2) - p^{T}_{k,j}(\xi))| +  e^{-\mu\xi_2}|\pa_{\xi_1}^p\pa_{\xi_2}^q(\tilde{v}^B_{k,j}(\xi_2) - p^{B}_{k,j}(\xi))| \leq C_{M} \, \text{ for } \xi\in\mathcal{V}_- \text{ with }\xi_2<-1,
    \end{align*}
    and
    \begin{align*}
        |\pa_{\xi_1}^p\pa_{\xi_2}^q(\tilde{u}^T_{k,j}(\xi) - p^T_{k,j}(\xi))| + |\pa_{\xi_1}^p\pa_{\xi_2}^q(\tilde{u}^B_{k,j}(\xi) - p^B_{k,j}(\xi))|\leq C_M(1+|\xi|)^{-M} \, \text{ for } \xi\in\mathcal{V}_+ \text{ with } |\xi|>1.
    \end{align*}
    \end{prop}
\begin{proof1}{Proposition \ref{prop:asymp}}
We will prove the proposition inductively, constructing values $\lambda_{k,j}$, admissible functions $v_{k,j}(x_2)$, $p_{k,j}^T(\xi)$, $p_{k,j}^B(\xi)$, and coefficients $a_{k,\ell,m}^T$, $a_{k,\ell,m}^B$ such that, for $\xi\in\mathcal{V}_-$ with $\xi_2<-1$,
\begin{align} \label{eqn:asymp0a}
        e^{-\mu\xi_2}|\pa_{\xi_1}^p\pa_{\xi_2}^q(\tilde{v}^T_{k,j}(\xi_2) - p^{T}_{k,j}(\xi))| +  e^{-\mu\xi_2}|\pa_{\xi_1}^p\pa_{\xi_2}^q(\tilde{v}^B_{k,j}(\xi_2) - p^{B}_{k,j}(\xi))| \leq C_{k}
    \end{align}
    and, for $\xi\in\mathcal{V}_+$ with $|\xi|>1$,
    \begin{align} \label{eqn:asymp0b}
        |\pa_{\xi_1}^p\pa_{\xi_2}^q(\tilde{u}^T_{k,j}(\xi) - p^T_{k,j}(\xi))| + |\pa_{\xi_1}^p\pa_{\xi_2}^q(\tilde{u}^B_{k,j}(\xi) - p^B_{k,j}(\xi))|\leq C_k(1+|\xi|)^{-100M+10k},
    \end{align}
    for all $0\leq j\leq k$, for each $0 \leq k \leq M$. 

    We first establish the base case, $k=0$. The value $\lambda_{0,0} = \lambda$ and function $v_{0,0}(x_2)$ are as given in Definition \ref{defn:evalue-ansatz}, and, as noted in Remark \ref{rem:VQ}, this gives $\tilde{v}^T_{0,0}(\xi_2) = v_{0,0}(1)=\psi(p\sub{T})$, $\tilde{v}^B_{0,0}(\xi_2) = 0$. We therefore define $p_{0,0}^T(\xi) = \psi(p\sub{T})$, $p_{0,0}^B(\xi) = 0$. These are admissible definitions for $p_{0,0}^T(\xi)$ and $p_{0,0}^B(\xi)$ since they are harmonic, and so we have the desired bounds for $k=0$.

    We now assume that $\lambda_{i,\ell}$, admissible $v_{i,\ell}(x_2)$, and $a^T_{i,\ell,m}$,  $a^B_{i,\ell,m}$ have been defined for $0 \leq i \leq k-1$ and all appropriate $\ell$, $m$, and that there exist admissible functions $p_{i,\ell}^T(\xi)$, $p_{i,\ell}^B(\xi)$, for $0\leq i \leq k-1$, $0\leq \ell\leq i$, with the required bounds from \eqref{eqn:asymp0a} and \eqref{eqn:asymp0b}. First, for $2\leq i \leq k$, $0\leq j\leq i-2$, the coefficients $a_{i,j}^T$ appearing in $V_T(x)$ are chosen as described in Definition \ref{defn:efn-ansatz}.

    We will use Proposition \ref{prop:match} as a starting point to defining the various quantities when $0 \leq j \leq k-1$, and then consider the case $j=k$.

    From Lemmas \ref{lem:VQ} and \ref{lem:VT}, the functions $\tilde{v}^{T,*}_{k,j}(\xi_2)$, $\tilde{v}^{B,*}_{k,j}(\xi_2)$ and $\tilde{u}^{T,*}_{k,j}(\xi)$, $\tilde{u}^{B,*}_{k,j}(\xi)$ are known for $0 \leq j \leq k-1$. For $0\leq j \leq k-1$, we define $F_{k,j}^T(\xi)$ to be the (known) $w^k(\ln(w))^j$ coefficient of $-w^2\lambda(w)P_T(\xi)$, and likewise for $F_{k,j}^B(\xi)$. Next we let $P_{k,j,+}^T(\xi) = \chi(|\xi|)\tilde{u}^{T,*}_{k,j}(\xi)$ and $P_{k,j,-}^T(\xi_2) = \chi(-\xi_2)\tilde{v}^{T,*}_{k,j}(\xi_2)$, and likewise for $P_{k,j,\pm}^B(\xi)$. 
    
    Note that by Remarks \ref{rem:VQ} and \ref{rem:VT}, $(\tilde{v}^{T,*}_{k,j})''(\xi_2)$ and $\Delta\tilde{u}_{k,j}^{T,*}(\xi)$ are equal to the $w^k(\ln(w))^j$ coefficient in $-w^2\lambda(w)\tilde{V}_Q^T(\xi_2)$ and $-w^2\lambda(w)\tilde{V}_T(\xi)$, and likewise with $T$ replaced by $B$. Therefore, setting
    \begin{align*}
        G^T_{k,j}(\xi) = \Delta (\chi(|\xi|)P_{k,j,+}^T(\xi) + \chi(-\xi_2)P_{k,j,-}^T(\xi_2)),
    \end{align*}
    by the inductive hypothesis, since $\chi(t) = 1$ for $t>2$, we have, for $\mu<\pi/2$,
\begin{align*}
        e^{-\mu\xi_2}|F_{k,j}^T(\xi) - G_{k,j}^T(\xi)| \leq \tilde{C}_{k} \, \text{ for } \xi\in\mathcal{V}_- \text{ with }\xi_2<-1,
    \end{align*}
    and
    \begin{align*}
        |F^T_{k,j}(\xi) - G^T_{k,j}(\xi))|\leq \tilde{C}_k(1+|\xi|)^{-100M+10(k-1)} \, \text{ for } \xi\in\mathcal{V}_+ \text{ with }|\xi|>1,
    \end{align*}
    for some constant $\tilde{C}_k$. The same estimates hold with $T$ replaced by $B$. We now apply Proposition \ref{prop:match} with the above $P^T_{k,j,\pm}(\xi)$ and $P^B_{k,j,\pm}(\xi)$, to define the functions $p^{T,*}_{k,j}(\xi)$ and $p^{B,*}_{k,j}(\xi)$, for $0\leq j \leq k-1$, satisfying $\Delta p^{T,*}_{k,j}(\xi) = F^T_{k,j}(\xi)$ and $\Delta p^{B,*}_{k,j}(\xi) = F^B_{k,j}(\xi)$. From Proposition \ref{prop:match}, we have, for $|\xi|>2$, with $\xi_2\geq0$,
    \begin{align} \label{eqn:asymp1}
        p^{T,*}_{k,j}(\xi) = P^T_{k,j,+}(\xi) +C^T_{k,j,0}\ln(|\xi|) + \sum_{m=1}^{100M-10(k-1)-5}C^T_{k,j,m}\text{Re}(1/\xi^m) + O(|\xi|^{-100M+10(k-1)+4}),
    \end{align}
    for constants $C^T_{k,j,m}$, and an analogous expression for $p^{B,*}_{k,j}(\xi)$ in terms of constants $C^B_{k,j,m}$. We can now define $\lambda_{k,j}$, for $0\leq j \leq k-1$, via
    \begin{align} \label{eqn:asymp1a}
        \psi(p\sub{T})^{-1}\pi^{-1}\lambda_{k,j} = C^T_{k,j,0},
        \end{align}
        and $a^B_{k,j,0}$ is defined via $\pi^{-1}a^B_{k,j,0} = C^B_{k,j,0}$. In particular, the constants $A^T_{k,j}$, $A^B_{k,j}$ from Lemma \ref{lem:VT} are now defined for $0 \leq j \leq k-1$, and so we set
        \begin{align*}
p^T_{k,j}(\xi) = p^{T,*}_{k,j}(\xi) + A^T_{k,j},\qquad  p^B_{k,j}(\xi) = p^{B,*}_{k,j}(\xi) + A^B_{k,j}.
        \end{align*}
        Note that since $\Delta p^T_{k,j}(\xi) = \Delta p^{T,*}_{k,j}(\xi)$, by construction, the functions $p^T_{k,j}(\xi)$ are admissible, and likewise for $p^B_{k,j}(\xi)$. Moreover, as well as the estimates that $p^T_{k,j}(\xi)$ will attain from \eqref{eqn:asymp1}, Proposition \ref{prop:match} ensures that, for $\xi\in\mathcal{V}$ with $\xi_2<-2$,
     \begin{align} \label{eqn:asymp2}
        e^{-\mu\xi_2}|p^{T}_{k,j}(\xi) - P^T_{k,j,-}(\xi_2)-D^T_{k,j}| \leq C_{\mu,k},
    \end{align}
    for each $\mu<\pi/2$, for constants $C_{\mu,k}$, $D^T_{k,j}$, and the same estimates for the first derivatives. We have the analogous estimates for $p^{B}_{k,j}(\xi)$.

    We now use the constants $C_{k,j,m}^T$ from \eqref{eqn:asymp1} and $D_{k,j}^T$ from \eqref{eqn:asymp2} to define $v_{k,j}(1)$  and $a^T_{k,j,m}$ (for $j\leq k-1$). We first define $a^T_{k,j,m}$ via $-\pi^{-1}(m-1)!a^T_{k,j,m} = C^T_{k,j,m}$ for $1 \leq m \leq 100M-10(k-1)-5$. For all other $m\leq 100M$, we set $a^T_{k,j,m}=0$.  Note that using, for $|\xi|>2$ and $\xi_2\geq0$,
    \begin{align*}
        P_{k,j,+}^T(\xi) = \tilde{u}^{T,*}_{k,j}(\xi) =  \tilde{u}^{T}_{k,j}(\xi) -  \tilde{u}^{T,0}_{k,j}(\xi),
    \end{align*}
    these definitions of $a_{k,j,m}^T$, the estimates from \eqref{eqn:asymp1} and \eqref{eqn:asymp1a}, and Lemma \ref{lem:VT}, ensure that $\tilde{u}^T_{k,j}(\xi) - p^T_{k,j}(\xi)$  satisfies the required inductive estimate from \eqref{eqn:asymp0b} for $0\leq j \leq k-1$. 
    
    We next define $\tilde{v}^{T,0}_{k,j} = v_{k,j}(1)$ to be $D^T_{k,j}$. In this case, since, for $\xi_2<-2$, 
    \begin{align*}
        P_{k,j,-}^T(\xi_2) = \tilde{v}_{k,j}^{T,*}(\xi_2) = \tilde{v}_{k,j}^T(\xi_2) - v_{k,j}(1),
    \end{align*}
    the estimate from \eqref{eqn:asymp2} ensures that that $\tilde{v}^T_{k,j}(\xi_2) - p^T_{k,j}(\xi)$  satisfies the required inductive estimate from \eqref{eqn:asymp0a} for $0\leq j \leq k-1$. 

      We define $v_{k,j}(-1)$ and $a^B_{k,j,m}$ analogously. Then, as noted in Remark \ref{rem:VQ}, the function $v_{k,j}(x_2)$ will be defined to ensure it is admissible.

      We are left to define the functions $p_{k,k}^T(\xi)$, $p_{k,k}^B(\xi)$, and  $v_{k,k}(x_2)$. We set $p_{k,k}^T(\xi) = \tilde{u}^T_{k,k}(\xi)$, which by Lemma \ref{lem:VT} is now a known constant, and so we let $\tilde{v}^T_{k,k}(\xi_2) = v_{k,k}(1)$ be this same constant. Similarly we define, $v_{k,k}(-1) = \tilde{v}_{k,k}^B(\xi_2) = \tilde{u}_{k,k}^B(\xi) = p_{k,k}^B(\xi)$. This automatically ensures that the inductive estimates hold, and then $v_{k,k}(x_2)$ is again defined to ensure that it is admissible. This completes the inductive step, and finishes the proof of the proposition. 
\end{proof1}
\begin{remark}  \label{rem:asymp}
From the proof of this proposition, and Remarks \ref{rem:VQ} and \ref{rem:VT}, we see that
\begin{align*}
    P_{1,0,+}^T(\xi) &= \chi(|\xi|)\pa_{x_1}\psi(p\sub{T})\xi_1,\quad P_{1,0,-}^T(\xi_2) = \chi(-\xi_2)v'_{0,0}(1)\xi_2, \\ P_{1,0,+}^B(\xi) &= 0,\quad P_{1,0,-}^B(\xi_2) = -\chi(-\xi_2)v'_{0,0}(-1)\xi_2,
\end{align*}
with $F_{1,0}^T(\xi) = F_{1,0}^B(\xi) = 0$. Therefore, looking at the two special cases after Proposition \ref{prop:match}, we have
\begin{align*}
    p^{T,*}_{1,0}(\xi) = v'_{0,0}(1)Y_2(\xi) + \pa_{x_1}\psi(p\sub{T}) Y_1(\xi), \qquad p^{B,*}_{1,0}(\xi) = -v'_{0,0}(-1)Y_2(\xi),
\end{align*}
and $C^T_{1,0,0} = 2\pi^{-1}v'_{0,0}(1)$ and $C^B_{1,0,0} = -2\pi^{-1}v'_{0,0}(-1)$. This means that $\lambda_{1,0} = 2\psi(p\sub{T})v'_{0,0}(1)$ and $a_{1,0,0}^B = -2v'_{0,0}(-1)$, satisfying the hypotheses of Lemma \ref{lem:approx1}.
\end{remark}
 From now on, we work with the values $\lambda_{k,j}$ and admissible functions from Proposition \ref{prop:asymp}. Using the estimates from this proposition, together with Lemmas \ref{lem:VQ} and \ref{lem:VT} gives the following estimates for $V_Q$, $V_T$, and $V_B$.
\begin{cor} \label{cor:asymp}
    For $\xi\in\mathcal{V}_-$, with $-1/w<\xi_2<-1$, we have, for $0\leq p+q\leq 1$,
    \begin{align*}
        |\pa^p_{\xi_1}\pa^q_{\xi_2}(V_Q(1+w\xi_2) - P_T(\xi))| \leq C_M(w^{M+1}|\ln(w)|^{M+1}|\xi_2|^{M+1} + e^{-\mu|\xi_2|}),
    \end{align*}
    while for $\xi\in\mathcal{V}_+$, with $1<|\xi|<c/w$, we have, for $0 \leq p+q \leq 1$,
      \begin{align*}
        |\pa^p_{\xi_1}\pa^q_{\xi_2}(V_T(p\sub{T}+w\xi) - P_T(\xi))| \leq C_M(w^{M+1}|\ln(w)|^{M+1}|\xi|^{M+1} + |\xi|^{-M}).
    \end{align*}
    The analogous estimates hold between $V_Q(-1-w\xi_2)$ and $V_B(p\sub{\!B}-w\xi)$ with $P_B(\xi)$.
\end{cor}
We will now define a function $u_w(x)$ that satisfy the hypotheses of Lemma \ref{lem:approx}. This will be constructed using the functions $V_T(x)$, $V_Q(x)$, and $V_B(x)$ outside of a $w^{1/2}$ neighborhood of $p\sub{T}$ and $p\sub{\!B}$, and the functions $P_T(\xi)$, $P_B(\xi)$ within these neighborhoods.
\begin{defn} \label{defn:final}
Let $\chi_T(x) = \chi_T(x;w)$ be a smooth cut-off function, such that:
\begin{enumerate}
    \item[i)] for $x_2\geq1$, $\chi_T(x)$ depends only on $|x-p\sub{T}|$, and equals $1$ for $|x-p\sub{T}|<w^{1/2}$ and $0$ for $|x-p\sub{T}|>2w^{1/2}$, with $j$-th derivative bounded by $Cw^{-j/2}$;
    \item[ii)] for $x_2<1$, $\chi_T(x)$ depends only on $x_2$, and  equals $1$ for $|x_2-1|<w^{1/2}$ and $0$ for $|x_2-1|>2w^{1/2}$, with $j$-th derivative bounded by $Cw^{-j/2}$;
    \item[iii)] for $w>0$ sufficiently small, the support of $\chi_T(x)$ for $x=(x_1,1)$ contains $(-w,w)$, is contained in $(-a,a)$, and $\chi_T(x)$  satisfies Neumann boundary conditions on $\pa\Omega(w)$.
\end{enumerate}
A function $\chi_B(x)$ supported in a neighborhood of $p\sub{\!B}$ is defined analogously.

Then, for $x\in\Omega(w)$ and for $w>0$ sufficiently small, we define 
\begin{align*}
    u_w(x) =& \,(1-\chi_T(x))V_T(x)1_{\Omega\sub{T}}(x) +(1-\chi_B(x))V_B(x)1_{\Omega\sub{B}}(x) \\&+ (1-\chi_T(x))(1-\chi_B(x))V_Q(x_2)1_{Q(w)}(x) 
        + \chi_T(x)P_{T}((x-p\sub{T})/w) + \chi_B(x) P_{B}((p\sub{\!B}-x)/w).
\end{align*}
\end{defn}
We can use Corollary \ref{cor:asymp}, and the equations satisfied by $V_Q(x)$, $V_T(x)$, and $V_B(x)$, together with $P_T(\xi)$ and $P_B(\xi)$, to bound $(\Delta + \lambda(w))u_w$.
\begin{prop} \label{prop:final}
The function $u_w$ from Definition \ref{defn:final} satisfies $\norm{u_w}_{L^2(\Omega(w))} = 1+o(1)$, Neumann boundary conditions on $\pa\Omega(w)$, and
\begin{align*}
    (\Delta + \lambda(w))u_w = f_w \text{ in }\Omega(w),
\end{align*}
for $f_w$ with $\norm{f_w}_{L^2(\Omega(w))} = o(w^{M/2-2})$.
\end{prop}
\begin{proof1}{Proposition \ref{prop:final}}
Note that by construction, since $\psi$ is $L^2$-normalized in $\Omega\sub{T}$, we have $\norm{u_w}_{L^2(\Omega(w))} = 1 + o(1)$, and moreover $\pa_{\nu}u_{w}=0$ on $\pa\Omega(w)$. We write $(\Delta+\lambda(w))u_{w}= f_{w,1}+f_{w,2}$, where  $f_{w,1}$ is the contribution from when no derivatives are applied to the cut-off functions $\chi_T(x)$ or $\chi_B(x)$. By Remarks \ref{rem:evalue-ansatz}, \ref{rem:efn-ansatz}, and \ref{rem:join-ansatz}, we have
\begin{align*}
    (\Delta+\lambda(w))V_Q(x) & = o(w^{M}) \text{ on the support of }(1-\chi_T(x))(1-\chi_B(x))1_{Q(w)}(x), \\
    (\Delta+\lambda(w))V_T(x) & = 0 \text{ on the support of } (1-\chi_T(x))1_{\Omega\sub{T}}(x),\\
    (\Delta+\lambda(w))V_B(x) & = 0 \text{ on the support of } (1-\chi_B(x))1_{\Omega\sub{B}}(x),\\
    (\Delta_{\xi}+w^2\lambda(w))P_T(\xi) &= o(w^{M/2}) \text{ and } (\Delta_{\xi}+w^2\lambda(w))P_B(\xi) = o(w^{M/2})  \text{ for } \xi\in\mathcal{V} \text{ with }|\xi| \leq 2w^{-1/2}.
\end{align*}
Therefore, since $|x-p\sub{T}|/w <2w^{-1/2}$ on the support of $\chi_T(x)$ and $|p\sub{\!B}-x|/w <2w^{-1/2}$ on the support of $\chi_B(x)$, we have $f_{w,1}= o(w^{M/2-2})$. 

Applying Corollary \ref{cor:asymp}, on the support of the derivatives of $\chi_T(x)$ (so that $w^{-1/2}<|\xi|<2w^{-1/2}$), the functions $1_{\Omega\sub{T}}(x)(V_T(x)-P_T((x-p\sub{T})/w)))$ and $1_{Q(w)}(x)(V_Q(x_2)-P_T((x-p\sub{T})/w))$ are $O(w^{M/2})$, while their gradients are $O(w^{M/2-1})$. Each time we take a derivative in $x_1$ or $x_2$, we lose a factor of at most $O(w^{-1/2})$. We have the analogous statement on the support of the derivatives of $\chi_B(x)$, and so $f_{w,2} = O(w^{M/2-1})$, finishing the proof.
\end{proof1}
We now choose $M>8$, and apply Lemma \ref{lem:approx} for $N=\tfrac{1}{2}M-2>2$, with $\lambda = \lambda_j = \mu_k(\Omega\sub{T})$, so that
\begin{align} \label{eqn:final1}
    |\lambda(w)-\lambda_j^w| = O(w^2), \qquad \norm{\varphi^w_j - c_w u_w}_{H^1(\Omega(w))} = O(w^2),
\end{align}
for a constant $c_w = 1+o(1)$. As we noted in Remark \ref{rem:asymp}, the coefficients $\lambda_{1,0}$ and $a_{1,0,0}^B$ satisfy the hypotheses of Lemma \ref{lem:approx1}, and so combining \eqref{eqn:final1} with the estimates in Lemma \ref{lem:approx1} completes the proof of Proposition \ref{prop:approx1}.
\\
\\
\indent\textbf{Case 2. $\lambda_j=\tau_n$ and the proof of Proposition \ref{prop:approx2}.}  To simplify notation, we will write $\lambda = \tau_n$, and denote a $L^2([-1,1])$-Dirichlet eigenfunction  of eigenvalue $\lambda$ in $[-1,1]$ by $\gamma(x_2)$. The proof of Proposition \ref{prop:approx2} will follow closely the proof of Proposition \ref{prop:approx1} above. The differences in the proof come from the following two facts:
\begin{enumerate}
    \item[1)] If $\lambda(w) = \lambda + o(1)$, then for small $w$, $\lambda(w)$ is separated from the Neumann spectrum of $\Omega\sub{T}$ and $\Omega\sub{B}$, and so  the Green's functions $G^T(x,p\sub{T},\lambda(w))$ and $G^B(x,p\sub{\!B},\lambda(w))$ can both be written as in \eqref{eqn:GB-0}.

    \item[2)] The boundary value problem
    \begin{align*}
        v''(x_2) + \lambda v(x_2) = F(x_2) \text{ for } x_2\in[-1,1], \qquad v(-1) = a, \, v(1) = b
    \end{align*}
    has a solution only if
    \begin{align*}
        \int_{-1}^{1}F(x_2)\gamma(x_2)\,dx_2 = -[\gamma'(x_2)v(x_2)]_{-1}^{1}
    \end{align*}
    and this solution can be uniquely defined by requiring that $v(x_2)$ is orthogonal to $\gamma(x_2)$.
\end{enumerate}
We will use the solvability condition in 2) to determine the coefficients in the approximate eigenvalue $\lambda(w)$. We begin by recording the change in the eigenvalue and eigenfunction Ansatze in the neck, compared to Definition \ref{defn:evalue-ansatz}. As seen in Theorem 2.5(b) of \cite{arrieta95}, the eigenvalue $\lambda_j^w$ now satisfies $\lambda_j^w = \lambda + O(w|\ln(w)|)$. Therefore, as well as different leading order terms, a main difference is to include $w^k(\ln(w))^k$ terms in $\lambda(w)$. We also now incorporate the above solvability condition in order for a function $V_Q(x_2)$ to be admissible.
\begin{defn} \label{defn:evalue-ansatz-neck}
    Fixing a large integer $M$, we write
    \begin{align*}
         \lambda(w)  = \lambda_{0,0} + \sum_{k=1}^{M}\sum_{j=0}^{k}w^k(\ln(w))^j\lambda_{k,j} , \qquad
    V_Q(x_2)  = \sum_{k=0}^{M}\sum_{j=0}^{k}w^k(\ln(w))^jv_{k,j}(x_2),
    \end{align*}
    for $\lambda_{k,j}\in\R$ and $v_{k,j}\in C^{\infty}([-1,1])$. We set $\lambda_{0,0}=\lambda$ and $v_{0,0}(x_2) = \gamma(x_2)$.
     For $(k,j)\neq(0,0)$, we call the function $v_{k,j}(x_2)$ admissible if it is orthogonal to $\gamma(x_2)$ and the coefficient of $w^k(\ln(w))^j$ in $V_Q''(x_2)+\lambda(w)V_Q(x_2)$ is equal to $0$, and call $V_Q(x_2)$ admissible if $v_{0,0}(x_2)=\gamma(x_2)$ and $v_{k,j}(x_2)$ is admissible for all $0\leq k \leq M$, $0\leq j \leq k$ with $(k,j)\neq (0,0)$.
    \end{defn}
\begin{remark} \label{rem:evalue-ansatz-neck}
   For the function $v_{k,j}(x_2)$ to be admissible, it will need to satisfy 
    \begin{align*}
        v_{k,j}''(x_2)+\lambda v_{k,j}(x_2) = -\lambda_{k,j}\gamma(x_2) +H_{k,j}(x_2)
    \end{align*}
    for a function $H_{k,j}(x_2)$ depending on $v_{m,\ell}$ and $\lambda_{m,\ell}$ for $m\leq k-1$. In particular, using 2) above,  in order for $v_{k,j}(x_2)$ to be admissible, we require that
    \begin{align} \label{eqn:compatible1}
       \lambda_{k,j} = [\gamma'(x_2)v_{k,j}(x_2)]_{-1}^{1} + \int_{-1}^{1}\gamma(x_2)H_{j,k}(x_2)\,dx_2.
    \end{align}
   Moreover, for $V_Q(x_2)$ admissible, we have $(\Delta +\lambda(w))V_Q(x_2) = o(w^{M})$, where here and throughout, the implicit constant depends on the $\lambda_{k,j}$ and $v_{k,j}(x_2)$.
\end{remark}
Since $\lambda$ is not in the Neumann spectrum of $\Omega\sub{T}$ or $\Omega\sub{B}$, we use the same structure of Ansatz in $\Omega\sub{T}$ and $\Omega\sub{B}$ that we did for $V_B(x)$ in the previous case. That is,
\begin{align} \label{eqn:efn-ansatz-neckT}
V_T(x) = \sum_{m=0}^{100M}\sum_{k=1}^{M}\sum_{\ell=0}^{k-1}w^{k+m}(\ln(w))^\ell a^T_{k,\ell,m}\pa_{y_1}^mG^T(x,p\sub{T},\lambda(w)), \text{ for } x\in\Omega\sub{T}\backslash\{p\sub{T}\}, \\ \label{eqn:efn-ansatz-neckB}
V_B(x) = \sum_{m=0}^{100M}\sum_{k=1}^{M}\sum_{\ell=0}^{k-1}w^{k+m}(\ln(w))^\ell a^B_{k,\ell,m}\pa_{y_1}^mG^B(x,p\sub{\!B},\lambda(w)), \text{ for } x\in\Omega\sub{B}\backslash\{p\sub{\!B}\} .   
\end{align}
As in Lemma \ref{lem:approx1}, we can reduce the proof of Proposition \ref{prop:approx2} to constructing $\lambda(w)$, $V_Q$, $V_T$, $V_B$ with appropriate leading order terms. The proof of the following lemma follows exactly as for the proof of Lemma \ref{lem:approx1}
\begin{lemma} \label{lem:approx1-neck}
Suppose that $V_Q$, $V_T$, $V_B$ are admissible with the coefficient $\lambda_{1,1}$ in $\lambda(w)$  equal to $n^2\pi$, and the coefficients $a_{1,0,0}^T=2{b}_1$ and $a_{1,0,0}^B=2\tilde{b}_1$. Here $b_1$ and $\tilde{b}_1$ are as in Definition \ref{defn:approx2}. Then,
    \begin{align*}
        & |\lambda(w) - \tau|  = O(w), \quad \norm{(2w)^{-1/2}V_Q - X_Q}_{H^1(Q(w)\backslash (D_\eta(p\sub{T})\cup D_\eta(p\sub{\!B})))} = O(w|\ln(w)|),  \\
   & \norm{(2w)^{-1/2}V_T - X_T}_{H^1(\Omega\sub{T}\backslash D_\eta(p\sub{T}))}  + \norm{(2w)^{-1/2}V_B - X_B}_{H^1(\Omega\sub{B}\backslash D_\eta(p\sub{\!B}))} = O(w^{3/2}|\ln(w)|),
    \end{align*}
   where the implicit constants depend on $\eta$.
\end{lemma}
We will construct admissible $V_Q$, $V_T$, $V_B$ again using matched asymptotics and Proposition \ref{prop:match}. To do this matching, we will use functions $P_T(\xi)$ and $P_B(\xi)$ from Definition \ref{defn:join-ansatz}.

Writing $\tilde{V}_Q^T(\xi) = V_Q(1+w\xi_2)$, $\tilde{V}_Q^B(\xi) = V_Q(-1-w\xi_2)$, as in Lemma \ref{lem:VQ}, we can construct functions $\tilde{v}_{k,j}^T(\xi_2)$ and $\tilde{v}_{k,j}^B(\xi_2)$. This time, using the properties of $v_{0,0}(x_2) = \gamma(x_2)$, we have
\begin{align} \label{eqn:neck-leading1}
          \tilde{v}^T_{0,0}(\xi_2) = 0,\quad \tilde{v}^{T,*}_{1,0}(\xi_2) = v_{0,0}'(1)\xi_2 =b_1\xi_2,\quad \tilde{v}^B_{0,0}(\xi_2) = 0, \quad \tilde{v}^{B,*}_{1,0}(\xi_2) = -v_{0,0}'(-1)\xi_2=\tilde{b}_1\xi_2.
          \end{align}
The functions $\tilde{v}_{k,j}^T(\xi_2)$ and $\tilde{v}_{k,j}^B(\xi_2)$ determine the values of ${v}_{k,j}(\pm1)$. Therefore, for $v_{m,\ell}$, $\lambda_{m,\ell}$ given for $m\leq k-1$, and by choosing $\lambda_{k,j}$ to satisfy the compatibility condition in \eqref{eqn:compatible1}, the functions $\tilde{v}_{k,j}^T(\xi_2)$ and $\tilde{v}_{k,j}^B(\xi_2)$ can be used to define an admissible function $v_{k,j}(x_2)$.

We can also write $\tilde{V}_T(\xi) = V_T(p\sub{T}+w\xi)$, $\tilde{V}_B(\xi) = V_B(p\sub{\!B}-w\xi)$, and construct functions $\tilde{u}_{k,j}^T(\xi)$ and $\tilde{u}_{k,j}^B(\xi)$ both with the same properties as in the analysis of $\tilde{V}_B(\xi)$ in Lemma \ref{lem:VT}. In particular,
\begin{align}  \label{eqn:neck-leading2}
    \tilde{u}_{0,0}^T(\xi) = \tilde{u}_{1,0}^{T,*}(\xi)= 0, \, \tilde{u}_{k,k}^T(\xi) = \pi^{-1}a^T_{k,k-1,0} \text{ for } 1 \leq  k\leq M,\\  \label{eqn:neck-leading3}
     \tilde{u}_{0,0}^B(\xi) = \tilde{u}_{1,0}^{B,*}(\xi)= 0, \, \tilde{u}_{k,k}^B(\xi) = \pi^{-1} a^B_{k,k-1,0} \text{ for } 1 \leq  k\leq M.
\end{align}
We now describe how to adapt the proof of Proposition \ref{prop:asymp} to construct $\lambda(w)$, and $V_Q(x)$, $V_T(x)$, $V_B(x)$, together with $P_T(\xi)$, $P_B(\xi)$, all as admissible functions, and with the same bounds as in the statement of Proposition \ref{prop:asymp}.

As in the proof of Proposition \ref{prop:asymp}, we define values $\lambda_{i,j}$, and admissible functions $v_{i,j}(x_2)$, $p_{i,j}^T(\xi)$, $p_{i,j}^B(\xi)$, and coefficients $a_{i,j,m}^T$, $a_{i,j,m}^B$ inductively on $i$. The base case, $i=0$ is now given by $\lambda_{0,0} = \lambda$, $v_{0,0}(x_2) = \gamma(x_2)$ as in Definition \ref{defn:evalue-ansatz-neck} above, and $p_{0,0}^T(\xi) = p_{0,0}^B(\xi) = 0$.

Assuming that everything has been defined for $0 \leq i \leq k-1$, we will now show to construct these coefficients and functions for $i=k$, again first considering $0 \leq j \leq k-1$, and then the case $j=k$. For $0 \leq j \leq k-1$, following the proof of Proposition \ref{prop:asymp} up to \eqref{eqn:asymp1}, we find that \eqref{eqn:asymp1} still holds for the functions $p^{T,*}_{k,j}(\xi)$, and analogously for $p^{B,*}_{k,j}(\xi)$.

Now, unlike before, we do not determine $\lambda_{k,j}$ at this stage (as we did in \eqref{eqn:asymp2}), but instead define the coefficients $a_{k,j,0}^T$ and $a_{k,j,0}^B$ as $\pi C_{k,j,0}^T$ and $\pi C_{k,j,0}^B$ respectively. We next continue to follow the proof of Proposition \ref{prop:asymp}, defining $v_{k,j}(\pm1)$ for $0\leq j\leq k-1$, and all remaining coefficients $a^T_{k,j,m}$ and $a^B_{k,j,m}$. We are yet to define $\lambda_{k,j}$, and so we define it using \eqref{eqn:compatible1}, to obtain an admissible function $v_{k,j}(x_2)$. 

Finally, we define (using \eqref{eqn:neck-leading2} and \eqref{eqn:neck-leading3}),
\begin{align} \label{eqn:neck-leading4}
    v_{k,k}(1) =\tilde{u}^T_{k,k}(\xi) = \pi^{-1}a^T_{k,k-1,0},\quad v_{k,k}(-1) =\tilde{u}^B_{k,k}(\xi)=\pi^{-1} a^B_{k,k-1,0}
\end{align} 
as in the original proof. Then, we again use \eqref{eqn:compatible1} to define $\lambda_{k,k}$, and hence obtain a function $v_{k,k}(x_2)$ that is admissible. This completes the inductive step, and allows us to construct $\lambda(w)$, and admissible $V_Q(x)$, $V_T(x)$, $V_B(x)$, with the desired bounds.

We also have a version of Remark \ref{rem:asymp} in this case, which gives the values of some of the leading order terms: From \eqref{eqn:neck-leading1}, we now have
\begin{align*}
        p^{T,*}_{1,0}(\xi) = v'_{0,0}(1)Y_2(\xi) = b_1Y_2(\xi) , \qquad p^{B,*}_{1,0}(\xi) = -v'_{0,0}(-1)Y_2(\xi) = \tilde{b}_1Y_2(\xi),
\end{align*}
and so $a^T_{1,0,0} = 2b_1$ and $a^B_{1,0,0} = 2\tilde{b}_1$. Therefore, from \eqref{eqn:neck-leading4} (with $k=1$), we have
\begin{align*}
    v_{1,1}(1) = \pi^{-1}a_{1,0,0}^T =  2\pi^{-1}b_1 = 2\pi^{-1}\gamma'(1),\quad v_{1,1}(-1) = \pi^{-1}a_{1,0,0}^B =  2\pi^{-1}\tilde{b}_1=-2\pi^{-1}\gamma'(-1).
\end{align*}
Since $v_{1,1}(x_2)$ satisfies the equation $v''_{1,1}(x_2) + \lambda v_{1,1}(x_2)=-\lambda_{1,1}\gamma(x_2)$, from \eqref{eqn:compatible1} we therefore have
\begin{align*}
    \lambda_{1,1} = [\gamma'(x_2)v_{1,1}(x_2)]_{-1}^1 = 2\pi^{-1}(\gamma'(1))^2+2\pi^{-1}(\gamma'(-1))^2.
\end{align*}
Since $\gamma(x_2) = \sin(n\pi (x_2+1)/2)$, the above formula gives $\lambda_{1,1} = n^2\pi$.

We can now define $u_w(x)$ as in Definition \ref{defn:final}, but now with each term multiplied by a factor of $(2w)^{-1/2}$ to ensure that $\norm{u_w}_{L^2(\Omega(w))} = 1+o(1)$. Then, using the proof of Proposition \ref{prop:final}, we get that
\begin{align*}
    (\Delta + \lambda(w))u_w = f_w \text{ in }\Omega(w), \quad \pa_{\nu}u_w = 0 \text{ on }\pa\Omega(w),
\end{align*}
with $\norm{f_w}_{L^2(\Omega(w))} = o(w^{M/2-2-1/2})$. We choose $M>8$ and apply Lemma \ref{lem:approx} for $N=\tfrac{1}{2}M-2-\tfrac{1}{2}>\tfrac{3}{2}$, with $\lambda = \lambda_j = \tau_n$, so that
\begin{align} \label{eqn:final1-neck}
    |\lambda(w)-\lambda_j^w| = O(w^{3/2}), \qquad \norm{\varphi^w_j - c_w u_w}_{H^1(\Omega(w))} = O(w^{3/2}),
\end{align}
for a constant $c_w = 1+o(1)$. Since the coefficients $\lambda_{1,1}$, $a_{1,0,0}^T$, and $a_{1,0,0}^B$ satisfy the hypotheses of Lemma \ref{lem:approx1-neck}, we can combine \eqref{eqn:final1-neck} with the estimates in Lemma \ref{lem:approx1-neck} to complete the proof of Proposition \ref{prop:approx2}.

\section{Numerical Experiments}\label{sec:numExp}

We conclude with three numerical experiments illustrating the nodal phenomena
studied above and the role of the graph structure suggested by
Conjecture \ref{conj:1}. In both experiments, the chain consists of two
rectangular subdomains: $\Omega\sub{T}$ has width $1$ and height $2$, while
$\Omega\sub{B}$ has width $1$ and height $1.8$. We compute the Neumann
eigenfunctions using the {\it Matlab} PDE toolbox.

\begin{figure}[!htbp]
\includegraphics[scale=0.4]{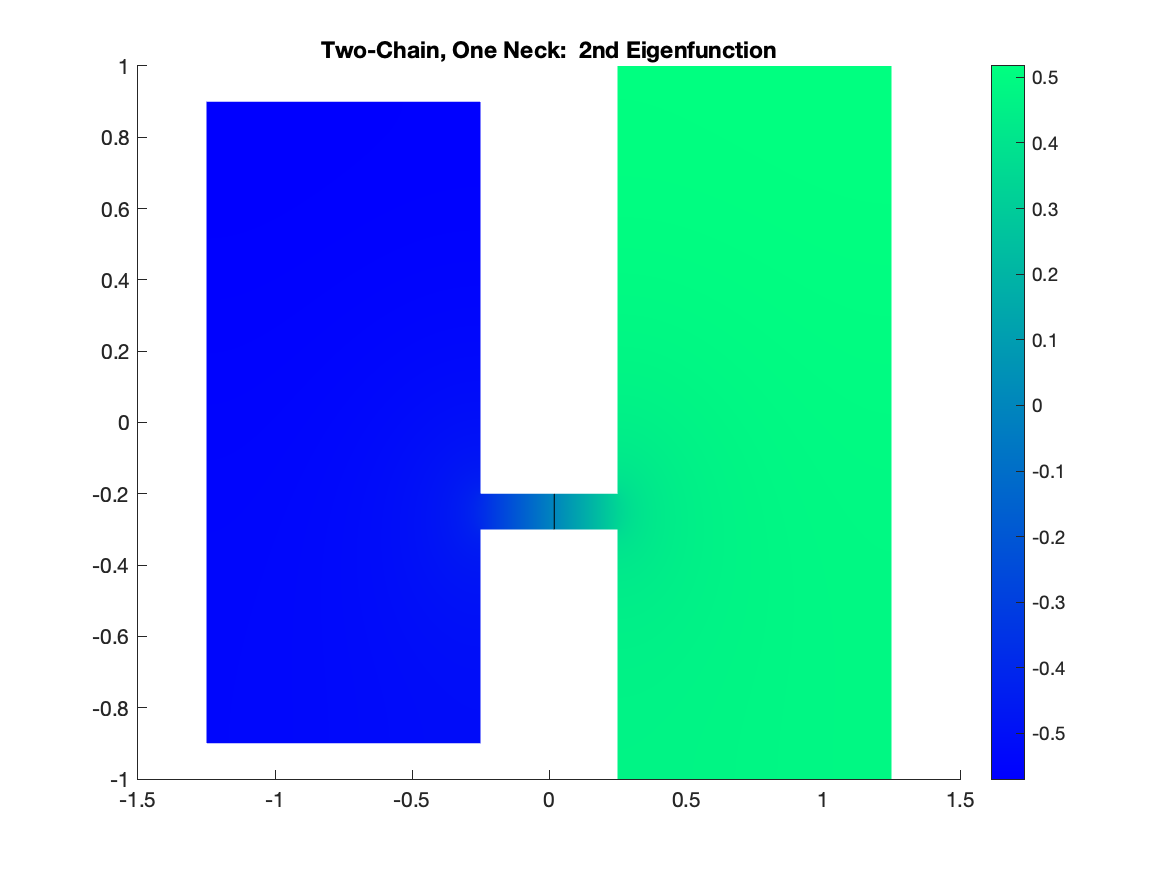}\hspace{-.25in}
\includegraphics[scale=0.4]{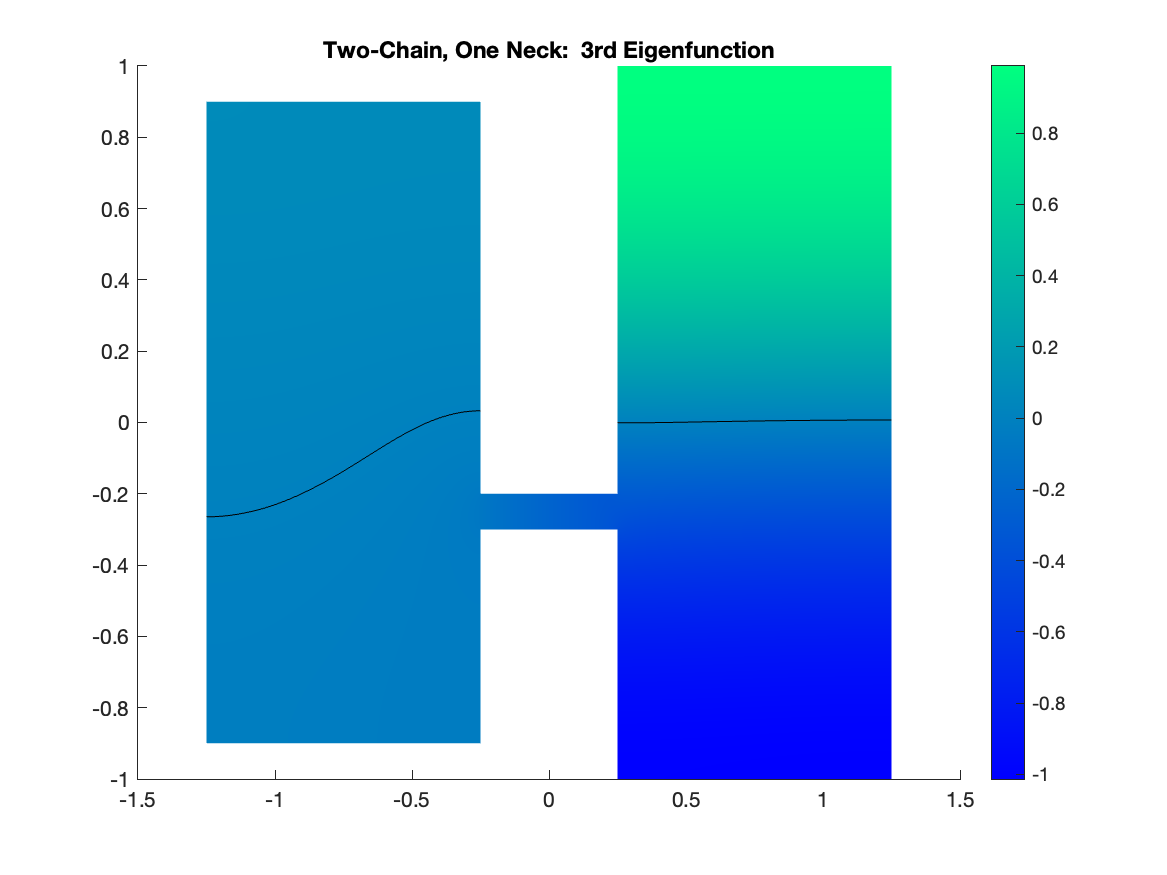} \\
\includegraphics[scale=0.4]{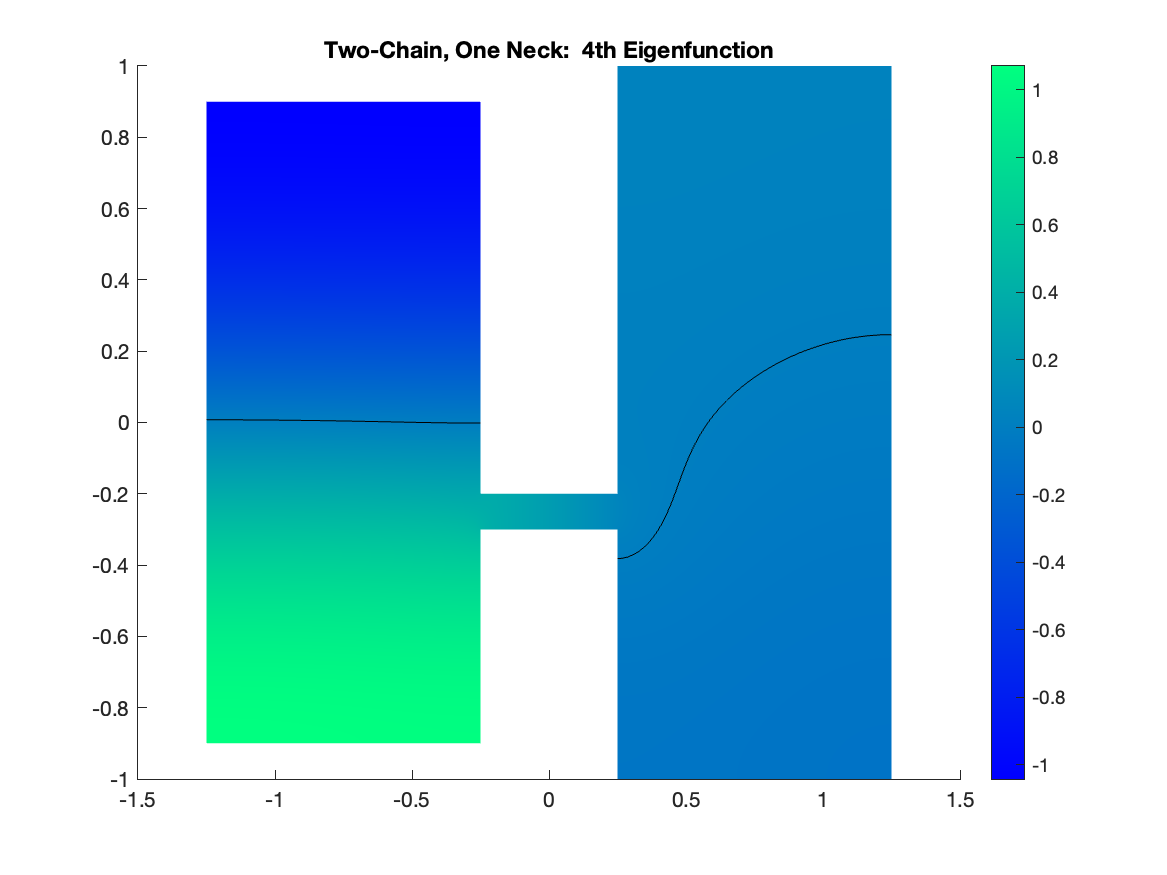}\hspace{-.25in}
\includegraphics[scale=0.4]{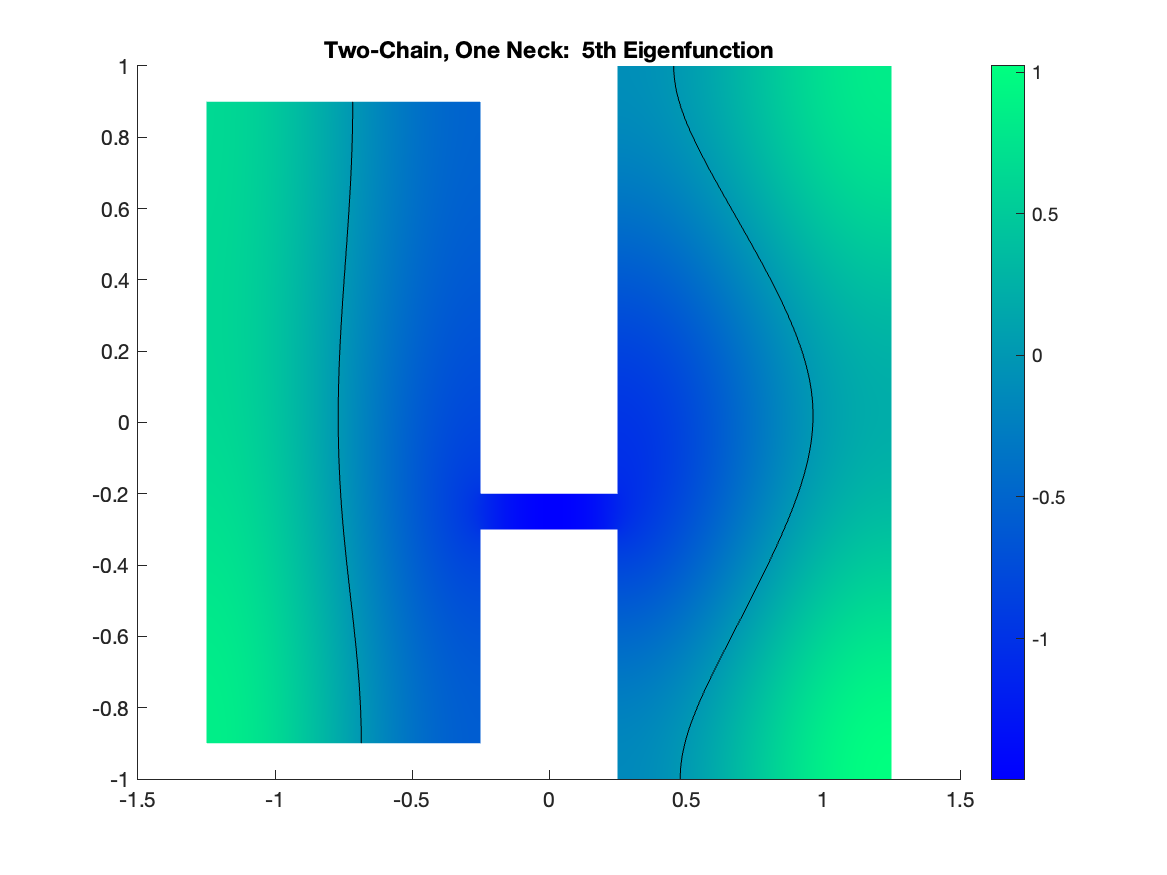}
\caption{The $2$nd (top left), $3$rd (top right), $4$th (bottom left), and
$5$th (bottom right) Neumann eigenfunctions of the $2$-chain with one neck.
The associated graph is a tree. The $2$nd and $3$rd eigenfunctions are
Courant sharp, while the $4$th and $5$th are not.
}
\label{fig:num1}
\end{figure}

The first experiment considers a single neck joining the two subdomains.
The associated graph has two vertices joined by one edge and is therefore a
tree. The neck has length $0.5$ and width $0.1$, and intersects each rectangle
at height $0.75$ from its lower edge. The finite element mesh size is $0.005$;
refining the mesh did not change the nodal counts reported below.
Figure \ref{fig:num1} shows the second through fifth eigenfunctions. The
second and third eigenfunctions are Courant sharp, in agreement with the
low-energy rigidity predicted for the tree case. The fourth and fifth eigenfunctions are
not Courant sharp, showing that this low-energy rigidity has already broken
down by the fourth mode in this example.

As an indication that the neck width does not change the outcome, in a second experiment, we fix the box dimensions and neck location as in Figure \ref{fig:num1} and consider the two chain, one neck set-up with fixed neck length $.5$ but different neck widths $.2, .1, .05,.025$ and plot the 3rd eigenfunction.  It is clear from this figure that the nodal domain structure is quite static with respect to neck width. We again use a finite element mesh size of $.005$.

\begin{figure}[!htbp]
\includegraphics[width=.45\linewidth]{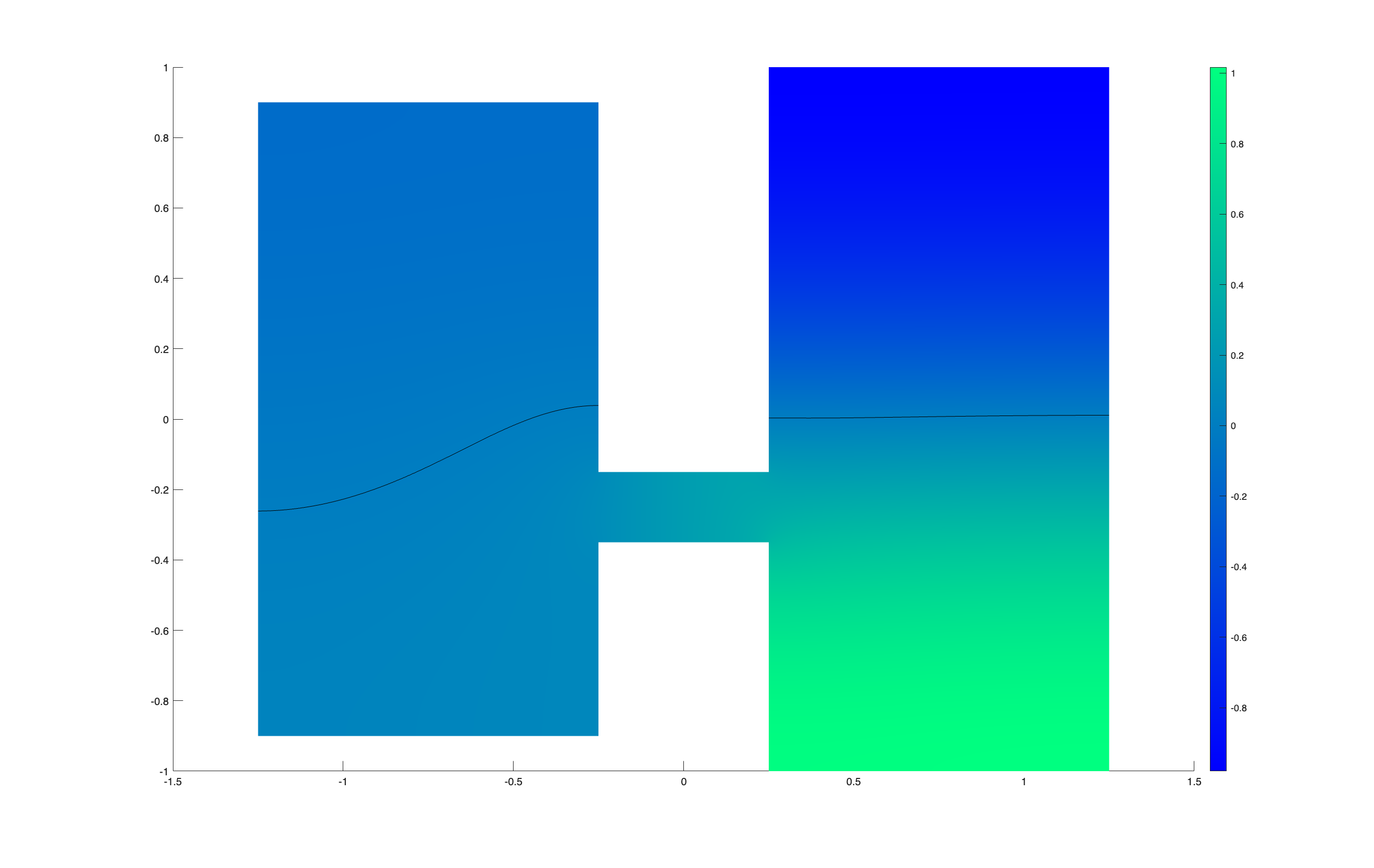}
\includegraphics[width=.45\linewidth]{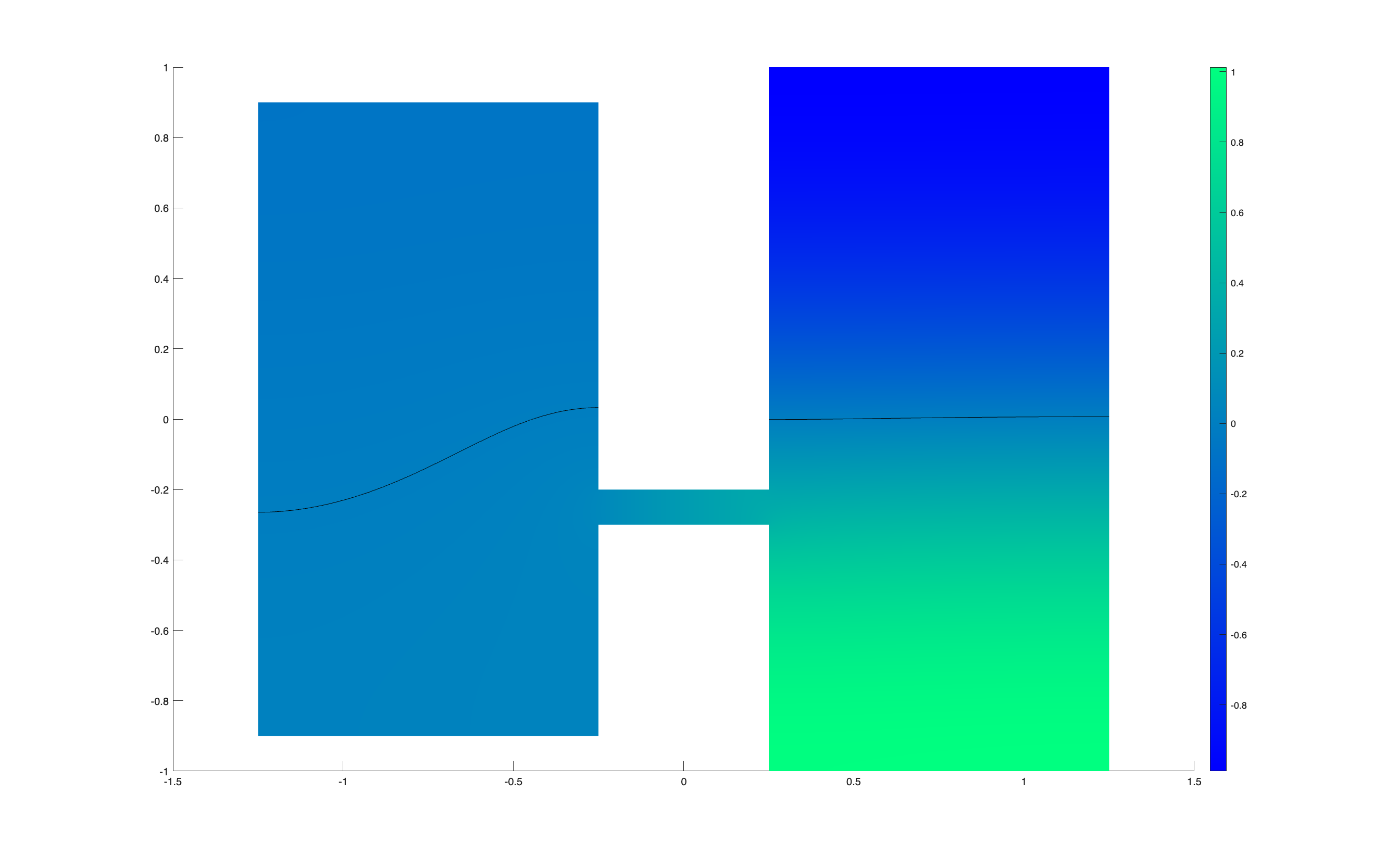} \\
\includegraphics[width=.45\linewidth]{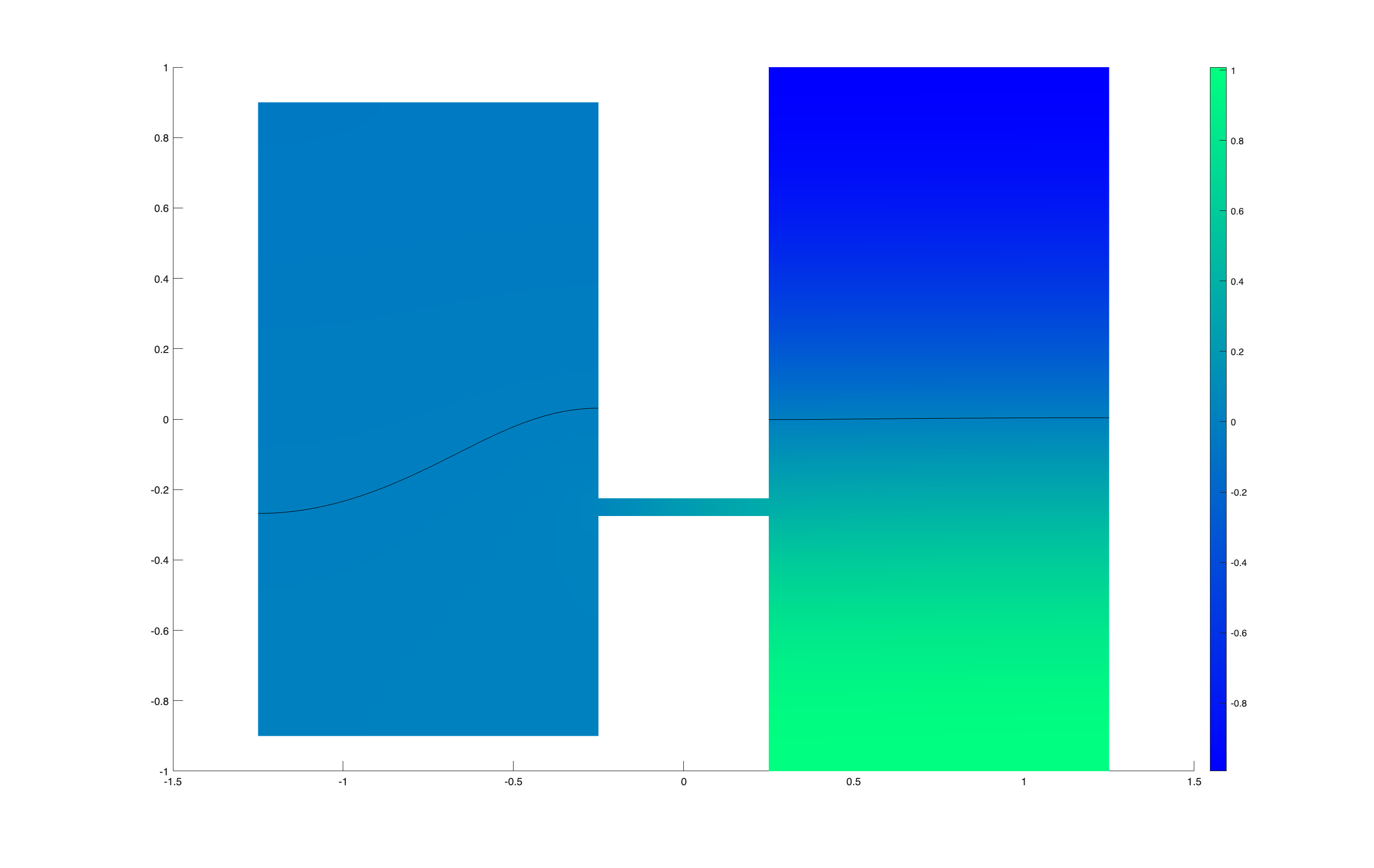}
\includegraphics[width=.45\linewidth]{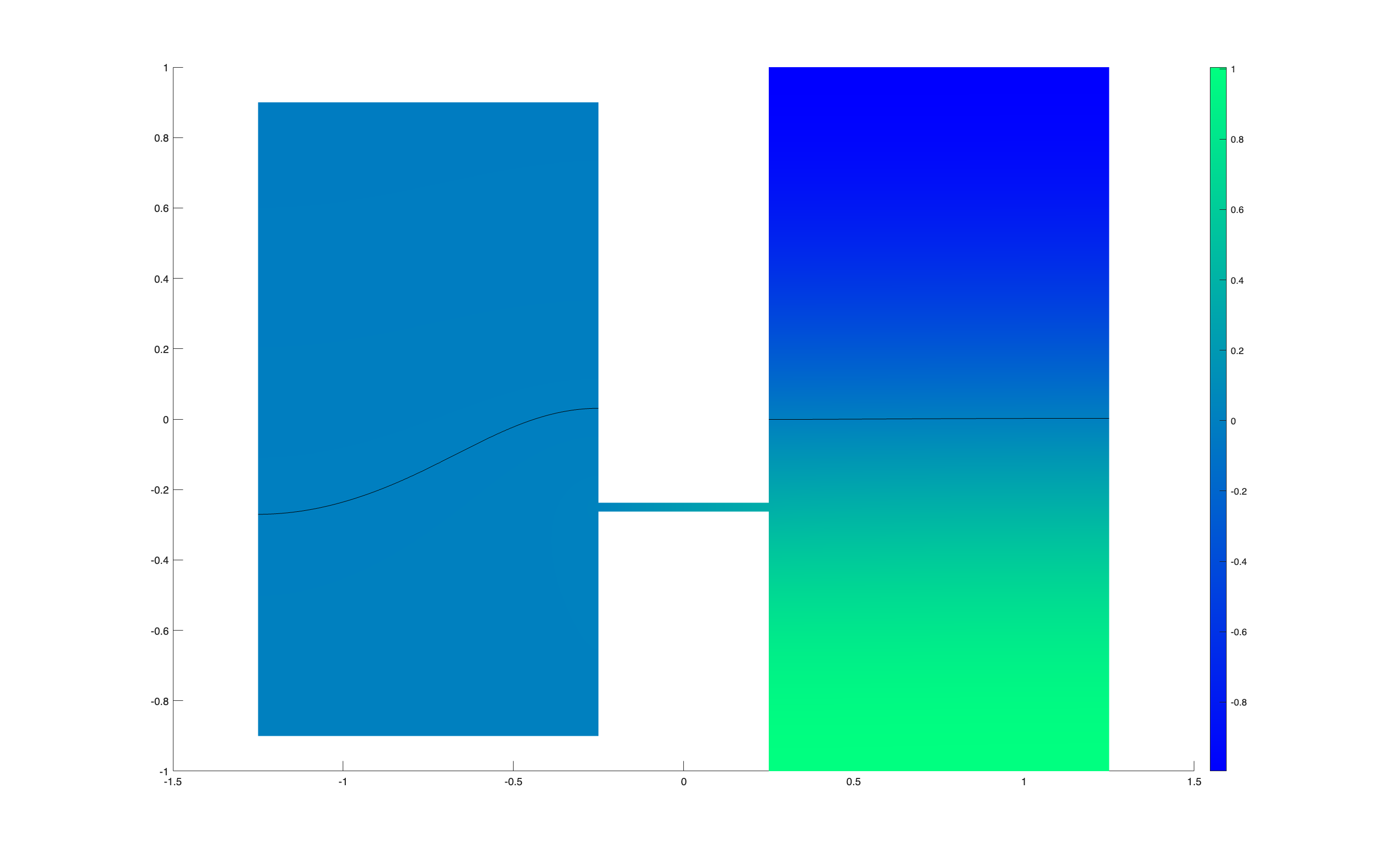}
\caption{The $3$rd Neumann eigenfunction of the $2$-chain with one neck associated to neck widths (Top Left) $.2$, (Top Right) $.1$, (Bottom Left) $.05$ and (Bottom Right) $.025$.
The associated graph is a tree. $3$rd eigenfunctions are
Courant sharp and have stable nodal domain configurations as the neck-width decreases.
}
\label{fig:num3}
\end{figure}

The third experiment is designed to test the role of the tree hypothesis.
We join the same two subdomains by two necks. The associated graph now has
two vertices joined by two edges, and hence has cycle rank one. The necks
have length $1$ and width $.1$, and intersect the rectangles at heights
$.5$ from the lower edges and $.4$ from the upper edges. We again use a
finite element mesh size of $.005$, and refinement does not change the
nodal counts.

\begin{figure}[!htbp]
\includegraphics[scale=0.4]{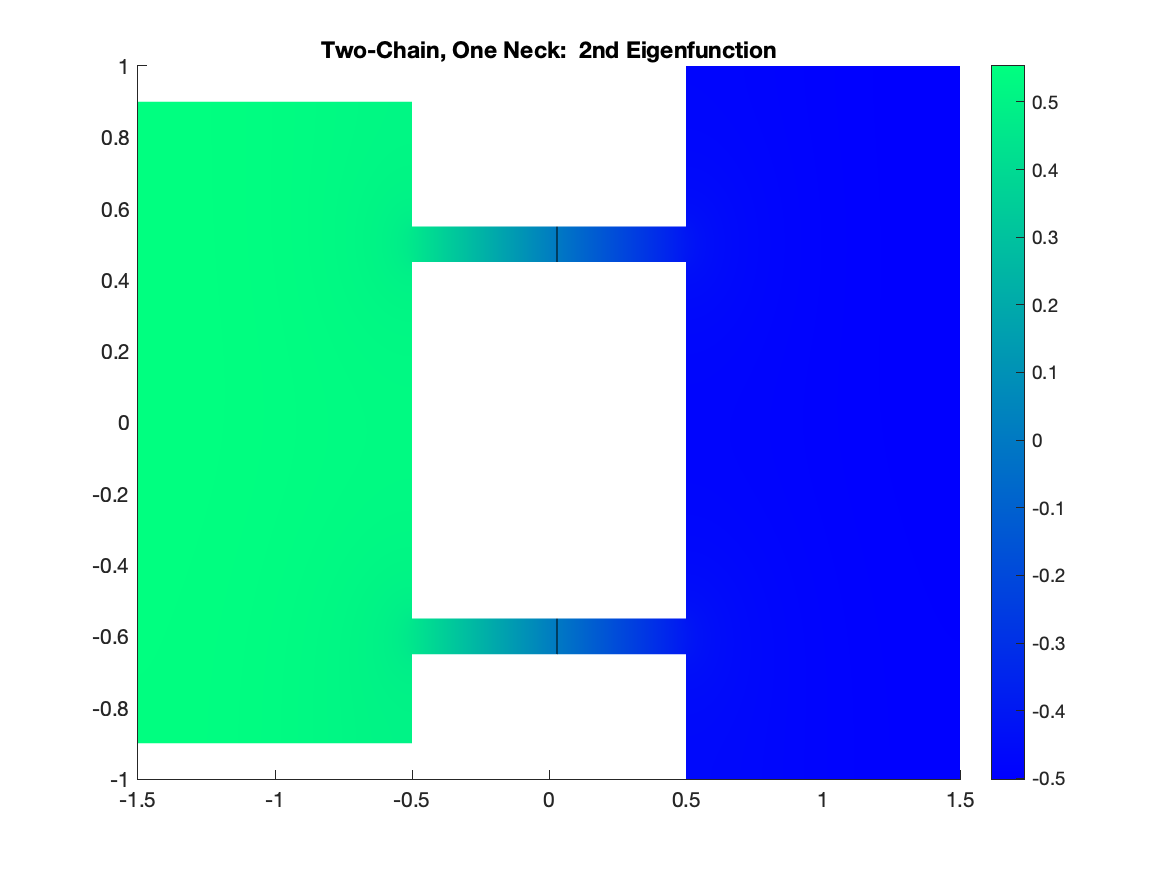}\hspace{-.25in}
\includegraphics[scale=0.4]{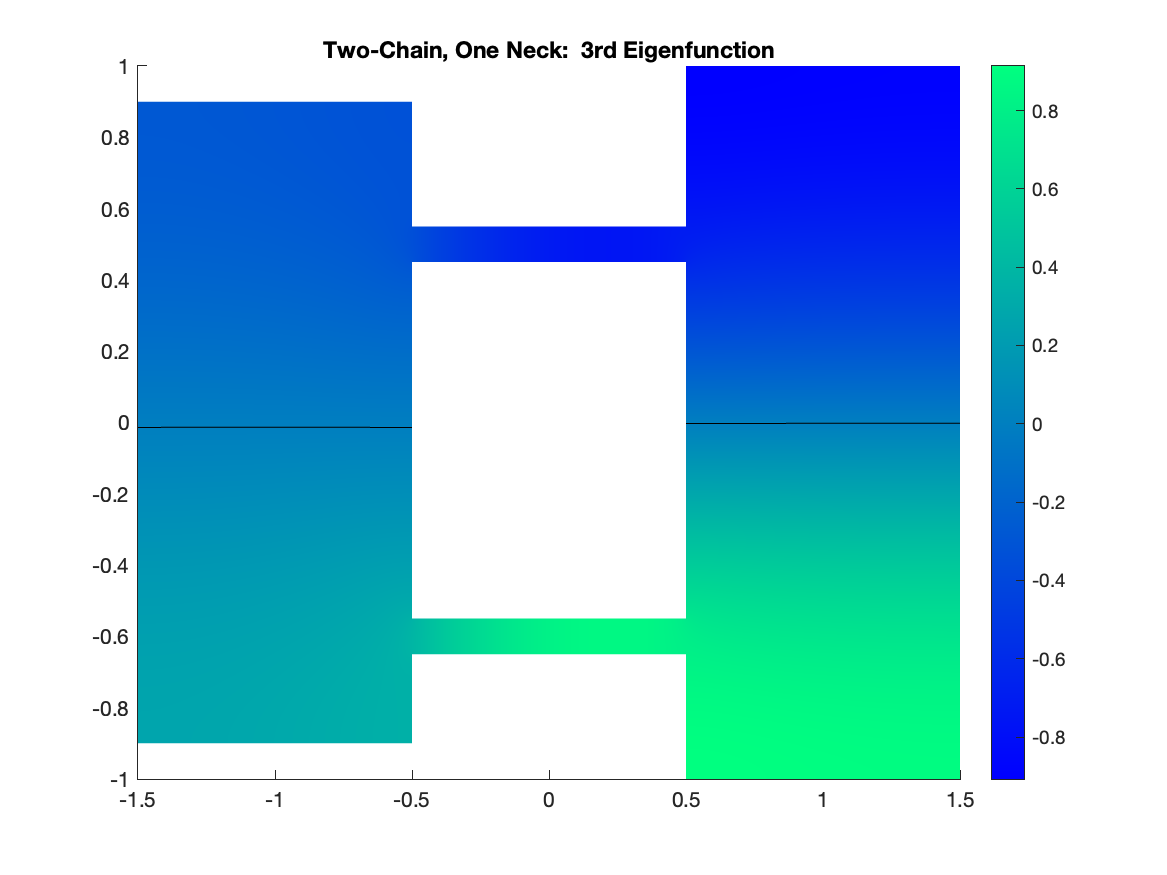} \\
\includegraphics[scale=0.4]{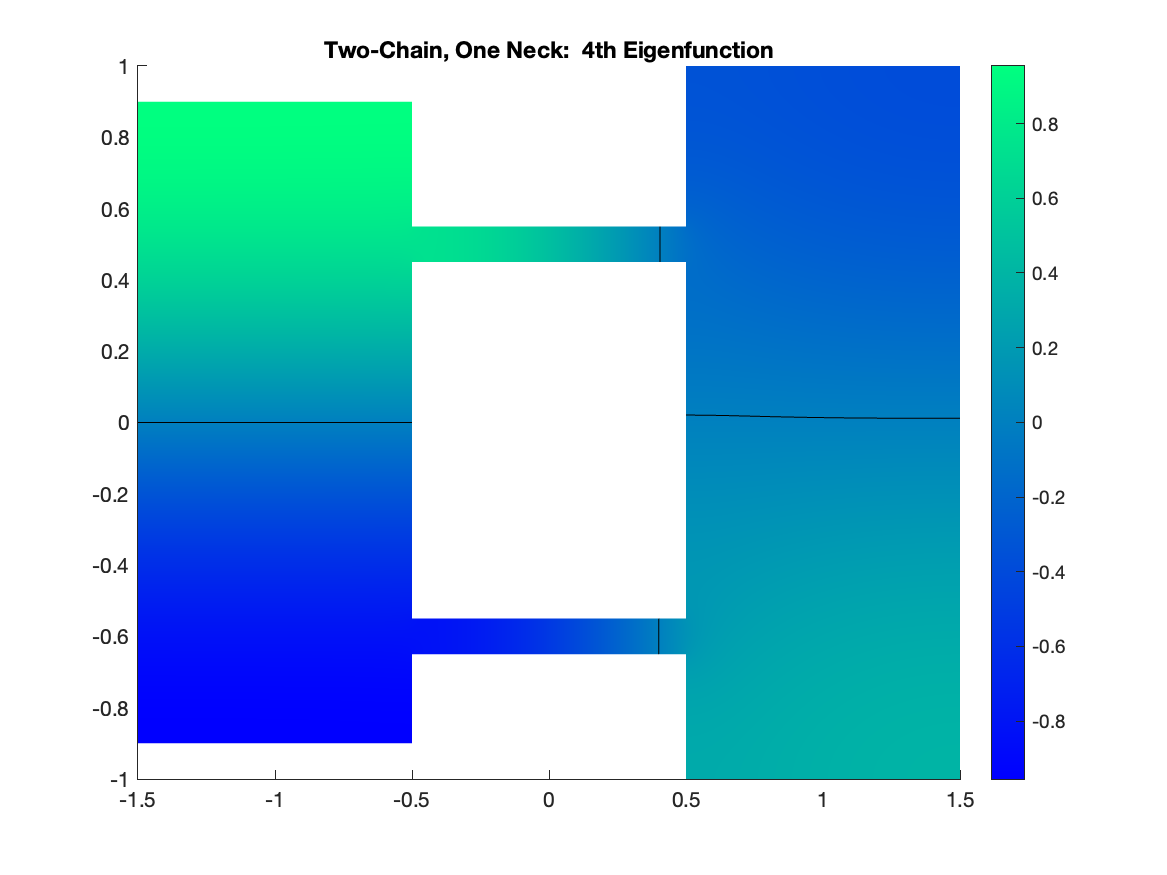}\hspace{-.25in}
\includegraphics[scale=0.4]{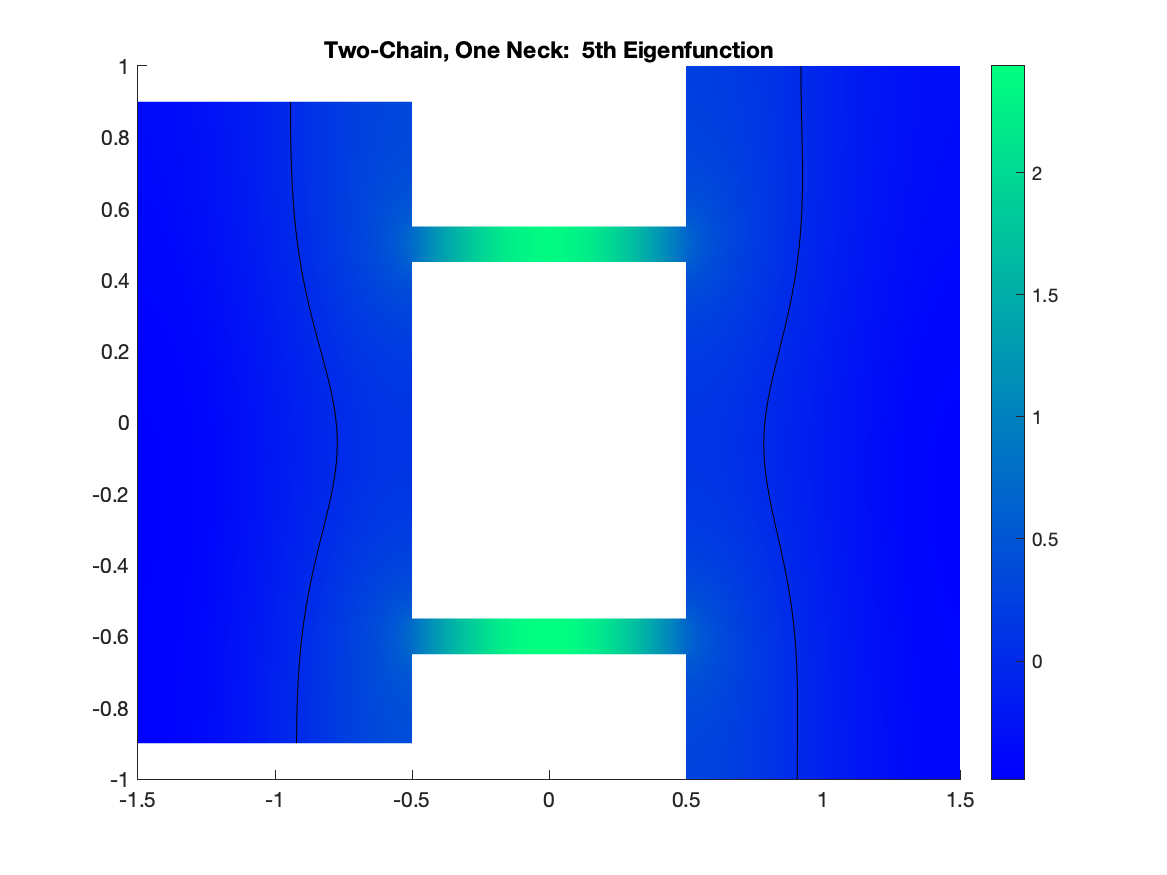}
\caption{The $2$nd (top left), $3$rd (top right), $4$th (bottom left), and
$5$th (bottom right) Neumann eigenfunctions of the $2$-chain with two
necks. The associated graph has one cycle. The $2$nd and $4$th
eigenfunctions are Courant sharp, while the $3$rd and $5$th are not.}
\label{fig:num2}
\end{figure}

In this case, the third and fifth eigenfunctions are not Courant sharp,
while the second and fourth are; see Figure \ref{fig:num2}. In particular,
the failure of Courant sharpness already at the third eigenfunction is the
feature relevant to Conjecture \ref{conj:1}. For a chain with $M=2$
subdomains, the conjecture predicts Courant sharpness through the
$(M+1)$-st eigenfunction when the associated graph is a tree. The two-neck
example shows that this conclusion can fail when the associated graph has a
cycle.

These experiments illustrate the distinction underlying
Conjecture \ref{conj:1}. In the one-neck case, the associated graph is a
tree and the first modes exhibit the expected low-energy nodal rigidity.
In the two-neck case, where the associated graph has a cycle, Courant
sharpness can fail already at the third eigenfunction. The computations
therefore suggest that low-energy nodal counts retain information about how
the subdomains of the chain are connected.

\clearpage

\bibliographystyle{plain}
\bibliography{nodal}

\end{document}